\pdfoutput=1
\RequirePackage{ifpdf}
\ifpdf 
\documentclass[pdftex]{sigma}
\else
\documentclass{sigma}
\fi

\numberwithin{equation}{section}

\newtheorem{Theorem}{Theorem}[section]
\newtheorem*{Theorem*}{Theorem}
\newtheorem{Corollary}[Theorem]{Corollary}
\newtheorem{Lemma}[Theorem]{Lemma}
\newtheorem{Proposition}[Theorem]{Proposition}
 { \theoremstyle{definition}
\newtheorem{Definition}[Theorem]{Definition}

\newtheorem{Example}[Theorem]{Example}
\newtheorem{Remark}[Theorem]{Remark}
\newtheorem{notate}[Theorem]{Notation}
}
\newtheorem*{proc contact}{Procedure Contact (Procedure A)}

\usepackage{amsxtra,amscd}

\usepackage{tikz}
\usepackage{tikz-cd}
\usepackage{comment}

\usepackage{mathtools}

\usepackage{relsize}
\usepackage[mathscr]{eucal}
\usepackage{bm}

\newcommand{\ann}{\operatorname{ann} }
\newcommand{\mcal}[1]{\mathcal{#1}}
\newcommand{\wh}[1]{\widehat{#1}}
\newcommand{\vel}[1]{\operatorname{vel}(\mathcal{#1})}
\newcommand{\dec}[1]{\operatorname{decel}(\mathcal{#1})}
\newcommand{\der}[2]{\mathcal{#1}^{(#2)}}
\newcommand{\coder}[2]{{#1}^{(#2)}}

\newcommand{\Char}[2]{\operatorname{Char} {\mathcal{#1}}^{(#2)}}

\newcommand{\inCharOne}[1]{\operatorname{Char} {\mathcal{#1}}^{(1)}_{0}}

\newcommand{\CommentBlock}[1]{\ignorespaces}

\newcommand{\ceil}[1]{ \lceil {#1} \rceil}
\newcommand{\floor}[1]{ \lfloor {#1} \rfloor }

\def\mb #1{\mathbf{#1}}
\def\wh{\widehat}
\def\bsy{\boldsymbol}

\def\wt{\widetilde}
\def\ch #1{\operatorname{Char}#1}

\def\sfl{static feedback linearizable}
\def\dfl{dynamic feedback linearizable}
\def\stf{static feedback}
\def\df{dynamic feedback}

\def\cfl{cascade feedback linearizable}

\def\arccot{\operatorname{arccot}}

\def\CV{{\mathcal{V}}}
\def\CH{{\mathcal{H}}}

\def\CB{{\mathcal{B}}}

\def\CH{{\mathcal{H}}}

\def\CT{{\mathcal{T}}}

\def\mcal #1{\mathcal{#1}}

\def\wt{\widetilde}

\def\a{\alpha}
\def\b{\beta}
\def\D{\Delta}
\def\e{\varepsilon}
\def\g{\gamma}
\def\G{\Gamma}

\def\l{\lambda}
\def\L{\Lambda}
\def\k{\kappa}

\def\vf{\varphi}

\def\o{\omega}

\def\r{\rho}

\def\t{\tau}
\def\u{\nu}
\def\th{\theta}

\def\B#1{\mathbb{#1}}

\def\P #1{\partial_{#1}}

\def\B#1{\mathbb{#1}}

\def\P #1{\partial_{#1}}

 	\definecolor{apricot}{rgb}{0.98, 0.81, 0.69}
	\definecolor{aqua}{rgb}{0.0, 1.0, 1.0}
	\definecolor{aquamarine}{rgb}{0.5, 1.0, 0.83}	
	\definecolor{blond}{rgb}{0.98, 0.94, 0.75}	
	\definecolor{blizzardblue}{rgb}{0.67, 0.9, 0.93}	
	\definecolor{cambridgeblue}{rgb}{0.64, 0.76, 0.68}	
	\definecolor{carolinablue}{rgb}{0.6, 0.73, 0.89}	
	\definecolor{champagne}{rgb}{0.97, 0.91, 0.81}	
	\definecolor{cream}{rgb}{1.0, 0.99, 0.82}	
	\definecolor{cornsilk}{rgb}{1.0, 0.97, 0.86}
	\definecolor{deepchampagne}{rgb}{0.98, 0.84, 0.65}
	\definecolor{britishracinggreen}{rgb}{0.0, 0.26, 0.15}
	\definecolor{arsenic}{rgb}{0.23, 0.27, 0.29}
	\definecolor{darksienna}{rgb}{0.24, 0.08, 0.08}
	\definecolor{black}{rgb}{0.00, 0.00, 0.00}
	\definecolor{lavenderpink}{rgb}{0.98, 0.68, 0.82}
	\definecolor{maroon(html/css)}{rgb}{0.5, 0.0, 0.0}
	\definecolor{onyx}{rgb}{0.06, 0.06, 0.06}
	\definecolor{paleblue}{rgb}{0.69, 0.93, 0.93}
	\definecolor{pearl}{rgb}{0.94, 0.92, 0.84}
	\definecolor{paleblue}{rgb}{0.69, 0.93, 0.93}
	\definecolor{smokyblack}{rgb}{0.06, 0.05, 0.03}
	\definecolor{springgreen}{rgb}{0.0, 1.0, 0.5}
	\definecolor{tan}{rgb}{0.82, 0.71, 0.55}
	\definecolor{sunset}{rgb}{0.98, 0.84, 0.65}
	\definecolor{timberwolf}{rgb}{0.86, 0.84, 0.82}
	\definecolor{vanilla}{rgb}{0.95, 0.9, 0.67}
	\definecolor{whitesmoke}{rgb}{0.96, 0.96, 0.96}
	\definecolor{prussianblue}{rgb}{0.0, 0.19, 0.33}
	\definecolor{peach}{rgb}{1.0, 0.9, 0.71}
	\definecolor{white}{rgb}{1.0,1.0,1.0}
	\definecolor{peterGold}{rgb}{1.0,0.95,0.78}
	\definecolor{peterBlue}{rgb}{0.89,1.0,1.0}
	\definecolor{peterRose}{rgb}{1.0,0.94,0.84}
	\definecolor{peterBuff}{rgb}{1.0,1.0,0.84}

\newcommand{\LieD}{{\EuScript L}}

\makeatletter
\DeclareRobustCommand\widecheck[1]{{\mathpalette\@widecheck{#1}}}
\def\@widecheck#1#2{%
 \setbox\z@\hbox{\m@th$#1#2$}%
 \setbox\tw@\hbox{\m@th$#1%
 \widehat{%
 \vrule\@width\z@\@height\ht\z@
 \vrule\@height\z@\@width\wd\z@}$}%
 \dp\tw@-\ht\z@
 \@tempdima\ht\z@ \advance\@tempdima2\ht\tw@ \divide\@tempdima\thr@@
 \setbox\tw@\hbox{%
 \raise\@tempdima\hbox{\scalebox{1}[-1]{\lower\@tempdima\box
\tw@}}}%
 {\ooalign{\box\tw@ \cr \box\z@}}}
\makeatother

\begin{document}
\allowdisplaybreaks

\renewcommand{\thefootnote}{}

\newcommand{\arXivNumber}{2103.05078}

\renewcommand{\PaperNumber}{058}

\FirstPageHeading

\ShortArticleName{Dynamic Feedback Linearization of Control Systems with Symmetry}

\ArticleName{Dynamic Feedback Linearization of Control Systems\\ with Symmetry\footnote{This paper is a~contribution to the Special Issue on Symmetry, Invariants, and their Applications in honor of Peter J.~Olver. The~full collection is available at \href{https://www.emis.de/journals/SIGMA/Olver.html}{https://www.emis.de/journals/SIGMA/Olver.html}}}

\Author{Jeanne N. CLELLAND~$^{\rm a}$, Taylor J.~KLOTZ~$^{\rm b}$ and Peter J.~VASSILIOU~$^{\rm c}$}

\AuthorNameForHeading{J.N.~Clelland, T.J.~Klotz and P.J.~Vassiliou}

\Address{$^{\rm a)}$~Department of Mathematics, 395 UCB, University of Colorado,\\
\hphantom{$^{\rm a)}$}~Boulder, CO 80309-0395, USA}
\EmailD{\href{mailto:Jeanne.Clelland@colorado.edu}{Jeanne.Clelland@colorado.edu}}

\Address{$^{\rm b)}$~Department of Mathematics, University of Hawaii at Manoa,\\
\hphantom{$^{\rm b)}$}~2565 McCarthy Mall (Keller Hall 401A), Honolulu, Hawaii 96822, USA}
\EmailD{\href{mailto:klotz@hawaii.edu}{klotz@hawaii.edu}}

\Address{$^{\rm c)}$~Mathematical Sciences Institute, Australian National University,\\
\hphantom{$^{\rm c)}$}~Canberra, ACT, 2601 Australia}
\EmailD{\href{mailto:Peter.Vassiliou@anu.edu.au}{Peter.Vassiliou@anu.edu.au}}

\ArticleDates{Received July 29, 2023, in final form May 30, 2024; Published online July 01, 2024}

\Abstract{Control systems of interest are often invariant under Lie groups of transformations. For such control systems, a geometric framework based on Lie symmetry is formulated, and from this a sufficient condition for dynamic feedback linearizability obtained. Additionally, a systematic procedure for obtaining all the smooth, generic system trajectories is shown to follow from the theory. Besides smoothness and the existence of symmetry, no further assumption is made on the local form of a control system, which is therefore permitted to be fully nonlinear and time varying. Likewise, no constraints are imposed on the local form of the dynamic compensator. Particular attention is given to the consideration of geometric (coordinate independent) structures associated to control systems with symmetry. To show how the theory is applied in practice we work through illustrative examples of control systems, including the vertical take-off and landing system, demonstrating the significant role that Lie symmetry plays in dynamic feedback linearization. Besides these, a number of more elementary pedagogical examples are discussed as an aid to reading the paper. The constructions have been automated in the {\sc Maple} package {\sc Differential\-Geo\-metry}.\looseness=1}


\Classification{53A55; 58A17; 58A30; 93C10}

\renewcommand{\thefootnote}{\arabic{footnote}}
\setcounter{footnote}{0}

\section{Introduction}\label{intro}

A very interesting question in geometric control theory is how to determine whether a given control system is dynamic feedback linearizable, and if so how to construct such a linearization.
The notion of dynamic feedback linearization is closely related to the notion of {\it explicit integrability}. A smooth control system
\begin{equation}\label{controlSystem}
\dot{x}=f(t, x, u),\qquad x\in\B R^n,\quad u\in\B R^m
\end{equation}
is said to be {\it explicitly integrable} if the set of all its generic, time-parametrized, smooth trajectories can be locally expressed in the form
\begin{align}
x(t) & ={A}\bigl(t, {z^a}(t), ({z^a})'(t),\dots, ({z}^a)^{(r_a)}(t)\bigr),\nonumber\\
{u}(t) & ={B}\bigl(t, {z^a}(t), ({z^a})'(t),\dots, ({z}^a)^{(s_a)}(t)\bigr),\label{explicitlyIntegrable}
\end{align}
where $a$ ranges over the integers $1\leq a\leq m$ and for some integers $r_a>0$, $s_a>0$, where the~$z^a(t)$ are arbitrary, smooth, vector-valued functions of time $t$, and ${A}$ and ${B}$ are assumed to be smooth functions over an open subset of a finite jet space that have been composed with a~section of this space over its source manifold, $\B R$. Note that the highest-order derivatives~$r_a$,~$s_a$ of the functions~$z^a(t)$ appearing in \eqref{explicitlyIntegrable} need not all be equal; in general they depend on $a$.
The expression \eqref{explicitlyIntegrable} is also referred to as a {\it Monge parametrization} for the system \eqref{controlSystem}; see, e.g., \cite{AP07}.

The appeal of explicit integrability is evident: An explicit Monge pa\-ra\-me\-tri\-za\-tion of the form~\eqref{explicitlyIntegrable} allows the general solution of the system \eqref{controlSystem} to be expressed in terms of $m$ arbitrary functions $z^1(t), \dots, z^m(t)$ and their derivatives, with no integrations required. This is extremely useful for control engineering problems such as motion planning; see, e.g., \cite{Levine}.
And it turns out that explicit integrability is equivalent to dynamic feedback linearizability. We will give a~detailed explanation of this equivalence in Section \ref{terminology-sec}.

Apart from in a dedicated review, it is difficult to do justice to all the advances that have been made in this area, but a non-exhaustive list might include \cite{ArandaMoogPomet, AP, AP07, BC3control, CharletLevineMarino2, Chetverikov, ClellandHu, FLRM1, FLRM2, Guay, Levine, Nic_Res_2017_Flatness, Pomet95, Pomet97, Shadwick90, SPR, SluisTilbury, NMR}. So far as we are aware, references particularly concerned with methods for the explicit construction of dynamic feedback linearizations of various classes of control systems include
\cite{AP07,BC3control, CharletLevineMarino2, Chetverikov, Guay, LNR, Pomet95, R2, Pomet97, R2}.
Despite this progress, however, it is fair to say that the theory of \df\ linearizability is still under development and numerous open questions remain before the class of dynamic feedback linearizable systems can be said to be well understood.

The primary aim of this paper is to set out a new, widely applicable theory of dynamic feedback linearization~-- and hence, of explicit integrability~-- in the presence of symmetry. Our approach is independent of the local form of a given control system (driftless, control affine, fully nonlinear, time-varying, etc), including the numbers of states and inputs. It relies on the existence of a static feedback linearizable quotient, or ``{\it symmetry reduction}'', which turns out to be very powerful.

While the existence of such a quotient is a geometric (i.e., coordinate independent) phenomenon, one might suppose that it is a rare event. Surprisingly, this {\it does not} appear to be the case. We have by now examined dozens of control systems, both dynamic feedback linearizable and not. It is very rare to come across a control system with symmetry that does not possess a~static feedback linearizable quotient control system, and there are usually numerous such quotients for a given control system. By itself, existence of a \sfl\ quotient does not imply dynamic feedback linearizability, but it is an important first step.

Let us briefly describe our approach in broad terms. Suppose that a control system is invariant under a Lie group $G$ acting by {\it control symmetries} (cf.\ Definition~\ref{admissibleSymmetries}). Very often, there exist subgroups $K\subseteq G$ with the property that the {\em quotient} control system
by the action of $K$ is static feedback linearizable \cite{DTVa}. In fact,
checking for subgroups $K$ that lead
to \sfl\ quotients\footnote{For control systems with symmetry, the notion of a {\it \sfl\ quotient} \cite{DTVa} is a refinement and generalization of the notion of {\it partial feedback linearization} studied in \cite{MarinoBoothby85} and \cite{Marino} with the role played by \sfl \ ``subsystem'' introduced in these references being replaced
by a group theoretic quotient control system which is \sfl.}
is an algorithmic procedure which can often be accomplished quickly and efficiently once the infinitesimal generators of the Lie transformation group $G$ are known (see~\cite{DTVa} for further details and examples).

We have found that, in practice, there is usually a plentiful supply of subgroups leading to \sfl\ quotients.
Once such a \sfl\ quotient is constructed, a further step
answers the questions of the {\em existence} and subsequently {\it construction} of a dynamic feedback linearization and its explicit solution. (See Section \ref{dynamicExtensionsSection} for details.)

Let us now give an outline of this paper, highlighting the main results. Section \ref{terminology-sec} sets out our basic terminology and the notion of a regular dynamic feedback linearization satisfying the solution correspondence condition. For such control systems it is proven that explicit integrability and dynamic feedback linearization are equivalent. For completeness, we also comment on the relationship between dynamic feedback linearization and the notion of (differential) flatness. Section 3 begins to set out the technical tools we use to achieve our results.
We introduce the notion of a Goursat bundle and its invariants, an exposition of the main procedure we use for determining \stf\ linearizing transformations, and the notion of control admissible symmetries. The latter is the main tool we need to determine dynamic feedback linearizations. An important ingredient is a means of very quickly identifying \sfl\ quotients of an invariant control system; this is discussed in
Section \ref{relativeGoursatSect}. The fundamental notion of a contact sub-connection is introduced in Section \ref{ConnectionSect}, and in Theorem \ref{ConstructingConnectionThm}, we prove a normal form result for any invariant control system that has a \sfl\ quotient.

Section \ref{contactReduction} gives a brief exposition of cascade feedback linearization and the notion of a partial contact curve reduction that we use in our construction of dynamic feedback linearization of control systems with symmetry.
Section~\ref{CommutativitySection} proves the preliminary results we need for Section~\ref{dynamicExtensionsSection} in which we apply all previously discussed notions to
establish our theory of dynamic feedback linearization. There are two main results: Theorem \ref{prolongationPredictor}, which proves that if a control system is cascade feedback linearizable then its contact sub-connection has a canonical partial prolongation by differentiation that is static feedback linearizable, and Theorem \ref{cfl to dfl}, which gives a~dynamic feedback linearization constructed from an application of Theorem \ref{prolongationPredictor}. As part of this we determine which inputs are to be differentiated and obtain a bound on the number of differentiations that must be performed to obtain a \sfl\ dynamic extension.
In the case of differentially flat invariant control systems, our procedure canonically determines the flat outputs.

Throughout we discuss numerous local coordinate examples in some detail, including pedagogical examples that illustrate the various constructions. We believe the examples are an important part of the paper, and have been included with the aim of making the paper easier to follow.

Finally, we remark that the approach of this paper is not restricted to static feedback equivalence. The theory of Goursat bundles allows one to replace all of the static feedback equivalences by orbital feedback equivalences \cite{FLRM2}.

\section[Background on explicit integrability, dynamic feedback linearizability, and flatness]{Background on explicit integrability, \\ dynamic feedback linearizability, and flatness}
\label{terminology-sec}

Throughout this paper, all assumptions and results are local, and singularities will not be considered. Statements such as ``assume that $X \neq 0$'' should be interpreted as ``assume that $X$ is not identically zero and restrict to the open set where $X \neq 0$''. All functions are assumed to be smooth (i.e., $C^\infty$) unless otherwise indicated.

\begin{Definition}
A smooth control system \eqref{controlSystem} on the manifold $M = \B R\times \B R^n \times \B R^m$ satisfying the requirement that in local coordinates $(t, x, u)$, the map $u\mapsto f(t, x, u)$ has maximal rank for all $(t, x)$ will be called {\it regular}, or {\it nonsingular}.
\end{Definition}

Throughout this paper, we will assume that all control systems under consideration are regular.

\begin{Definition}\label{SFLdefn0}
A control system of the form \eqref{controlSystem} is called {\em static feedback linearizable $($SFL$)$} if there exists a transformation of the form
\begin{equation}\label{SFL-transformation}
x = \th(t, z), \qquad u = \psi(t, z, v), \qquad z \in \B R^n,\quad v\in\B R^m
\end{equation}
such that $(t, z,v) \mapsto (t, \th(t, z), \psi(t, z, v))$ is a local diffeomorphism that transforms the system~\eqref{controlSystem} to a controllable linear system of the form
\[ \dot{z} = Az + Bv, \]
where $A$, $B$ are constant matrices.
\end{Definition}

It was shown by Brunovsk\'y in \cite{Brunovsky} that every controllable linear system is locally static feedback equivalent to a system of the form
\begin{equation}\label{Brunovsky-form-eq}
\dot{z}^a_\ell = z^a_{\ell+1}, \qquad 0 \leq \ell \leq r_a-1, \quad 1 \leq a \leq m
\end{equation}
for some integers $r_1, \dots, r_m \geq 1$. The system \eqref{Brunovsky-form-eq} is said to be in {\em Brunovsk\'y normal form}; it models a system of $m$ smooth functions $z^1(t), \dots, z^m(t)$ with unconstrained dynamics, where each jet coordinate $z^a_{\ell}$ represents the $\ell$th derivative $(z^a)^{\ell}(t)$ of the function $z^a(t)$ and the highest-order coordinates $z^a_{r_a}$ are the control inputs. Thus, an equivalent formulation of Definition \ref{SFLdefn0} is the following:

\begin{Definition}\label{SFLdefn}
A control system \eqref{controlSystem} is {\em static feedback linearizable} (SFL) if
there exists a~transformation of the form \eqref{SFL-transformation} such that $(t, z,v) \mapsto (t, \th(t, z)$, $\psi(t, z, v))$ is a local diffeomorphism that transforms the system \eqref{controlSystem} to a linear system in Brunovsk\'y normal form.
\end{Definition}

\begin{Remark}
In Definitions \ref{SFLdefn0} and \ref{SFLdefn}, we allow the possibility that the transformation \eqref{SFL-transformation} depends nontrivially on the time parameter $t$, even when the underlying system \eqref{controlSystem} is {\it autonomous}, or time-independent. This flexibility will be important for some of our results, and it does not impact well-known results regarding static feedback linearizability, such as those due to Jakubczyk and Respondek \cite{RespondekLin}, van der Schaft \cite{vdS84}, and Gardner and Shadwick \cite{GSalgorithm}. See \cite{PdS08} for a full discussion.
\end{Remark}

\begin{Definition}\label{dynamic-feedback-def}
A {\em dynamic feedback} (or {\em dynamic compensator}) for the system \eqref{controlSystem} is an augmented control system of the form
\begin{align}
\dot{x} & = f(t, x, \beta(t, x, y, w)), \qquad
\dot{y} = g(t, x, y, w), \qquad
u = \beta(t, x, y, w),\label{dynamic-feedback-eq}
\end{align}
where $x \in \B R^n$, $y \in \B R^k$, $w \in \B R^{m'}$ with $m' \geq m$.

The dynamic feedback \eqref{dynamic-feedback-eq} is {\em regular} if $m' = m$ and it satisfies the {\em solution correspondence condition:} For every smooth solution $(x(t), u(t))$ to the original system \eqref{controlSystem}, there exist (not necessarily unique) smooth functions $(y(t), w(t))$ that, together with the given functions~$(x(t), u(t))$, identically satisfy the system \eqref{dynamic-feedback-eq}.

\end{Definition}

\begin{Definition}\label{DFLdefn1}
The system \eqref{controlSystem} is {\em regularly dynamic feedback linearizable} if it admits a regular dynamic feedback \eqref{dynamic-feedback-eq} with the property that the augmented system \eqref{dynamic-feedback-eq} is static feedback linearizable.
\end{Definition}

It is important to note that the dimension $(n+k)$ of the state space for the augmented system is, in principle, unbounded. Most existing results characterizing dynamic feedback linearizable systems rely on assumptions restricting either the state dimension $(n+k)$ or the particular form of the dynamic compensator, leaving open the question of whether a control system that cannot be linearized via a dynamic feedback of a particular state dimension and/or particular form might still be linearizable via a dynamic feedback of some higher state dimension and/or some more general form.

\begin{Definition}\label{EIdefn}
The system \eqref{controlSystem} is {\em explicitly integrable} if it admits a Monge parametrization of the form \eqref{explicitlyIntegrable} that describes the set of all generic trajectories of the control system \eqref{controlSystem} in terms of $m$ arbitrary functions $z^1(t), \dots, z^m(t)$ and their derivatives. We will also refer to the Monge parametrization \eqref{explicitlyIntegrable} as an {\em explicit solution} of the control system \eqref{controlSystem}.
\end{Definition}

For purposes of this paper, we will often (though not exclusively) represent a control system \eqref{controlSystem} as a codistribution $\bsy{\o}$ on the manifold $M = \B R\times \B R^n \times \B R^m$, with local coordinates~${\bigl(t, x^1, \dots, x^n, u^1, \dots, u^m\bigr)}$ by
\[
\bsy{\o} = \operatorname{span}\bigl\{\xi^1, \dots, \xi^n \bigr\},
\]
where
\begin{equation}\label{define-control-1-forms-eq}
\xi^i = {\rm d}x^i - f^i(t, x, u) {\rm d}t , \qquad 1 \leq i \leq n .
\end{equation}
The codistribution $\bsy{\o}$ will also be referred to as a {\em Pfaffian system,} as it represents a set of Pfaffian equations that must be satisfied by solutions to the control system \eqref{controlSystem}.
In this language, the Brunovsk\'y system \eqref{Brunovsky-form-eq} is represented by the codistribution $\bsy{\beta}$ on the manifold
\[ J_{\bsy{\beta}} = \B R\times \B R^{r_1+1} \times \cdots \times \B R^{r_m+1} \]
generated by the 1-forms
\begin{equation}\label{brun pfaff}
 \bsy{\beta} = \operatorname{span}\bigl\{ \eta^1_0, \dots, \eta^1_{r_1-1}, \dots, \eta^m_0, \dots, \eta^m_{r_m-1} \bigr\} ,
\end{equation}
where
\begin{equation*}
\eta^a_\ell = {\rm d}z^a_\ell - z^a_{\ell+1}{\rm d}t, \qquad 0 \leq \ell \leq r_a-1, \qquad 1 \leq a \leq m.
\end{equation*}

The following proposition shows that the notions of ``explicitly integrable'' and ``regularly dynamic feedback linearizable'' are equivalent for regular control systems.

\begin{Proposition}\label{AEI-DFL-prop}
A regular control system of the form \eqref{controlSystem} is explicitly integrable if and only if it is regularly dynamic feedback linearizable.
\end{Proposition}

\begin{proof}
Suppose that the control system \eqref{controlSystem} admits a regular dynamic feedback linearization of the form \eqref{dynamic-feedback-eq}. Then this augmented system is defined on the manifold $\widehat{M} = \B R\times \B R^n \times \B R^k \times \B R^m$ with local coordinates $(t, x, y, w)$, and is represented by the Pfaffian system
\[
\widehat{\bsy{\o}} = \operatorname{span}\bigl\{ \xi^1, \dots, \xi^{n+k} \bigr\} ,
\]
where
\begin{gather}
\xi^i = {\rm d}x^i - f^i(t, x, \beta(t, x, y, w))\, {\rm d}t , \qquad 1 \leq i \leq n ,\nonumber \\
\xi^{n+j} = {\rm d}y^j - g^j(t, x, y, w)\, {\rm d}t, \qquad 1 \leq j \leq k.\label{define-augmented-control-1-forms-eq}
\end{gather}
We have deliberately used the same notation for the 1-forms $\xi^1, \dots, \xi^n$ in equations \eqref{define-control-1-forms-eq} and~\eqref{define-augmented-control-1-forms-eq}, and this should be interpreted as follows: The solution correspondence condition for the dynamic feedback \eqref{dynamic-feedback-eq} implies that the map $\bar{\pi}\colon \widehat{M} \to M$ defined by
\[ (t, x, u) = \bar{\pi}(t, x, y, w) = (t, x, \beta(t, x, y, w)), \]
is a local submersion onto $M$. Consequently, the pullback map \smash{$\bar{\pi}^*\colon \Omega^*(M) \to \Omega^*\bigl(\widehat{M}\bigr)$} is an immersion, and so any differential form $\theta$ on $M$ may naturally be identified with its pullback~$\bar{\pi}^*(\theta)$ on \smash{$\widehat{M}$}. Using this identification, we regard the 1-forms $\xi^1, \dots, \xi^n$ defined by equation \eqref{define-control-1-forms-eq} as identical to those defined by the first equation in \eqref{define-augmented-control-1-forms-eq}, and the
Pfaffian system~${\bsy{\o} = \operatorname{span}\bigl\{ \xi^1, \dots, \xi^n \bigr\}}$ as a subsystem of the system $\widehat{\bsy{\o}} = \operatorname{span}\bigl\{ \xi^1, \dots, \xi^{n+k} \bigr\}$ on \smash{$\widehat{M}$}.

Now, the assumption that \eqref{dynamic-feedback-eq} is static feedback linearizable means that there exists a~diffeomorphism \smash{$\Phi\colon J_{\bsy{\beta}} \to \widehat{M}$} for some Brunovsk\'y normal form $\bsy{\beta}$, with the property that~${\Phi^*\bigl(\widehat{\bsy{\o}}\bigr) = \bsy{\beta}}$. By applying the composition $\pi^* = \Phi^* \circ \bar{\pi}^*$ (where $\pi = \bar{\pi} \circ \Phi$) to the Pfaffian system~$\bsy{\o}$ on~$M$, we may regard $\bsy{\o}$ as a subsystem of the Brunovsk\'y system $\bsy{\beta}$ on $J_{\bsy{\beta}}$. Furthermore, the map~${\pi\colon J_{\bsy{\beta}} \to M}$ is a local submersion that provides an explicit Monge parametrization
\begin{align}
&x^i = A^i\bigl(t, z^1_0, \dots, z^1_{r_1}, \dots, z^m_0, \dots, z^m_{r_m}\bigr), \nonumber\\
&u^a = B^a\bigl(t, z^1_0, \dots, z^1_{r_1}, \dots, z^m_0, \dots, z^m_{r_m}\bigr) .\label{Monge-param-eq-0}
\end{align}
Therefore, the system \eqref{controlSystem} is explicitly integrable.

Observe that equations \eqref{Monge-param-eq-0} imply that the codistribution
\[ \bsy{\o} \oplus \operatorname{span}\{ {\rm d}t \} = \operatorname{span}\bigl\{ {\rm d}t, {\rm d}x^1, \dots, {\rm d}x^n \bigr\} \]
is a subbundle of the codistribution
\[
\bsy{\beta} \oplus \operatorname{span}\{ {\rm d}t \} =\operatorname{span}\bigl\{ {\rm d}t, {\rm d}z^1_0, \dots, {\rm d}z^1_{r_1-1}, \dots, {\rm d}z^m_0, \dots, {\rm d}z^m_{r_m-1} \bigr\}.
\]
It follows that the functions $A^1\!, \dots,\! A^n$ are independent of the highest-order variables $z^1_{r_1}\!, \dots,\! z^m_{r_m}$, while the functions $B^1, \dots, B^m$ may depend nontrivially on these variables. Thus, since $\Phi$ is a~static feedback transformation, we can write
\begin{align}
&x^i = A^i\bigl(t, z^1_0, \dots, z^1_{r_1-1}, \dots, z^m_0, \dots, z^m_{r_m-1}\bigr),\nonumber \\
&u^a = B^a\bigl(t, z^1_0, \dots, z^1_{r_1}, \dots, z^m_0, \dots, z^m_{r_m}\bigr).\label{Monge-param-eq}
\end{align}

Conversely, suppose that the system \eqref{controlSystem} admits a Monge pa\-ra\-me\-tri\-za\-tion of the form~\eqref{Monge-param-eq}. From the definition of explicit integrability, it follows that the formulas \eqref{Monge-param-eq} define a~local submersion $\pi\colon J_{\bsy{\beta}} \to M$. This submersion may be extended to a~local diffeomorphism
\smash{$\Phi\colon J_{\bsy{\beta}} \to \widehat{M} \cong M \times \B R^k$}, where $k = r_1 + \dots +r_m - n$, as follows: Let $y^1, \dots, y^k$ be a relabeling of a~subset of the variables $z^1_0, \dots, z^1_{r_1-1}, \dots, z^m_0, \dots, z^m_{r_m-1}$ with the property that the map
\begin{equation}\label{local-diffeo-eq}
\bigl(t, z^1_0, \dots, z^1_{r_1-1}, \dots, z^m_0, \dots, z^m_{r_m-1}\bigr) \mapsto \bigl(t, x^1, \dots, x^n, y^1, \dots, y^k\bigr),
\end{equation}
with $x^1, \dots, x^n$ defined by the first set of equations in \eqref{Monge-param-eq},
is a local diffeomorphism. Then set
\begin{gather}\label{local-diffeo-eq2}
w^a = z^a_{r_a}, \qquad 1 \leq a \leq m.
\end{gather}
The derivatives $\dot{y}^1, \dots, \dot{y}^k$ are determined by the derivatives of the appropriate coordinates $z^a_j$ in the Brunovsk\'y system $\bsy{\beta}$ on $J_{\bsy{\beta}}$, while the
second set of equations in \eqref{Monge-param-eq}, together with the inverse of the local diffeomorphism defined by \eqref{local-diffeo-eq} and \eqref{local-diffeo-eq2}, determine the functions
\[ u^a = \beta^a(t, x, y, w). \]
Together with the original system \eqref{controlSystem}, these equations constitute a regular dynamic feedback linearization for the system \eqref{controlSystem} of the form \eqref{dynamic-feedback-eq}.
\end{proof}

\begin{Definition}\label{flatness-def}
A regular dynamic feedback linearization \eqref{dynamic-feedback-eq} for the control system \eqref{controlSystem} is called {\em endogenous} if there exists an integer $K \geq 0$ and functions $\varphi^1, \dots, \varphi^m$, each depending on the variables $t$ and $\bigl\{x^i, u^b_0, u^b_1, \dots, u^b_K \mid 1 \leq i \leq n,\, 1 \leq b \leq m \bigr\}$, such that for every solution~${(x(t), u(t))}$ of the system \eqref{controlSystem}, the functions
\begin{equation}\label{flat-outputs-eq}
z^a_0(t) = \varphi^a\bigl(t, x^i(t), u^b(t), \bigl(u^b\bigr)'(t), \dots, \bigl(u^b\bigr)^{(K)}(t)\bigr),\qquad 1\leq a\leq m,
\end{equation}
identically satisfy the Monge parametrization equations \eqref{Monge-param-eq} determined by the dynamic feedback linearization \eqref{dynamic-feedback-eq}.

The control system \eqref{controlSystem} is called {\em flat} if it admits an endogenous regular dynamic feedback linearization. In this case, the functions $\bigl(z^1_0, \dots, z^m_0\bigr)$ determined by equation \eqref{flat-outputs-eq} are called {\em flat outputs} for the system.
\end{Definition}

Flatness may be thought of as follows: Suppose that $\pi\colon J_{\bsy{\beta}} \to M$ is the local submersion corresponding to a regular dynamic feedback linearization. The solution correspondence condition implies that every solution curve
\begin{equation}\label{downstairs-soln-curve}
 t \mapsto (t, x(t), u(t))
\end{equation}
of the system \eqref{controlSystem} in $M$ lifts to at least one contact curve
\begin{equation}\label{upstairs-soln-curve}
 t \mapsto \bigl(t, z^a(t), (z^a)'(t), \dots, (z^a)^{(r_a)}(t)\bigr)
\end{equation}
in $J_{\bsy{\beta}}$. The dynamic feedback is endogenous~-- which implies that the system is flat~-- if and only if every solution curve \eqref{downstairs-soln-curve} can be lifted to a {\em unique} contact curve \eqref{upstairs-soln-curve}.

\section{Geometry of control systems }\label{cascadeIntegrationSect.}

In this section, we review a geometric formulation of control systems expressed in terms of differential geometry on {finite} smooth manifolds. The exposition emphasizes those aspects relevant to the applications that follow. More details can be found in \cite{DTVa,KlotzThesis,VassiliouGoursat, VassiliouGoursatEfficient, VassiliouCascade1}.

\begin{Definition}
A {\it control system} is a parametrized family of ordinary differential equations
\[
\dot{x}=f(t,x,u),\qquad x\in\B R^n,\quad u\in\B R^m
\]
in which the vector $x$ is comprised of the {\it state variables} taking values in some open set
$\mathbf{X}\subseteq \B R^n$ and the vector $u$ is comprised of the {\it inputs} or {\it controls} taking values in some open set $\mathbf{U}\subseteq \B R^m$. Time $t$ takes values in a connected subset of the real line.

Throughout, we very often invoke the {\it Pfaffian system} representation of a control system as the vanishing
of differential 1-forms
\[
\bsy{\o}=\operatorname{span}\bigl\{{\rm d}x^1-f^1(t,x,u){\rm d}t,\, {\rm d}x^2-f^2(t,x,u) {\rm d}t,\,\dots,\,
{\rm d}x^n-f^n(t,x,u) {\rm d}t\bigr\}
\]
defining a sub-bundle of the cotangent bundle $\bsy{\o}\subset T^*(\B R\times\mathbf{X}\times \mathbf{U})$, and we exploit the geometric properties of $\bsy{\o}$ under local changes of variable. By the same token, we often express our control systems dually as a sub-bundle of the tangent bundle
$\ker\bsy{\o}=\CV\subset T(\B R\times\mathbf{X}\times \mathbf{U})$
\[
\CV=\operatorname{span}\left\{\P t+\sum_{i=1}^n f^i(t,x,u)\P {x^i}, \P {u^1}, \P {u^2}, \dots, \P {u^m}\right\},
\]
and frequently switch between the two representations as the need arises. We often refer to
$\bsy{\o}$ and $\CV$ {\em themselves} as {control systems}.
\end{Definition}
\begin{Definition}\label{SFT def}
A diffeomorphism $\varphi\colon M\to N$ of the form
\[
\varphi\colon\ (t,x,u) \mapsto (t,\theta(t,x),\psi(t,x,u))
\]
is called a \textit{static feedback transformation} (SFT). Two control systems $(M,\boldsymbol{\omega})$ and $(N,\boldsymbol{\eta})$ are called \textit{static feedback equivalent} (SFE) if $\varphi^*\boldsymbol{\eta}=\boldsymbol{\omega}$ for some SFT $\varphi\colon M\to N$.
\end{Definition}
Static feedback linearizable control systems represent a special case of static feedback equivalence, namely, those that are SFE to the Brunovsky normal forms.

Let $M:=\B R\times\mathbf{X}\times\mathbf{U}$. When we wish to draw attention to the state space or control space factors of a manifold $M$ carrying a control system, we write $\mathbf{X}(M)$ or $\mathbf{U}(M)$, respectively.
\begin{Definition}\label{EI solSignature}
Suppose that $\bsy{\o}$ is a \dfl\ control system on a manifold~$M$ with explicit solution $s\colon\mathbb{R}\to M$. The number of arbitrary functions in an explicit solution is equal to the number of inputs associated to $\bsy{\o}$. To an explicit solution $s$ we can associate the notion of a {\it signature},
\[
\kappa=\langle \r_1, \r_2, \dots, \r_k\rangle,
\]
where $\r_j$ is the number of arbitrary functions occurring to highest order $j$ in the explicit solution~$s$, where $1\leq j\leq k$.
\end{Definition}

We now present some definitions and notation from the geometry of distributions that are used in this paper.
Let us denote by $\CV^{(j)}$ the $j^{\text{th}}$ derived bundle of $\CV:=\CV^{(0)}$,
defined recursively~by
\[\CV^{(j)}=\CV^{(j-1)}+\big[\CV^{(j-1)},\CV^{(j-1)}\big],\qquad j\geq 1.\]
The sequence
\[
\CV\subset\CV^{(1)}\subset\cdots\subset\CV^{(k)}\subseteq TM
\]
is the {\it derived flag} of $(M, \CV)$ and the integer $k$ is its {\it derived length}. This
is the smallest integer~$k$ such that $\CV^{(k)}=\CV^{(k+1)}$. Throughout we always assume that
$\CV^{(k)}=TM$.
\begin{Definition}\label{decel}
Let $\CV\subset TM$ be a sub-bundle of derived length $k>1$.
The {\it velocity} of $\CV$ is the ordered list of $k$ integers
\[
\text{\rm vel}(\CV)=\langle\D_1,\D_2,\dots,\D_k\rangle,\qquad \text{where}\quad \D_j=\dim\bigl(\CV^{(j)}\bigr)-\dim\bigl(\CV^{(j-1)}\bigr),\qquad 1\leq j\leq k.
\]
The {\it deceleration} of $\CV$ is the ordered list of $k$ integers
\[
\text{\rm decel}(\CV)=\bigl\langle -\D^2_2,-\D^2_3,\dots,-\D^2_k, \D_k\bigr\rangle,\qquad \text{where}\quad \D^2_j=\D_j-\D_{j-1}.
\]

\end{Definition}
Denote by $\ch\CV^{(j)}$ the
{\it Cauchy bundle} of $\CV^{(j)}$,
\[
\ch\CV^{(j)}=\operatorname{span}\bigl\{X\in\CV^{(j)}| [X, \CV^{(j)}]\subset\CV^{(j)}\bigr\},\qquad j\geq 0.
\]
We assume that for all $j\geq 0$, $\CV^{(j)}$ and $\ch\CV^{(j)}$ have constant rank; we refer to such sub-bundles as {\it totally regular}. In this totally regular case
$\ch\CV^{(j)}$ can be shown to be integrable for each~$j$. Define the {\it intersection bundles}
\begin{equation}\label{intersectionBundleDefn}
\ch\CV^{(i)}_{i-1}:=\CV^{(i-1)}\cap\ch\CV^{(i)},\qquad1\leq i\leq k-1.
\end{equation}
Unlike the Cauchy bundles, the intersection bundles $\ch\CV^{(i)}_{i-1}$ are not guaranteed to be integrable; however, this will arise as a condition in the definition of Goursat bundle (cf.\ Definition~\ref{Goursat}).
\begin{notate}\label{Notation3.1}
We will sometimes denote the codistribution version of Cauchy and intersection bundles respectively by
\begin{align*}
\Xi^{(i)}&=\operatorname{ann} \ch\CV^{(i)},\qquad
\Xi^{(i)}_{i-1}=\operatorname{ann}\ch\CV^{(i)}_{i-1}.
\end{align*}
\end{notate}
\begin{Definition}
Let
\begin{gather*}
m_i=\dim\CV^{(i)},\\
\chi^i=\dim\ch\CV^{(i)},\\
\chi^j_{j-1}=\dim\ch\CV^{(j)}_{j-1},\qquad 0\leq i\leq k,\quad 1\leq j\leq k-1,
\end{gather*}
called the {\it type numbers} of $(M, \CV)$. The list of lists of type numbers
\begin{equation}\label{refinedDerivedTypeDefn}
\mathfrak{d}_r(\CV)=\big[ \big[m_0,\chi^0\big],\big[m_1,\chi^1_0,\chi^1\big],\big[m_2,\chi^2_1,\chi^2\big],\dots,\big[m_{k-1}, \chi^{k-1}_{k-2},\chi^{k-1}\big],
\big[m_k,\chi^k\big] \big]
\end{equation}
is called the {\it refined derived type} of $(M, \CV)$.
\end{Definition}

\subsection{Brunovsk\'y normal forms}\label{brun-form-subsect}
In this section, we will give an exposition of the partial prolongations of jet spaces and the Brunovsk\'y normal form.
\begin{Definition}
Let $J^k(\mathbb{R},\mathbb{R}^m)$ be the standard jet space of order $k$ and let $\bsy{\b}^k_{m}$ be the standard contact system on $J^k(\mathbb{R},\mathbb{R}^m)$. We will drop the subscript $m$ in the notation $\bsy{\b}^k_m$ and use the shorthand notation $J^k$ when the integer $m$ is known by context.
\end{Definition}

We now introduce another notion of signature beside that of an explicit solution (see Definition~\ref{EI solSignature}); namely the signature of a {\em control system itself}. It is proven in \cite{VassiliouGoursat} that the deceleration,
$\text{decel}(\CV)$ (Definition \ref{decel}), is a diffeomorphism invariant that uniquely identifies the
Brunovsk\'y normal form of any linearizable control system $\CV$ and we therefore call $\text{decel}(\CV)$ the {\em signature of} $\CV$. As in the case of the signature of an explicit solution, the signature of a control system consists of $k$ non-negative integers
$\text{decel}(\CV)=\langle \r_1, \r_2,\dots, \r_k\rangle$, where $k$ is the derived length of $\CV$. While these are numerical invariants for any control system, if $\CV$ is diffeomorphic to a~Brunovsk\'y normal form
then $\r_j$ in $\text{decel}(\CV)$ is the number of sequences of differential forms of order $j$ in the
Brunovsk\'y normal form of
$\CV$. The signature of a control system is widely used in this paper and it is also convenient as a means of classifying Brunovsk\'y normal forms. Only when a control system is feedback linearizable are the two notions of signature equal.

\begin{Definition}\label{partial prolongation}
A \textit{partial prolongation} of the Pfaffian system $\bigl(J^1(\mathbb{R},\mathbb{R}^m),\bsy{\b}^1\bigr)$ is a Brunovsk\'y form defined in the sense of equation \eqref{brun pfaff}, i.e., the Pfaffian system associated to the Brunovsk\'y normal form of mixed orders. We use $(J^\kappa(\mathbb{R},\mathbb{R}^m),\bsy{\b}^\kappa)$ to refer to the partial prolongation of {signature} $\kappa=\langle \r_1,\dots,\r_k\rangle$ where $k$ is the derived length of $\bsy{\b}^\kappa$.
\end{Definition}

\begin{Remark}
In the case of a linear control system $\dot{x}=Ax+Bu$, the signature of its distribution representation
$\CV$ agrees precisely with the
collection of {\it Kronecker indices} of the pair of matrices $(A,B)$. However, the definition of signature is far more versatile since it can be found without putting the control system into the above linear form required to compute the Kronecker indices. Furthermore, it is defined for general control systems and by its definition is manifestly a diffeomorphism invariant. The two notions of signature play a central role in this paper.
\end{Remark}

It is helpful for us to arrange our Brunovsk\'y normal forms according to their signature. In particular, we will think of the new partially prolonged jet space $J^\kappa(\mathbb{R},\mathbb{R}^m)$ as being constructed from jet spaces $J^i(\mathbb{R},\mathbb{R}^{\r_i})$ of fixed order. However, we cannot use a strict product of jet spaces. We must identify the independent variables (i.e., the source) of each jet space together in a~product like so,
\begin{align}
J^\kappa(\mathbb{R},\mathbb{R}^m)&:=\left(\prod_{i\in I}J^i(\mathbb{R},\mathbb{R}^{\rho_i})\right)\Big/{\sim},\qquad
\bsy{\b}^\kappa_m:=\bigoplus_{i\in I}\bsy{\b}^i_{\rho_i},\label{jkappa}
\end{align}
with
\[
I=\{1\leq a\leq k \mid\rho_a\neq0 \}
\]
 and each $\bsy{\b}^i_{\r_i}$ is the Brunovk\'y form on jet space $J^i(\mathbb{R},\mathbb{R}^{\r_i})$. The equivalence relation `$\sim$' in \eqref{jkappa} is defined by
\begin{equation*}
\pi_i\left(J^i(\mathbb{R},\mathbb{R}^{\rho_i})\right)=\pi_j\left(J^j(\mathbb{R},\mathbb{R}^{\rho_j})\right),
\end{equation*}
for all $1\leq i$, $j\leq k$, where $\pi_i$, $\pi_j$ are source projection maps (i.e., projection on to the $t$-coordinate on $\mathbb{R}$).

\begin{Example}
A Brunovsk\'y normal form of signature $\kappa=\langle1,2,0,0,1\rangle$ on
\begin{equation*}
J^\kappa\bigl(\mathbb{R},\mathbb{R}^5\bigr)=\bigl(J^1(\mathbb{R},\mathbb{R})\times J^2\bigl(\mathbb{R},\mathbb{R}^2\bigr)\times J^5(\mathbb{R},\mathbb{R})\bigr)/{\sim}
\end{equation*}
is generated by the 1-forms
\begin{alignat*}{5}
& && && && \theta^4_4={\rm d}z^4_4-z^4_5\,{\rm d}t,&\\
& && && && \theta^4_3={\rm d}z^4_3-z^4_4 {\rm d}t, &\\
& && && && \theta^4_2={\rm d}z^4_2-z^4_3 {\rm d}t, &\\
& && \theta^2_1={\rm d}z^2_1-z^2_2{\rm d}t,\qquad && \theta^3_1={\rm d}z^3_1-z^3_2 {\rm d}t,\qquad&& \theta^4_1={\rm d}z^4_1-z^4_2 {\rm d}t,&\\
&\theta^1_0={\rm d}z^1_0-z^1_1 {\rm d}t,\qquad&& \theta^2_0={\rm d}z^2_0-z^2_1 {\rm d}t,\qquad&& \theta^3_0={\rm d}z^3_0-z^3_1 {\rm d}t,\qquad&&\theta^4_0={\rm d}z^4_0-z^4_1 {\rm d}t.
\end{alignat*}
In this example, one can say that $J^\kappa$ has one variable of order 1, two of order 2, zero of orders~3 and 4, and one of order 5. So the signature $\kappa$ represents a list of the local coordinates on $J^\kappa$ categorized by order.
\end{Example}
\begin{notate}\label{notation3.2}
For ease of notation, we will often use $\bsy{z}^\k$ to represent all the jet coordinates save $t$ on $J^\k$ space, $\bsy{z}^{\ceil{\k}}$ to refer to all highest order jet coordinates for $J^\k$, and $\bsy{z}^{\floor{\k}}$ to denote all jet coordinates of $J^\k$ that have order strictly smaller than maximal orders given by $\k$.
\end{notate}

\begin{notate}\label{notation3.3}
We will need to use the total derivative operator on each $J^\kappa$. Let $\mathbf{D}_{t,\r_i}$ be the usual total derivative on $J^i(\mathbb{R},\mathbb{R}^{\r_i})$. Then by the definition of $J^\kappa$ it easy to see that the total derivative on $J^\kappa$ is
\begin{equation}\label{totalDiffOp}
\mathbf{D}_t=\partial_t+\sum_{i\in I}(\mathbf{D}_{t,\r_i}-\partial_t).
\end{equation}
\end{notate}
The following proposition characterizes the refined derived type \eqref{refinedDerivedTypeDefn} of the
Brunovsk\'y normal forms.

\begin{Proposition}[\cite{VassiliouGoursatEfficient}]\label{refined derived type numbers}
Let $\mcal{B}_\k\subset TJ^\k$ be the distribution that annihilates the $1$-forms in a~Brunovsk\'y normal form $\bsy{\b}^\kappa$ with signature $\kappa=\langle \rho_1,\dots,\rho_k\rangle$. Then the entries in the refined derived type
\begin{equation*}
\mathfrak{d}_r(\mcal{B}_\k)=\big[\big[m_0,\chi^0\big],\big[m_1,\chi^1_0,\chi^1\big],\dots,\big[m_{k-1},\chi^{k-1}_{k-2},\chi^{k-1}\big],\big[m_k,\chi^k\big]\big]
\end{equation*}
satisfy the following relations:
\begin{align}
&\kappa=\operatorname{decel}(\mcal{B}_\k),\qquad \Delta_i=\sum_{l=i}^{k}\rho_\ell,\nonumber\\
& m_0=1+m,\qquad m_j=m_0+\sum_{l=1}^j\Delta_\ell,\qquad 1\leq j\leq k,\nonumber\\
&\chi^j=2m_j-m_{j+1}-1,\qquad 0\leq j\leq k-1,\nonumber\\
& \chi^i_{i-1}=m_{i-1}-1,\qquad 1\leq i\leq k-1,\label{linearTypeConstraints}
\end{align}
where $\Delta_j$ is given in Definition {\rm\ref{decel}}.
\end{Proposition}

\subsection{Goursat bundles}\label{Goursat-subsect}
Here we provide an brief exposition of the theory of Goursat bundles, \cite{VassiliouGoursat,VassiliouGoursatEfficient} used in this paper and a~discussion of the relevance of this topic to the present study. Certainly Brunovsk\'y normal forms
$\mcal{B}_\k\subset TJ^\k$ are the local normal forms of Goursat bundles. But in the first instance the theory of Goursat bundles handles the case in which a~distribution is equivalent to some Brunovsk\'y normal form via a {\em general diffeomorphism} of the ambient manifolds. That is, if a~distribution $\CV\subset TM$ determines a Goursat bundle, there will exist an equivalence of $\CV$ to a~Brunovsk\'y normal form but not necessarily by a SFT. In that case integral curves may \textit{not have} a parametrization by the time variable $t$, but by some other variable. There is a simple check for when a Goursat bundle $\CV$ has $t$ as a parameter for integral curves (see Theorem~\ref{Goursat SFL}), and if satisfied then $\CV$ is indeed SFE to a Brunovsk\'y normal form. A particular application of the theory in this paper concerns the intersection bundle
\smash{$\ch\CV^{(1)}_0$}, which is integrable and can be used to geometrically characterize the control variables. This fact is exploited in the construction of a dynamic extension in the proof of Theorem \ref{cfl to dfl}.
The theory is also used in the proof of Theorem \ref{prolongationPredictor} and elsewhere.
In particular, the procedure {\it contact}
\cite{VassiliouGoursatEfficient} for efficiently producing normal form coordinates for a Goursat bundle, has been adapted and automated to determine linearizing maps for SFL systems. It has a structure similar to the well-known GS algorithm \cite{GSalgorithm} but expressed in the language of tangent distributions rather than Pfaffian systems. A detailed description of {\it contact} together with application examples appear later in this paper.

A Goursat bundle is described as follows.
\begin{Definition}[\cite{VassiliouGoursat}]\label{Goursat}
A totally regular subbundle $\mathcal{V}\subset TM$ of derived length $k$ with $\Delta_k=1$ will be called a {\bf\it Goursat bundle $($of signature $\kappa)$} if:
\begin{enumerate}\itemsep=0pt
\item[(1)] the subbundle $\mathcal{V}$ has the refined derived type of a partial prolongation of $J^1(\mathbb{R},\mathbb{R}^m)$ (as characterized in Proposition \ref{refined derived type numbers}) whose signature
$\kappa=\operatorname{decel}(\mathcal{V})$;
\item[(2)] each intersection bundle \smash{$\Char{V}{i}_{i-1}:= \mathcal{V}^{(i-1)}\cap\Char{V}{i}$} is an integrable subbundle, the rank of which agrees with the corresponding rank of \smash{$\operatorname{Char}(\mcal{B}_\k)^{(i)}_{i-1}$}. That is, intersection bundle ranks satisfy equations \eqref{linearTypeConstraints}.
\end{enumerate}
\end{Definition}
In the case $\Delta_k>1$, the full theory of Goursat bundles in \cite{VassiliouGoursat,VassiliouGoursatEfficient} requires one to construct an additional bundle, which we can omit for the purposes of the present paper.

Theorem \ref{Generalized Goursat Normal Form} below asserts that Goursat bundles are locally diffeomorphic to the Brunovsk\'y normal forms at generic points; and conversely, every Brunovsk\'y normal form is a Goursat bundle. However, this theorem has nothing to say about the singularities of the related {\em Goursat structures} which have been the subject of recent work; see, for example, \cite{MontgomeryGoursatFlags2001} and citations therein. As is the case for the classical Goursat normal form, the generalized Goursat normal form is concerned with generic local behaviour, in terms of which it geometrically characterizes the partial prolongations of the contact system on $J^1(\B R, \B R^m)$ exclusively in terms of their {\it derived type} \cite{BryantThesis}.\looseness=-1
\begin{Theorem}[generalized Goursat normal form, \cite{VassiliouGoursat}]\label{Generalized Goursat Normal Form}
 Let $\mathcal{V}\subset TM$ be a Goursat bundle on a manifold $M$, with derived length $k>1$, and signature $\kappa=\operatorname{decel}(\mathcal{V})$. Then there is an open dense subset ${\rm U}\subset M$ such that the restriction of $\mathcal{V}$ to ${\rm U}$ is locally equivalent to $\mcal{B}_\k$ via a local diffeomorphism of $M$. Conversely, any partial prolongation of $\mcal{B}_{\langle m \rangle}$ is a Goursat bundle.
\end{Theorem}
The paper \cite{VassiliouGoursat} establishes the local normal form for Goursat bundles constructively. However, in \cite{VassiliouGoursatEfficient}, the construction of local coordinates is streamlined into a nearly algorithmic procedure. We'll next outline this procedure, often referred to as procedure {\it contact}, and apply it to an example in detail.

\begin{Definition}[\cite{VassiliouGoursatEfficient}]
Let $\mathcal{V}$ be a Goursat bundle of derived length $k$ with $\Delta_k=1$, $\tau$~a~first integral of $\Char{V}{k-1}$, and $Z$ any section of $\mathcal{V}$ such that $Z\tau=1$. Then the \textit{fundamental bundle}~${\Pi^{k}\subset\mathcal{V}^{(k-1)}}$ is defined inductively as
\begin{equation*}
\Pi^{\ell+1}=\Pi^\ell+\big[\Pi^\ell,Z\big], \qquad\Pi^0=\inCharOne{V}, \qquad 0\leq \ell\leq k-1.
\end{equation*}
\end{Definition}

The proof of \cite[Theorem 4.2]{VassiliouGoursat} shows that $\Pi^{k}$ is integrable and has corank 2 in $TM$ while in~\cite{VassiliouGoursatEfficient} it is proven that in any Goursat bundle, \smash{$\ch\CV^{(i)}_{i-1}$} and
$\Pi^k$ have the form
\begin{align}
&\Pi^{(k)}=\operatorname{span}\bigl\{\Pi^0,\operatorname{ad}(Z)\Pi^0,\dots,\operatorname{ad}^{k-1}(Z)\Pi^0 \bigr\},\nonumber\\
&\ch\CV^{(i)}_{i-1}=\operatorname{span}\bigl\{C_0,\operatorname{ad}(Z)C_0,\dots,\operatorname{ad}^{i-1}(Z)C_0 \bigr\},\qquad
C_0=\Pi^0,
\qquad 1\leq i \leq k-1,\label{ad char}
\end{align}
once $Z$ and $\tau$ are known.

We shall be making use of the forms \eqref{ad char} for
$\ch\CV^{(i)}_{i-1}$ and $\Pi^k$ in the proof of Theorem \ref{prolongationPredictor} in Section \ref{DFLbySymmetrySection}.

\begin{Definition}
The first integrals $\phi^{\ell_j,j}$ of the quotient bundles \smash{$\Xi^{(j)}_{j-1}/\Xi^{(j)}$} and are known as the \textit{fundamental functions of order} $j$. We may also refer to the non-$\tau$ first integral of $\Pi^k$ as a~fundamental function of order $k$.
\end{Definition}
\begin{proc contact}[\cite{VassiliouGoursatEfficient}]\label{procContact}
\label{proc A}
Let $\mathcal{V}\subset TM$ be a Goursat bundle of signature $\langle \r_1, \r_2, \dots, \r_k\rangle$ and derived length $k>1$ such that $\r_k:=\Delta_k=1$. Then one can do the following to produce local contact coordinates for $\mathcal{V}$:
\begin{enumerate}\itemsep=0pt
\item[$(1)$] Compute \smash{$\ch\CV^{(k-1)}$} and for $j<k-1$ such that $\r_j\neq 0$ compute \smash{$\ch\CV^{(j)}$} and \smash{$\ch\CV^{(j)}_{j-1}$}.
\item[$(2)$] Identify a first integral $\tau$ of \smash{$\Char{V}{k-1}$} and a section $Z$ of $\mathcal{V}$ with the property $Z\tau=1$. Then construct $\Pi^{k}$ as in \eqref{ad char}.
\item[$(3)$] For each $1\leq j\leq k-1$ such that $\rho_j\neq 0$, compute the integrable quotient bundle
\smash{$\Xi^{(j)}_{j-1}/\Xi^{(j)}$} using step $1$.
\item[$(4)$] Compute the first integrals $\phi^{\ell_j,j}$ of \smash{$\Xi^{(j)}_{j-1}/\Xi^{(j)}$} $($fundamental functions of order $j)$.
\item[$(5)$] Define $z^{1,k}_0=\phi^{1,k}$ to be any first integral of $\Pi^{k}$ such that ${\rm d}\tau\wedge {\rm d}\phi^{1,k}\neq0$.
\item[$(6)$] For each $1\leq j\leq k$ such that $\rho_j\neq 0$, define \smash{$z^{\ell_j,j}_0=\phi^{\ell_j,j}$}, $1\leq \ell_j\leq \rho_j$. The remaining contact coordinates are
\begin{equation}\label{contact coord A}
z^{\ell_j,j}_{s_j}=Zz^{\ell_j,j}_{s_j-1}=Z^{s_j}z^{\ell_j,j}_0,\qquad 1\leq s_j\leq j,\quad 1\leq \ell_j\leq \rho_j.
\end{equation}
\end{enumerate}
The local coordinates for $J^\kappa(\mathbb{R},\mathbb{R}^m)$ are given by $\tau $, \smash{$z^{\ell_j,j}_0$}, and \eqref{contact coord A}. In these coordinates $\mathcal{V}$ has the form $\mcal{B}_\k$.

\begin{Remark}
Note the restriction in above procedure to the case $\Delta_k=1$. There is a slightly different, {\it Procedure B}, for the remaining case $\Delta_k>1$, which for simplicity of presentation, we shall not discuss in this paper. Details of this can be found in \cite{VassiliouGoursatEfficient}.
In any case, if $\Delta_k>1$, then, in practice, it is often possible to convert it to the
$\Delta_{k+1}=1$ case by performing a partial prolongation.
\end{Remark}
\end{proc contact}

Procedure {\it contact} produces a local equivalence between a Goursat bundle and a contact system. In particular, the first integral $\tau$ in procedure {\it contact} plays the role of the source variable of some $J^\kappa$, so that ${\rm d}\tau$ forms the independence condition for the linear Pfaffian system~${(J^\kappa,\bsy{\b}^\kappa)}$. Therefore, if $\mathcal{V}$ represents a control system with ${\rm d}t$ as the independence condition, then integral curves of $\mathcal{V}$ may not be sent to integrals curves of $\bsy{\b}^\kappa$ that are parameterized by $t$. The following theorem gives additional conditions under which procedure {\it contact} produces a static feedback equivalence between a Goursat bundle $\mathcal{V}$ representing a control system and a Brunovsk\'y normal form.
\begin{Theorem}[\cite{DTVa}]\label{Goursat SFL}
Let $\mathcal{V}$ be a Goursat bundle of derived length $k>1$ that represents a~control system on the manifold $M\cong_{\rm loc}\mathbb{R}\times\mathbf{X}(M)\times\mathbf{U}(M)$. Then $\mathcal{V}$ is static feedback equivalent to a~Brunovsk\'y normal form if and only if
\begin{enumerate}\itemsep=0pt
\item[$(1)$] $\inCharOne{V}=\operatorname{span}\{\partial_{u^1}, \dots, \partial_{u^m}\}$,
\item[$(2)$] ${\rm d}t\in \Xi^{(k-1)}$ if $\Delta_k=1$.
\end{enumerate}
\end{Theorem}

\begin{Remark}
 If $\CV$ is SFL, then $\t=t$ will be a first integral of $\ch\CV^{(k-1)}$ and {\it contact}
will produce a static feedback equivalence.
Again, for this paper we need only be concerned with those examples in which
$\Delta_k=1$.
\end{Remark}

We illustrate procedure {\it contact} using an example of Hunt--Su--Meyer \cite{HuntSuMeyerLin}, which was linearized via the GS algorithm in \cite{GSalgorithmExample}. As an aid to the reader, the example is presented with almost no details suppressed. Most of the calculations to follow can be automated and executed algorithmically using a software package. Throughout this paper we have used the {\sc Maple} package {\sc DifferentialGeometry}.

\begin{Example}[\cite{GSalgorithmExample,HuntSuMeyerLin}]
\begin{align*}
&\frac{{\rm d}x^1}{{\rm d}t}=\sin\bigl(x^2\bigr), \qquad\frac{{\rm d}x^2}{{\rm d}t}=\sin\bigl(x^3\bigr), \qquad\frac{{\rm d}x^3}{{\rm d}t}=\bigl(x^4\bigr)^3+u^1, \\
&\frac{{\rm d}x^4}{{\rm d}t}=x^5+\bigl(x^4\bigr)^3-\bigl(x^1\bigr)^{10}, \qquad\frac{{\rm d}x^5}{{\rm d}t}=u^2.
\end{align*}
\end{Example}
First, we will rewrite the control system as the distribution $\mathcal{V}=\operatorname{span}\{X,\partial_{u^1},\allowbreak \partial_{u^2}\}$, where
\begin{gather*}
X=\partial_t+\sin\bigl(x^2\bigr) \partial_{x^1}+\sin\bigl(x^3\bigr) \partial_{x^2}+\bigl(\bigl(x^4\bigr)^3+u^1\bigr) \partial_{x^3}+\bigl(x^5+\bigl(x^4\bigr)^3-\bigl(x^1\bigr)^{10}\bigr) \partial_{x^4}+u^2\partial_{x^5}.
\end{gather*}

\emph{Step $1$.} The derived flag of $\mathcal{V}$ is given by
\begin{align*}
\der{V}{1} =\mathcal{V}+\{\partial_{x^3},\partial_{x^5}\},\qquad
\der{V}{2}=\der{V}{1}+\{\partial_{x^2},\partial_{x^4}\},\qquad
\der{V}{3}=\der{V}{2}+\{\partial_{x^1}\}=TM.
\end{align*}
Hence $\mathcal{V}$ has derived length 3, $\vel{V}=\langle 2,2,1 \rangle$, and $\dec{V}=\langle 0,1,1 \rangle$. Since $\Delta_k=1$, we will implement Procedure~A. Next we compute the Cauchy bundles for $\der{V}{1}$ and $\der{V}{2}$. Let
\begin{equation*}
C=TX+a^1\partial_{u^1}+a^2\partial_{u^2}+b^1\partial_{x^3}+b^2\partial_{x^5}\in \der{V}{1}
\end{equation*}
be a section of the Cauchy bundle of $\der{V}{1}$, where $T$, $b^1$, $b^2$, $c^1$, and $c^2$ are smooth functions. Then
\begin{equation*}
[C,Y]\in \der{V}{1} \qquad\text{for all}\quad Y\in\der{V}{1}.
\end{equation*}
It is enough to check the Lie derivative $\LieD_C$ applied to the linearly independent sections generating~$\der{V}{1}$. Doing so, and solving for the coefficients in $C$ yields
\begin{equation*}
\Char{V}{1}=\operatorname{span}\{\partial_{u^1},\,\partial_{u^2}\}.
\end{equation*}
We now repeat this calculation for $\der{V}{2}$ and obtain
\begin{equation*}
\Char{V}{2}=\operatorname{span}\{ \partial_{u^1},\,\partial_{u^2},\,\partial_{x^3},\,\partial_{x^4},\,\partial_{x^5}\}.
\end{equation*}
From here, it is easily deduced that
\begin{align*}
\inCharOne{V}&=\operatorname{span}\{\partial_{u^1},\,\partial_{u^2}\},\qquad
\operatorname{Char}\der{V}{2}_1=\operatorname{span}\{\partial_{u^1},\,\partial_{u^2},\,\partial_{x^3},\,\partial_{x^5}\}.
\end{align*}
Thus the refined derived type of $\mathcal{V}$ is
\begin{equation*}
\mathfrak{d}_r(\mathcal{V})=[[3,0],[5,2,2],[7,4,5],[8,8]].
\end{equation*}
Checking that the relations in Proposition \ref{refined derived type numbers} are satisfied and seeing that all the Cauchy and intersection bundles are integrable, we see that $\mathcal{V}$ is a Goursat bundle. Furthermore, since~${{\rm d}t\in \ann\Char{V}{2}}$; by Theorem \ref{Goursat SFL} we deduce that $\mathcal{V}$ is SFL. Constructing the filtration of $T^*M$ (excluding the fundamental bundle) induced by $\mathcal{V}$, we find
\begin{align*}
\coder{\Xi}{2}=\operatorname{span}\bigl\{{\rm d}t,\,{\rm d}x^1,\,{\rm d}x^2\bigr\}\subset \coder{\Xi}{2}_1=&\operatorname{span}\bigl\{{\rm d}t,\,{\rm d}x^1,\, {\rm d}x^2,\,{\rm d}x^4\bigr\}\\
\subset\coder{\Xi}{1}=&\operatorname{span}\bigl\{{\rm d}t,\,{\rm d}x^1,\,{\rm d}x^2,{\rm d}x^3,\,{\rm d}x^4,\,{\rm d}x^5\bigr\}=\coder{\Xi}{1}_0.
\end{align*}

\emph{Step $2$.} Notice that $t$ is a first integral of $\Char{V}{2}$ and that $X(t)=1$. Now the fundamental bundle $\Pi^2$ is given by
\begin{equation*}
\Pi^2=\operatorname{span}\{\partial_{u^1},\,\partial_{u^2},\,\partial_{x^2},\,\partial_{x^3},\,\partial_{x^4},\,\partial_{x^5}\}.
\end{equation*}

\emph{Steps $3$ and $4$.} There is only one non-trivial quotient bundle to be computed for this step,
\begin{equation*}
\coder{\Xi}{2}_1/\coder{\Xi}{2}=\operatorname{span}\bigl\{{\rm d}x^4\bigr\},
\end{equation*}
and therefore $z_0^{1,2}=x^4$ is the fundamental function of order 2.

\emph{Step $5$.} From $\Pi^2$, we deduce the fundamental function of (highest) order 3 given by $z_0^{1,3}=x^1$ since ${\rm d}t\wedge {\rm d}x^1\neq0$. For simplicity, we shall relabel these fundamental functions as $z_0^{1,2}=z_0^1$ and~${z_0^{1,3}=z_0^2}$.

\emph{Step $6$.} Applying the final step of the procedure, we conclude that the remaining contact coordinates are
\begin{align}\label{HSM contact coords start}
z^1_1&=X\bigl(z_0^1\bigr)=x^5+\bigl(x^4\bigr)^3-\bigl(x^1\bigr)^{10},\\
z^1_2&=X\bigl(z^1_1\bigr)=u^2+3\bigl(x^4\bigr)^2\bigl(x^5+\bigl(x^4\bigr)^3-\bigl(x^1\bigr)^{10}\bigr)-10\bigl(x^1\bigr)^9\sin\bigl(x^2\bigr), \\
z^2_1&=X\bigl(z^2_0\bigr)=\sin\bigl(x^2\bigr),\qquad
z^2_2=X\bigl(z^2_1\bigr)=\cos\bigl(x^2\bigr)\sin\bigl(x^3\bigr), \\
z^2_3&=X\bigl(z^2_2\bigr)=-\sin\bigl(x^2\bigr)\sin^2\bigl(x^3\bigr)+\bigl(\bigl(x^4\bigr)^3+u^1\bigr)\cos\bigl(x^2\bigr)\cos\bigl(x^3\bigr).\label{HSM contact coords end}
\end{align}
Thus $t$, $z^1_0=x^4$, $z^2_0=x^1$, and \eqref{HSM contact coords start}--\eqref{HSM contact coords end} define a static feedback transformation of $\mathcal{V}$ to the Brunovsk\'y normal form $\bsy{\b}^{\langle 0,1,1\rangle}$.

\subsection{Control admissible symmetries}
Let $\mu\colon M\times G\to M$ be a smooth, regular right action of a Lie group $G$ on a smooth manifold~$M$ \cite{OlverLieBook,palais}. Thus the orbit space $M/G$ is a smooth manifold of dimension $\dim M-\dim G$, and~${\pi\colon M\to M/G}$ denotes the natural projection.\footnote{While regularity guarantees the quotient is a smooth manifold, it may nevertheless not have the Hausdorff separation property. In this case, we restrict to open sets where this holds. For more details, see \cite[Section~3.4]{OlverLieBook}.} The {\em quotient} of $\CV$ is
{$\CV/G:={\rm d}\pi(\CV)$}, where~${\rm d}\pi$ is the differential of $\pi$. The latter is a distribution on $M/G$, but not necessarily a~control system. One can therefore ask:
when is the $G$-quotient of $\CV$ also a control system?
To answer this we describe the appropriate Lie group action. Let $\bsy{\G}$ be the Lie algebra of infinitesimal generators of the action of $G$ on $M$. See \cite{OlverLieBook, OlverSymmetryBook} for information on Lie symmetry.

\begin{Definition}[control symmetries, \cite{DTVa}]\label{admissibleSymmetries}
Let $\mu:M\times G\to M$ be a Lie transformation group with Lie algebra $\bsy{\G}$ leaving the control system $(M, \CV)$ invariant, i.e.,
${\mu_g}_*\CV=\CV$, $\forall g\in G$ where $\mu_g(x)=\mu(x,g)$. Then $G$ is a {\it control symmetry group} if
\begin{enumerate}\itemsep=0pt
\item[(1)] $G$ acts regularly on $M$,
\item[(2)] the function $t$ is invariant, i.e.,
$\mu_g^*t=t$ for all $g\in G$, and
\item[(3)] $\operatorname{rank} \bigl({\rm d}\bsy{p}(\bsy{\G})\bigr)=\dim G$, where $\bsy{p}$ is the projection $\bsy{p}\colon M\to \B R\times \mathbf{X}(M)$ given by $\bsy{p}(t,x,u)=(t, x)$.\footnote{This ensures that the quotient control system will have $\dim\mb X(M)-\dim G$ state variables and that the number of controls will be preserved in the quotient.}
\end{enumerate}
\end{Definition}

The elements of a control symmetry group are static feedback transformations. That is, they have the form \cite[Theorem 4.9]{DTVa}
\begin{equation*}
\bar{t}=t,\qquad \bar{x}=\th(t,x), \qquad \bar{u}=\psi(t,x,u).
\end{equation*}
The class of control symmetries is essential for studying the general properties of smooth control systems under the action of a Lie group.

\begin{Lemma}
Let $\bsy{\o}$ be a smooth control system on a manifold $M$ and $\bsy{\D}$ the Lie algebra of all infinitesimal symmetries of $\bsy{\o}$. The subset $\bsy{\G}\subseteq\bsy{\D}$ of infinitesimal control symmetries is a Lie subalgebra of $\bsy{\D}$.
\end{Lemma}
\begin{proof}
A control symmetry is a static feedback transformation that is also a self equivalence. Consequently, an infinitesimal control symmetry has the form \cite[Theorem 4.9]{DTVa}
\[
X=\xi^i(t,x)\P {x^i}+\chi^a(t,x,u)\P {u^a}\in\bsy{\G}.
\]
If $Y\in\bsy{\G}$ is another infinitesimal control symmetry, then it is easy to deduce that the Lie bracket~${[X,Y]}$ also belongs to $\bsy{\G}$.
\end{proof}

\begin{Remark}
{ There is a further subalgebra
$\bsy{\Sigma}\subset\bsy{\G}$ of {\it state-space symmetries} which is better known \cite{Elkin, Grizzle}. This is the case $\chi^a\equiv 0$ in the infinitesimal generators of $\bsy{\G}$. But the restriction to $\bsy{\Sigma}$ is both unnecessary and inadequate for studying the full range of phenomena presented by control systems.}
\end{Remark}
We can now give criteria whereby the quotient (symmetry reduction) of a control system by a control symmetry group $G$ is also a control system on the quotient manifold $M/G$.

\begin{Theorem}[\cite{DTVa}]
Let $\mu\colon M\times G\to M$ be a group of control symmetries of control system~$(M, \CV)$ defined by \eqref{controlSystem}. Let $\bsy{\G}$, the Lie algebra of
$G$ satisfy
\[ \bsy{\G}\cap\CV^{(1)}=\{0\}, \]
with
$\dim G<\dim\mathbf{X}(M)$. Then the quotient $(M/G, \CV/G)$ is a control system in which
\[ \dim\mathbf{X}(M/G)=\dim\mathbf{X}(M)-\dim G, \qquad \dim\mathbf{U}(M/G)=\dim\mathbf{U}(M). \]
\end{Theorem}

\begin{Definition}\label{controlAdmissibleDefn}
An action of a Lie group $G$ on the manifold $M$ is {\bf\it control admissible} for a~control system $(M, \CV)$ defined by \eqref{controlSystem} if
\begin{enumerate}\itemsep=0pt
\item[(1)] $G$ is a control symmetry group of $\CV$,
\item[(2)] $\dim G<\dim\mb X(M)$,
\item[(3)] the action of $G$ is {\it strongly transverse}, meaning $\bsy{\G}\cap\CV^{(1)}=0$.
\end{enumerate}

\begin{Corollary}\label{controlAdmissibleCorollary}
Let $(M, \CV)$ be a control system defined by \eqref{controlSystem}. If a Lie group $G$ is control admissible for $(M, \CV)$, then the quotient $\CV/G$ is a smooth control system on the smooth quotient manifold $M/G$ in which $\dim\mb X(M/G)=\dim M-\dim G$ and $\dim \mb U(M/G)=\dim\mb U(M)$.
\end{Corollary}

\end{Definition}

\begin{Example}\label{MarinoExample} The control system on $\B R^8$ in 5 states and 2 controls,
\begin{equation*}
\CV=\operatorname{span}\bigl\{\P t+\bigl(x^5x^3+x^2\bigr)\P {x^1}+\bigl(x^5x^1+x^3\bigr)\P {x^2}+u^1\P {x^3}+x^5\P {x^4}+u^2\P {x^5} ,
\P {u^1} ,\P {u^2}\bigr\},
\end{equation*}
has a 5-dimensional Lie group of control symmetries (calculated via {\sc Maple}). It is easy to check that the subgroup $G$ generated by the Lie algebra
\[
\bsy{\G}=\operatorname{span}\bigl\{X:=x^1\P {x^1}+x^2\P {x^2}+x^3\P {x^3}+u^1\P {u^1}\bigr\}
 \]
 is control admissible (Definition \ref{controlAdmissibleDefn}).
On the $G$-invariant open set ${\rm U}\subset \B R^8$ where $x^1\neq 0$, the functions
\[ t,\quad q_1=x^2/x^1,\quad q_2=x^3/x^1,\quad q_3=x^4,\quad q_4=x^5,\quad v^1=u^1/x^1,\quad v^2=u^2 \]
are invariant under the action of $G$, which is given by
\begin{gather*}
 \bar{x}_1=\e x^1, \quad \bar{x}_2=\e x^2,\quad \bar{x}_3=\e x^3,\quad \bar{x}_4=x^4, \quad \bar{x}_5=x^5,\quad
\bar{u}^1=\e u^1,\quad \bar{u}^2=u^2,
\end{gather*}
where $\e\in G$ is an element of the multiplicative group of positive real numbers. If we denote these transformations by $\mu_\e$, then for all $\e\in G$, \smash{${\mu_\e}_*\CV_{|_x}=\CV_{|_{\mu_\e(x)}}$}. The functions
$(t, q_i, v^a)$ form a local coordinate system on an open subset of the quotient manifold
$M/G$. Furthermore, these functions on ${\rm U}\subset\B R^8$ are the components of a local representative of the projection ${\pi\colon \B R^8\to \B R^8/G}$ given by,
$\pi_{|_{{\rm U}}}(t,x,u)=(t, q_i(t,x,u), v^a(t,x,u))$. A computation then gives
\begin{gather}
{\rm d}\pi_{|_{{\rm U}}}(\CV)=\CV/G_{|_{\pi({\rm U})}}=
\operatorname{span}\bigl\{\P t-\bigl(q_1q_2q_4+q_1^2-q_2-q_4\bigr)\P {q_1}\nonumber\\
\hphantom{{\rm d}\pi_{|_{{\rm U}}}(\CV)=\CV/G_{|_{\pi({\rm U})}}=\operatorname{span}\bigl\{}{} {-}\bigl(q_2^2q_4+q_1q_2-v^1\bigr)\P {q_2}+q_4\P {q_3}+v^2\P {q_4},\, \P {v^1},\, \P {v^2}\bigr\},\label{MarinoQuotient}
\end{gather}
a smooth control system on $\pi({\rm U})\subset M/G$ in accordance with Corollary \ref{controlAdmissibleCorollary}. While $\CV/G$ has~4 states compared to the 5 states of $\CV$, its local form is more complicated, which is typical of a~symmetry reduction. On the other hand, while $\CV$ is not \sfl, it turns out that \eqref{MarinoQuotient} is \sfl. We will see later that this property of control systems $\CV$ having \sfl\ quotients has very significant consequences for the
\df\ linearizability of~$\CV$.
\end{Example}

\subsubsection{Relative Goursat bundles}\label{relativeGoursatSect}
Each Brunovsk\'y normal form $\CB_\k$ has trivial Cauchy bundle, $\ch\CB_\k$ $=\{0\}$. However, there is an important situation in which a sub-bundle can satisfy all the constraints of a Goursat bundle except for the triviality of its Cauchy bundle.

\begin{Definition}
A totally regular sub-bundle $\CV\subset TM$ is a {\it relative Goursat bundle} if it satisfies the requirements of a Goursat bundle (see Definition \ref{Goursat}) except for the triviality of its Cauchy bundle. That is, the type number $\chi^0$ need not be equal to zero in a relative Goursat bundle.
\end{Definition}
It is important to note that a relative Goursat bundle has refined derived type satisfying equations \eqref{linearTypeConstraints}.

The utility of relative Goursat bundles stems from the ability to very quickly determine the existence of a linearizable quotient $\CV/G$ of a $G$-in\-var\-i\-ant control system $\CV$. In this we always assume that the action of $G$ is control admissible. If its Lie algebra of infinitesimal generators of the action is denoted by $\bsy{\G}$,
then our assumption of strong transversality implies that $\CV\cap\bsy{\G}=0$. We often denote the direct sum $\CV\oplus\bsy{\G}$ by $\wh{\CV}$.

\begin{Theorem}[{\cite[Theorem 4.5]{DTVa}}]\label{relGoursatThm}
Suppose that the control system $\CV\subset TM$ has the control admissible Lie group $G$ with Lie algebra of infinitesimal generators $\bsy{\G}$ and satisfies $\ch\CV=0$.
If $(M, \CV\oplus\bsy{\G})$ is a relative Goursat bundle of derived length $k>1$ and signature
$\k=\operatorname{decel}( \CV\oplus\bsy{\G})$, then there is a local diffeomorphism $\varphi\colon M/G\to J^\k$ such that
$\vf_*\left(\CV/G\right)=\CB_\k$.
\end{Theorem}

An important observation is that even if $\CV$ is not a Goursat bundle, it happens very often that $\CV\oplus\bsy{\G}$ {\it is} a relative Goursat bundle and this can be very significant. However, the local diffeomorphism $\vf$ guaranteed by Theorem \ref{relGoursatThm} may not be a static feedback transformation. To guarantee the existence of such a transformation one imposes slightly more constraints on the relative Goursat bundle. The following generalizes Theorem \ref{Goursat SFL} to the case of {\it \stf}\ relative Goursat bundles and may be regarded as an {\em ``infinitesimal test''} for the existence of \sfl\ quotient systems.

\begin{Theorem}[\cite{DTVa}]\label{relStatFeedbackLin}\samepage
Let the control system $\CV\subset TM$ admit the control admissible Lie group~$G$ with Lie algebra of infinitesimal generators $\bsy{\G}$, and satisfy $\ch\CV=0$. Set $\wh{\CV}:=\CV\oplus\bsy{\G}$ and
suppose that $\bigl(M, \wh{\CV}\bigr)$ is a relative Goursat bundle of derived length $k>1$ and signature
${\k=\operatorname{decel}\bigl(\wh{\CV}\bigr)}$. Then the local diffeomorphism $\varphi\colon M/G\to J^\k$ that identifies
$\CV/G$ with its Brunovsk\'y normal form
 can be chosen to be a static feedback transformation if and only if
\begin{enumerate}\itemsep=0pt
\item[$(1)$] $\{\P {u^1}, \dots, \P {u^m}\}\subset\ch\wh{\CV}^{(1)}_0$,
\item[$(2)$] ${\rm d}t\in\operatorname{ann} \ch\wh{\CV}^{(k-1)}:=\wh{\Xi}^{(1)}_0$ if $\Delta_k=1$.
\end{enumerate}
\end{Theorem}
Note: Again, the case where $\Delta_k>1$ is not discussed.

Theorem \ref{relStatFeedbackLin} \cite[Theorem 4.12]{DTVa} is a geometric characterization of static feedback linearizable quotients of an invariant control system. It is the relative version of Theorem \ref{Goursat SFL}, applied to group quotients. In practice, an invariant control system has many SFL quotients depending on the number of subgroups that satisfy Theorem \ref{relStatFeedbackLin}.

\begin{Definition}\label{StFrelGoursatDefn}
If a relative Goursat bundle satisfies the hypotheses of Theorem \ref{relStatFeedbackLin}, then we call it a {\it \stf\ relative Goursat bundle}.
\end{Definition}
An important point is that \stf\ relative Goursat bundles are very easy to identify in practice, once the Lie algebra $\bsy{\G}$ is known. Hence, \sfl\ quotients $\CV/G$ are similarly quickly identified.

\begin{Example}\label{ExmpResolventRelativeGoursat}
This example illustrates the various constructions encountered so far. Firstly, we look at a system which is linearizable but {\it not} by \stf\ transformations.

The control system,
\begin{equation}\label{CharletExample}
\CV=\operatorname{span}\bigl\{\P t+x^2\P {x^1}+u^1\P {x^2}+u^2\P {x^3}+x^3\bigl(1-u^1\bigr)\P {x^4}, \P {u^1}, \P {u^2}\bigr\},
\end{equation}
occurs in \cite{CharletLevineMarino2}. Applying one of the well-known tests, for instance, \cite[Theorem 1]{SluisTilbury} or Theorem~\ref{Goursat SFL}, or otherwise, shows that $\CV$ is not \sfl.

However, it happens in this case that by augmenting
\eqref{CharletExample} by an integrator leads to an \sfl\ system, revealing the hidden simplicity of \eqref{CharletExample}.

But rather than performing this partial prolongation along $u^1$ \big(i.e., differentiating~$u^1$\big), let us instead reconsider the system \eqref{CharletExample} in the light of one of its {\em symmetries}. This calculation is quite revealing.

It is easy to see that $X=\P {x^4}$ is an infinitesimal symmetry of $\CV$ defined by \eqref{CharletExample}; that is, $\LieD_X\CV\subseteq\CV$. In fact, it is an infinitesimal control admissible symmetry. Let us then consider the {\it augmented} distribution $\wh{\CV}:=\CV\oplus\operatorname{span}\{X\}$,
as in Theorems \ref{relGoursatThm} and \ref{relStatFeedbackLin}. That is, $\wh{\CV}$ consists of~$\CV$ ``extended'' by its infinitesimal symmetry~$X$. The notation is justified since $\CV\cap\{X\}=\{0\}$.

While $\wh{\CV}$ is {not a control system}, the generalized Goursat normal form can be applied to any smooth sub-bundle $\CV\subset TM$ whose generic solutions are smoothly immersed curves, as in this case. Indeed, we find that the refined derived type of $\wh{\CV}$ is
\[
\mathfrak{d}_r\bigl(\wh{\CV}\bigr)=[[4, 1], [6, 3, 4], [7, 7]].
\]
This is not the refined derived type of a Brunovsk\'y normal form, however it satisfies equations~\eqref{linearTypeConstraints} with signature $\text{decel}\bigl(\wh{\CV}\bigr)=\langle 1, 1\rangle$ and derived length $k=2$. Furthermore, we have
exactly one non-trivial intersection bundle \smash{$\ch\wh{\CV}^{(1)}_0=\operatorname{span}\{\P {u^1}, \P {u^2}, X\}$} which is integrable, and
\[
\ch\wh{\CV}^{(1)}=\ch\wh{\CV}^{(1)}_0\oplus\{\P {x^3}\}.
\] Since $\r_2=1$, we check that $t$ is an invariant of
$\ch\wh{\CV}^{(1)}$, and conclude by Theorem~\ref{relStatFeedbackLin} that~$\wh{\CV}$ is a {\em \stf\ relative Goursat bundle} of signature $\langle 1, 1\rangle$ which proves that
$\CV/G$ is static feedback equivalent to $\CB_{\langle 1, 1\rangle}$, where $G$ is the Lie transformation group generated by $X=\P {x^4}$. This is an important observation to which we shall return.
\end{Example}
\begin{Remark}
While it is very easy to compute the dynamic feedback linearization
of the particular system \eqref{CharletExample}, this task is generally very difficult. However, much progress has been made recently for deriving dynamic feedback linearizations of classes of control systems by restricting the number of controls and the fiber dimension of the dynamic extension
\cite{TVFL1,Gstottner_Kolar_Schoberl_2022_Flat, Gstottner_Kolar_Schoberl_2023_Flat, Nic_Res_2016_Flatness, Nic_Res_2017_Flatness}. Our present aim is to develop a canonical, coordinate independent framework that provides tools for deriving dynamic feedback linearizations in the case of invariant control systems. Thus symmetry plays a commanding role. These tools permit one to study control systems in which the number of controls and size of the dynamic extension are not restricted. To illustrate this, we work through a number of elementary examples throughout the paper. In particular, it will be illuminating to continue using \eqref{CharletExample} as a simple running example that encapsulates the basic ideas.
\end{Remark}

\section{The contact sub-connection}\label{ConnectionSect}

In previous sections, we have seen how one can often perform a symmetry reduction of a control system $(M, \CV)$ so that the resulting quotient is again a control system and is, importantly, \sfl. Ultimately, we wish to use this property to construct an explicit solution for the originally given control system $\CV$, which generally will not itself be \sfl. We call this procedure of constructing an explicit solution for $\CV$ from that of its quotient, {\em cascade feedback linearization}.
In this section, we begin to explain how this can be achieved by describing a key geometric object in the theory.

\begin{Definition}\label{subconnection}
Let $G$ act regularly on $M$ on the right, $\pi\colon M\to M/G$ be the right principal~$G$-bundle and
$\text{V}M$ the vertical bundle
$\ker {\rm d}\pi$. Let a given sub-bundle $\Pi^G\subseteq T(M/G)$ together with a constant rank distribution
$M\ni p\mapsto H_{p}\subset T_p M$ satisfy
\begin{enumerate}\itemsep=0pt
\item[(1)] $H_p\cap \text{V}_p M= \{0\}$,
\item[(2)] ${\rm d}\pi(H_p)=\Pi_{\pi(p)}^G$,
\item[(3)] ${\mu_g}_*H_p=H_{p\cdot g}$, with $\mu_g(p):=p\cdot g$ being the right $G$-action,
\item[(4)] $p\mapsto H_p$ is smooth,
\end{enumerate}
$\forall g\in G$, $p\in M$. Then $H$ will be called a {\em $($right$)$ principal sub-connection relative to} {$\Pi^G$.}
\end{Definition}

Evidently, this is the usual definition of a principal connection (see \cite[Chapter II]{KobNomizu}) when
${\Pi^G=T(M/G)}$. A curve in $M/G$ passing through the point $q\in M/G$ and all of whose tangent vectors belong to $\Pi^G$ has a unique lifting to
$M$ passing through a prescribed point $p\in\pi^{-1}(q)$ such that the lifted curve is an integral submanifold of $H$.
This is all we shall require for the present time.

Let $\bsy{\G}$ denote the Lie algebra of infinitesimal generators for the action of control admissible symmetries $G$ on $M$. By Theorem \ref{relStatFeedbackLin}, there is a quotient control system
$(M/G, \bsy{\o}/G)$ that is static feedback equivalent to a Brunovsk\'y normal form if
$\CV\oplus\bsy{\G}\subset TM$ is a static feedback relative Goursat bundle (cf.\ Definition \ref{StFrelGoursatDefn}).
We will now explain how this can, in principle, be used to compute the integral manifolds of $\bsy{\o}$ from those of
$\bsy{\o}/G$.

If $G$ is a control admissible symmetry group of $\bsy{\o}$, then finding an SFL quotient $\bsy{\o}/G$ in fact determines a local isomorphism of principal bundles as follows. Since $\bsy{\o}/G$ is SFL, it can be identified by a static feedback transformation
$\varphi\colon{\rm U}\subseteq M/G\to J^\k$ with a Brunovsk\'y normal form $\bsy{\b}^\k$ of some signature $\k$ via the pullback
$\bsy{\o}/G = \varphi^*\bsy{\b}^\k$.

The local diffeomorphism $\varphi$ can be lifted to a local \stf\ transformation
$\wt{\varphi}\colon M\to J^\k\times G$ that realizes
the $G$-principal bundle $\pi\colon M\to M/G$ as locally isomorphic to the trivial bundle
$\pi'\colon J^\k\times G\to J^\k$ via a principal bundle isomorphism $(\varphi, \wt{\varphi})$ as depicted in the commuting diagram of Figure \ref{diagramBasic}. That $\wt{\bsy{\varphi}}$ is a static feedback transformation follows from the fact that $\varphi$ is a static feedback transformation and that $G$ acts by static feedback transformations. In coordinates, the map $\wt{\varphi}$ is determined by
\begin{equation}\label{bundle sfl}
\wt{\varphi}\colon \ (t,\bsy{x},\bsy{u})\mapsto (\varphi(t,\bsy{q}(t,\bsy{x}),\bsy{v}(t,\bsy{x},\bsy{u})),\e(t,\bsy{x})),
\end{equation}
where $\e(t,\bsy{x})$ are local coordinates about the identity in $G$ chosen so that the infinitesimal generators $\{X_a\}_{a=1}^r$ for the action of $G$ on $M$ may be expressed as vector fields on $G$ with local coordinates $\e^a$. To achieve this, functions $\e(t,\bsy{x})$ may be chosen by solving for the
$\bigl(\e^1,\dots,\e^r\bigr)$ parameters that appear in the action of $G$ on $M$ when restricted to an appropriate choice of section of $\pi\colon M\to M/G$. For this one first computes the transformation group by solving the ODE systems defined by the infinitesimal generators so that each local Lie group coordinate $\e^a$ is the flow coordinate for a corresponding infinitesimal generator.
By construction, there exists a unique Pfaffian system $\bsy{\g}^G$ on $J^\k\times G$ with the property that $\wt{\varphi}^*\bsy{\g}^G=\bsy{\o}$.
The dual sub-bundle
\[ \CH_G = \wt{\varphi}_*(\CV) \subset T(J^\k\times G) \]
is a principal sub-connection relative to the contact distribution $\mcal{B}_\k$ on $J^\k$.
\begin{Definition}
Let $(\varphi, \wt{\varphi})$ be the principal bundle isomorphism for the principal bundle $\pi\colon M\!\to M/G$. We call the principal sub-connection
\[ \CH_G = \wt{\varphi}_*(\CV) \subset T(J^\k\times G) \]
\big(and its dual $\bsy{\g}^G$\big) the {\it contact sub-connection} on $J^\k\times G$ of $(M, \CV)$.
\end{Definition}

Now we present a theorem which gives a coordinate normal form for a contact sub-connection arising from a control admissible symmetry group $G$ and control system $(M,\mcal{V})$.

\begin{Theorem}\label{ConstructingConnectionThm}
Suppose $G$ is a Lie group of control admissible symmetries of a control system~${(M, \CV)}$ acting on the right and $\bsy{\G}$ is its Lie algebra of infinitesimal generators. Assume~$\CV/G$ is \sfl\ of signature $\k$. Let the principal bundle equivalence $(\varphi, \wt{\varphi})$ be a local trivialization of $\pi\colon M\to M/G$ to
$\pi'\colon J^\k\times G\to J^\k$. Then the contact sub-connection $\CH_G=\wt{\varphi}_*\CV$ has the following local normal form:
\begin{equation}\label{contactConnection}
{\CH}_G=\operatorname{span}\left\{\mathbf{D}_t+\sum_{a=1}^r\lambda^a(t,z^\k)R_a,\,
 \P{\bsy{z}^{\ceil{\k}}}\right\},
\end{equation}
for some functions $\lambda^a$ on $J^\k$, where $\mathbf{D}_t$ is the total differential
operator \eqref{totalDiffOp} on $J^\k$, $\{R_a\}_{a=1}^r$ is a basis for the Lie algebra of right-invariant vector fields on $G$ and $\operatorname{span}\left\{\mathbf{D}_t, \P {z^{[\k]} }\right\}=\CB_\k$ is the Brunovsk\'y normal form of signature $\k$.
\end{Theorem}
\begin{proof}
Let $\CV=\operatorname{span}\bigl\{\P t+f^i(t,x,u)\P {x^i}, \P {u^1},\dots,\P {u^m}\bigr\} $ be a smooth control system, invariant under a control admissible Lie group $G$ with SFL quotient control system $\CV/G$ via a~map ${\varphi\colon M/G \to J^\k}$. Let $(\varphi, \wt{\varphi})$ be a trivializing bundle isomorphism between $\pi\colon M \to M/G$ and $\pi'\colon J^\k\times G \to J^\k$. Since ${\rm d}\pi'(\CH_G)=\CB_\k$, and because the rank of a bundle is a diffeomorphism invariant, the principal sub-connection $\CH_G$ on $J^\k\times G$ must be of the form
\begin{equation*}
{\CH}_G=\{\mathbf{D}_t+X,
\P {\bsy{z}^{\ceil{\k}}}+Y_{\bsy{z}^{\ceil{\k}}}\},
\end{equation*}
where $\{X,Y_{\bsy{z}^{\ceil{\k}}}\}$ are vector fields tangent to the fibers of $\pi'$.
Furthermore, since the intersection bundle \smash{$\ch\CV^{(1)}_{0}=\operatorname{span}\{\P{u^1},\dots,\P{u^m}\}$} is integrable and a diffeomorphism invariant, then
\[
\wt{\varphi}_*\ch\CV^{(1)}_{0}=\ch\left(\CH_G\right)^{(1)}_0=\operatorname{span}\{ \P {\bsy{z}^{\ceil{\k}}}+Y_{\bsy{z}^{\ceil{\k}}}\}
\]
must also be integrable. In particular, since $\wt{\varphi}$ is a static feedback transformation,
\begin{align*}
\bigl(\wt{\varphi}^{-1}\bigr)^*\Xi^{(1)}_0&=\operatorname{span}\bigl\{{\rm d}t,\,{\rm d}\bsy{x}\bigl(t,\bsy{z}^{\floor{\k}},\bsy{\e}\bigr)\bigr\}
 =\operatorname{span} \bigl\{{\rm d}t,\,{\rm d}\bsy{z}^{\floor{\k}},\,{\rm d}\e\bigr\},
=\operatorname{ann} \ch(\CH_G)_0^{(1)}.
\end{align*}
Therefore, \smash{$\ch(\CH_G)_0^{(1)}=\operatorname{span}\{\P {\bsy{z}^{\ceil{\k}}}\}$} and each of the $Y_{\bsy{z}^{\ceil{\k}}}$ must be zero.
Lastly, since $X\in \ker({\rm d}\pi)$, then we can write
\begin{equation*}
X=\sum_{a=1}^r\lambda^a(t,\bsy{z},\bsy{\e})R_a,
\end{equation*}
where $\{R_a\}_{a=1}^r$ form the basis of right invariant vector fields on $G$.
The natural action of $G$ on $J^\k\times G$ is precisely $\wt{\mu}_g(t,\bsy{z},\bsy{\e})=(t,\bsy{z},\e\cdot g)$ where $\e\cdot g$ is multiplication on $G$ with $g$ on the right. The infinitesimal generators of the action of $G$ on $M$ must push forward under $\wt{\varphi}$ to the infinitesimal generators of the action of $G$ on $J^\k\times G$. As such, the infinitesimal generators~$\bsy{\wt{\Gamma}}$ of~$\wt{\mu}_G$ must be precisely the left-invariant vector fields $\bsy{\wt{\Gamma}}=\operatorname{span}\{L_1,\dots, L_r\}$. The contact sub-connection $\CH_G$ is invariant under the infinitesimal action $\bsy{\wt{\G}}$ of $G$ on $J^\k\times G$ and as such we must have that $[L_b, \CH_G]\subset \CH_G$ implying that
\begin{equation*}
[L_b, \mathbf{D}_t+X] =\sum_{a=1}^rL_b\left(\lambda^a(t,\bsy{z},\bsy{\e})\right)R_a+\lambda^a(t,\bsy{z},\bsy{\e})[L_b,R_a]
 =\sum_{a=1}^rL_b\left(\lambda^a(t,\bsy{z},\bsy{\e})\right)R_a,
\end{equation*}
since $[L_b,R_a]=0$ for all $1\leq a, b\leq r$. Therefore, in order that $[L_b, \mathbf{D}_t+X]\in \CH_G$ we must have
\begin{equation*}
L_b\left(\lambda^a(t,\bsy{z},\bsy{\e})\right)=0
\end{equation*}
for all $1\leq a, b\leq r$ and therefore each $\lambda^a$ has no dependence on $\e$. This provides the claimed normal form.
\end{proof}

Theorem~\ref{ConstructingConnectionThm} plays an important role in Section~\ref{DFLbySymmetrySection}, where a method for constructing dynamic feedback linearizations using symmetry is described. More generally, quite apart from
\df\ linearization, it exhibits, in group theoretic terms, an obstruction to the \stf\ linearizability of control systems with symmetry.
\begin{notate}
It is sometimes convenient to use the differential forms version of the contact sub-connection. We may write $\bsy{\gamma}^G=\operatorname{ann} \CH_G$ as
\begin{equation*}
\bsy{\gamma}^G=\bsy{\b}^\k\oplus \Theta^G
\end{equation*}
with
\begin{equation*}
\Theta^G=\operatorname{span}\{ \eta^a - \lambda^a(t,\bsy{z}^\k) {\rm d}t\}_{a=1}^r,
\end{equation*}
where each $\eta^a$ is an entry of the right-invariant Maurer--Cartan form on $G$, i.e., $\eta^a(R_b)=\delta^a_b$, where $\delta^a_b$ is the Kronecker delta.
\end{notate}

In order to construct solution curves for the original system $(M, \bsy{\o})$, we start with an arbitrary {\it contact curve} $c\colon \B R\to J^\k$; that is, an {integral curve of the} contact system $\bsy{\b}^\k$ on $J^\k$. The {\it horizontal lift} $\wt{c}\colon\B R\to J^\k\times G$ of any contact curve to an integral curve of the contact sub-connection $\CH_G$ on $J^\k\times G$ which corresponds to an integral curve (i.e., a solution) of the system~${(M, \bsy{\o})}$ via composition with the map $\wt{\varphi}^{-1}$.

In general, constructing the lifted curve $\wt{c}\colon \B R\to J^\k\times G$ might be expected to require the solution of a system of ODE, as in the reconstruction theorem of \cite{AF}. Even in the best-case scenario where $G$ is a solvable Lie group, the solution would be expected to require quadratures. But, as we will show in the following sections, if the contact sub-connection admits a {\em reduction by partial contact curves} to a system which is also \stf\ linearizable, then it will be possible to obtain an explicit solution for the system $(M, \bsy{\o})$, which is quadrature free by definition.

\begin{figure}[ht]
\centering
\begin{tikzcd}[every label/.append style = {font = \normalsize}]
{(M,\boldsymbol{\omega})} \arrow[d, "\pi"'] \arrow[r, "\widetilde{\varphi}"] & {(J^\kappa\times G,\boldsymbol{\gamma}^G)} \arrow[d, "\pi'"] & \\
{(M/G,\boldsymbol{\omega}/G)} \arrow[r, "\varphi"] & {(J^\kappa,\bsy{\b}^{\kappa})} & \mathbb{R} \arrow[l, "c"] \arrow[lu, "\widetilde{c}"']
\end{tikzcd}
\caption{Principal bundle equivalence where
$\bsy{\g}^G=\operatorname{ann}\CH_G$ and $c\colon\mathbb{R}\to J^\k$ is an integral curve of $\bsy{\b}^\k$.}
\label{diagramBasic}
\end{figure}

\begin{Example}\label{contactConnectionExample}
Example \ref{ExmpResolventRelativeGoursat} shows that
$\wh{\CV}:=\CV\oplus\{\P {x^4}\}$, consisting of the control system $\CV$ defined by \eqref{CharletExample} and the symmetry $X=\P {x^4}$, is a \stf\ relative Goursat bundle of signature
$\langle 1, 1\rangle$, and hence $\CV/G$ is static feedback equivalent to $\CB_{\langle 1, 1\rangle}$ by Theorem \ref{relStatFeedbackLin}.

Here we construct the principal bundle equivalence $(\varphi, \wt{\varphi})$ of Figure \ref{diagramBasic} and the contact sub-connection $\CH_G$ of $\CV$, in local coordinates. As indicated above this can be done by first computing~$\CV/G$ and then applying {\it contact} to the result. In this case, since $G$ is generated by $\partial_{x^4}$, it is easy to see by inspection that
\begin{equation}\label{quotientBasis11}
\CV/G=\operatorname{span}\bigl\{\P t+x^2\P {x^1}+u^1\P {x^2}+u^2\P {x^3},\, \P {u^1},\, \P {u^2}\bigr\},
\end{equation}
in which, for simplicity, we have used identical labels for the coordinates on $\B R^7/G$ and on $\B R^7$. We compute
\[
\ch(\CV/G)^{(1)}_0=\operatorname{span}\{\P {u^1},\, \P {u^2}\},\qquad \ch (\CV/G)^{(1)}=\operatorname{span}\{\P {x^3}, \,\P {u^1},\, \P {u^2}\},
\]
and the latter is annihilated by ${\rm d}t$. Hence $\CV/G$ is a \stf\ Goursat bundle by Theorem~\ref{Goursat SFL}, as predicted in Example
\ref{ExmpResolventRelativeGoursat}. Procedure {\it contact} relies on finding the fundamental functions of all relevant orders; in this case one fundamental function of order 1 and one of order~2, in accordance with the signature
$\langle 1, 1\rangle$. The order 2 fundamental function is any invariant (first integral) $\phi^2$ of the fundamental bundle $\Pi^2$ of \eqref{quotientBasis11} such that
${\rm d}t\wedge {\rm d}\phi^2\neq 0$.
Since the derived length of $\CV/G$ is $k=2$, we have
\[
\Pi^2=\operatorname{span}\{\P {u^1}, \,\P {u^2},\, [Z, \P {u^1}],\, [Z,\P {u^2}]\}=\operatorname{span}\{ \P {u^1},\, \P {u^2},\,\P {x^2},\, \P {x^3}\},
\]
the invariants of which are spanned by $t$ and $x^1$; hence $\phi^2=x^1$. The order 1 fundamental function $\phi^1$, is any invariant of the integrable quotient bundle
\[ \Xi^{(1)}_0/\Xi^{(1)}=\operatorname{span}\bigl\{{\rm d}t,\, {\rm d}x^1,\, {\rm d}x^2,\, {\rm d}x^3\bigr\}\slash\bigl\{{\rm d}t,\, {\rm d}x^1,\, {\rm d}x^2\bigr\}=\operatorname{span}\bigl\{ {\rm d}x^3\bigr\}, \]
where
\[ \Xi^{(1)}_0=\operatorname{ann}\ch(\CV/G)^{(1)}_0, \qquad \Xi^{(1)}=\operatorname{ann}\ch(\CV/G)^{(1)}. \]
Thus a representative basis of this quotient bundle is $\bigl\{{\rm d}x^3\bigr\}$, and hence $\phi^1$ can be taken to be equal to $x^3$. Setting $z=x^3$ and $w=x^1$ and using $t$, $z$, $z_1$, $w$, $w_1$, $w_2$ for the standard coordinates on
$J^{\langle 1, 1\rangle}$, we generate the remaining contact coordinates for $\CV/G$ as
\[
z_1=\LieD_Zz=u^2,\qquad w_1=\LieD_Zw=x^2,\qquad w_2=\LieD^2_Zw=u^1.
\]
These functions are components of a local diffeomorphism $\varphi\colon M\to J^{\langle 1, 1\rangle}$, where $Z$ can be taken to be the first element in the basis \eqref{quotientBasis11} and $M\subset\B R^7/G$ is some open set.
Explicitly, $\varphi$ has the form
\[ t=t,\qquad z=x^3,\qquad z_1=u^2,\qquad w=x^1,\qquad w_1=x^2,\qquad w_2=u^1. \]
The function $\e=x^4$ is a local coordinate on $G$, which together with the components of $\varphi$ define a local diffeomorphism\footnote{We show in Section \ref{dynamicExtensionsSection} how $\varphi$ is lifted to
$\wt{\varphi}$ in the general case.}
$\wt{\varphi}\colon M\to J^{\langle 1, 1\rangle}\times G$ given by
\begin{equation}\label{sigmaEq-Exmp}
t=t,\qquad z=x^3,\qquad z_1=u^2,\qquad w=x^1,\qquad w^1=x^2,\qquad w^2=u^1, \qquad \e=x^4.
\end{equation}
Thus, we get the contact sub-connection for system \eqref{CharletExample} in the form
\begin{equation*}
\CH_G:=\wt{\varphi}_*\CV=\operatorname{span}\{\P t+z_1\P z+w_1\P w+w_2\P {w_1}+z(1-w_2)\P {\e},\,\P {z_1}, \,\P {w_2}\},
\end{equation*}
which may be compared to \eqref{contactConnection}. Dually, we have $\bsy{\g}^G:=\operatorname{ann}\CH_G$ is given by
\begin{align}
\bsy{\g}^G&=\operatorname{span}\left\{{\rm d}z-z_1{\rm d}t,\, {\rm d}w-w_1{\rm d}t,\, {\rm d}w_1-w_2{\rm d}t\right\}\oplus\{{\rm d}\e-z(1-w_2){\rm d}t\}\nonumber\\
&=\bsy{\b}^{\langle 1, 1\rangle}\oplus\{{\rm d}\e-z(1-w_2){\rm d}t\}.\label{connectionFormCharlet}
\end{align}
Here, $G$ is isomorphic to $(\B R,+)$, so $r:=\dim G=1$,\ and
\[
\lambda^1\bigl(t,z^{\langle 1, 1\rangle}\bigr)=z(1-w_2).
\]
\end{Example}

While the local normal form \eqref{contactConnection} is \stf\ equivalent to the originally given control system $(M, \CV)$, it embodies very significant advantages over the latter. One of these is that it permits a reduction by partial contact curves that induces the dynamic feedback linearization of $\CV$, as we will demonstrate in subsequent sections.

\section{Cascade feedback linearization}\label{contactReduction}

So far, we have presented the construction of the integral submanifolds of control system $\bsy{\o}$ as a two-step process, as in Figure \ref{diagramBasic}:
\begin{enumerate}\itemsep=0pt
\item[(1)] Find a static feedback linearizable quotient $(M/G, \bsy{\o}/G)$ of $(M, \bsy{\o})$ by a control admissible subgroup $G$ of its Lie group of control symmetries, and construct the corresponding diffeomorphism
\[ \wt{\varphi}\colon\ (M, \bsy{\o}) \to \bigl(J^\k \times G, \bsy{\g}^G\bigr). \]
\item[(2)] For a general contact curve $c\colon\B R \to (J^\k, \bsy{\beta}^\k)$, construct the horizontal lift
\[\wt{c}\colon \ \B R \to \bigl(J^\k \times G, \bsy{\g}^G\bigr).\]
\end{enumerate}

This section will be devoted to the construction in the second step, originally introduced in~\cite{VassiliouCascade1}. No novel results appear in this section and the reader may reference \cite{KlotzThesis} for additional results and exposition.

For $m\geq 2$, we can always rewrite the Brunovsk\'y normal form $\bsy{\b}^\k$ as $\bsy{\b}^\k=\bsy{\b}^\nu\oplus\bsy{\b}^{\nu^\perp}$, where~${\kappa=\nu+\nu^\perp}$ and $m=m_\nu+m_{\nu^\perp}$, so that $\bsy{\b}^\nu$ and $\bsy{\b}^{\nu^\perp}$ are the canonical contact systems on $J^\nu(\mathbb{R},\mathbb{R}^{m_\nu})$ and $J^{\nu^\perp}(\mathbb{R},\mathbb{R}^{m_{\nu^\perp}})$, respectively. For example, take $\k=\langle 1,2 \rangle=\langle 1,1 \rangle +\langle 0, 1 \rangle$ with $\nu=\langle 1,1 \rangle$ and $\nu^\perp=\langle 0,1 \rangle$. Then the Brunovsk\'y normal form $\bsy{\beta}^\k$ may be decomposed as
\begin{align*}
\bsy{\b}^\k&=\operatorname{span}\bigl\{ {\rm d}z^1_0-z^1_1 {\rm d}t,\,{\rm d}z^2_0-z^2_1 {\rm d}t,\,{\rm d}z^2_1-z^2_2 {\rm d}t,\,{\rm d}z^3_0-z^3_1 {\rm d}t,\, {\rm d}z^3_1-z^3_2 {\rm d}t\bigr\}\\
 &=\operatorname{span}\bigl\{ {\rm d}z^1_0-z^1_1 {\rm d}t,\, {\rm d}z^2_0-z^2_1 {\rm d}t,\, {\rm d}z^2_1-z^2_2 {\rm d}t\bigr\} \oplus \operatorname{span}\bigl\{ {\rm d}z^3_0-z^3_1 {\rm d}t,\, {\rm d}z^3_1-z^3_2 {\rm d}t\bigr\}\\
 &=\bsy{\b}^\nu\oplus \bsy{\b}^{\nu^\perp}.
\end{align*}

\begin{Definition}\label{contactReductionDef}
We say that a submanifold $\Sigma^\nu_f\subset J^\kappa\times G$ is a {\bf\it codimension $s$ partial contact curve} of \smash{$\bsy{\b}^\k
=\bsy{\b}^\nu\oplus\bsy{\b}^{\nu^\perp}$} if $\Sigma^\nu_f$ is an integral manifold of $\bsy{\b}^\nu$ and $s$ is the sum of the entries in $\nu$. It is the image of a map
\[\mathbf{c}_f^{\nu}\colon\ J^{\nu^{\perp}} \times G \to J^\k \times G \]
defined by
\[ \mathbf{c}_f^{\nu}\bigl(t, \bsy{z}^{\nu^{\perp}}, \e\bigr) = \bigl(t, j^\nu f(t), \bsy{z}^{\nu^{\perp}}, \e\bigr), \]
where $\bsy{z}^{\nu^{\perp}}$ represents the local contact coordinates on $J^{\nu^{\perp}}$, $j^\nu f(t)$ represents the integral curve of $\bsy{\beta}^\nu$ corresponding to the jet of some smooth function $f\colon \mathbb{R}\to\mathbb{R}^{m_\nu}$, and $\e$ represents local coordinates on $G$.

In particular, we refer to a system $\bsy{\g}^G$ restricted to the {\it family} of such submanifolds as $f$ ranges over the space of {\it generic} smooth functions as a {\bf\it partial contact curve reduction} of $\bsy{\g}^G$ and denote it by $\bsy{\bar{\g}}^G$.
\end{Definition}

In general, the local coordinate expression for $\bsy{\bar{\g}}^G$ will be $t$-dependent via explicit dependence on the arbitrary function $f$ and its derivatives up to some finite order, and its geometry will be quite different from that of $\bsy{\g}^G$. In particular, $\bsy{\bar{\g}}^G$ may turn out to be static feedback linearizable even when the original connection $\bsy{\g}^G$ is not.
Matters being so, the second step of the cascade feedback linearization process takes place under the assumption that the partial contact curve reduction can be chosen in such a way that $\bsy{\bar{\g}}^G$ is feedback linearizable. In this case, an explicit linearization for this system provides a path to an explicit formula for the desired lifted curves~${\wt{c}\colon \B R \to \bigl(J^\k \times G, \bsy{\g}^G\bigr)}$ that requires no integration.

We remark that the generic family of smooth functions that define a partial contact curve are truly generic in the sense that certain ODE solutions must be avoided. See \cite{KlotzThesis} for more.

\begin{Definition}[cascade feedback linearization]\label{cfl-def}
 A control system $(M, \bsy{\o})$ invariant under a~Lie group of control admissible transformations $G$ is called {\it cascade feedback linearizable} if
\begin{enumerate}\itemsep=0pt
\item[(1)] $(M/G, \bsy{\o}/G)$ is \stf\ linearizable, and
\item[(2)] the contact sub-connection $\bsy{\g}^G$ admits a static feedback linearizable partial contact curve reduction $\bsy{\bar{\g}}^G$.
\end{enumerate}
\end{Definition}

The above definition of cascade feedback linearization is visualized in Figure~\ref{CFL diagram}.

\begin{figure}[ht]
\hspace*{1mm}\begin{tikzcd}[every label/.append style = {font = \normalsize}]
{(M,\boldsymbol{\omega})} \arrow[r, "\wt{\varphi}"] \arrow[d, "\pi"'] \arrow["G"', loop, distance=2em, in=215, out=145] & {\bigl(J^\kappa\times G,\boldsymbol{\gamma}^G\bigr)=\bigl(J^{\nu+\nu^\perp}\times G,\bsy{\gamma}^G\bigr)} \arrow[d, "\pi'", shift right=13] \arrow[d, "\pi'", shift left=10] & {\bigl(J^{\nu^\perp}\times G,\bar{\boldsymbol{\gamma}}^G\bigr)} \arrow[l, "\boldsymbol{c}^\nu_{\boldsymbol{f}}", hook'] \arrow[r, "\bar{\varphi}"] \arrow[r, "(2)", phantom, shift right=3] & {\bigl(J^{\bar{\kappa}},\boldsymbol{\beta}^{\bar{\kappa}}\bigr)} \\
{(M/G,\boldsymbol{\omega}/G)} \arrow[r, "\varphi"] \arrow[r, "(1)", phantom, shift right=3] & { \bigl(J^\kappa, \boldsymbol{\beta}^\kappa \bigr) = \bigl(J^{\nu+\nu^\perp},\bsy{\b}^{\nu+\nu^\perp}\bigr)} & &
\end{tikzcd}
\caption{Items (1) and (2) in Definition \ref{cfl-def} are labelled by (1) and (2) in this diagram.}
\label{CFL diagram}
\end{figure}

It follows by construction that a \cfl\ control system is explicitly integrable.
A detailed illustration of the cascade feedback linearization process is given in the following example.

\begin{Example}\label{cflExample} Example \ref{contactConnectionExample} computed the contact connection $\CH_G:=\wt{\varphi}_*\CV$ of system \eqref{CharletExample}. Using this we now give a simple example of cascade feedback linearization.

Having computed the sub-connection $\CH_G$, the next step is to check for the existence of a~\sfl\ partial contact curve reduction of $\CH_G$.
By inspection, we find that~\smash{$\mb c^{\langle 0, 1\rangle}_f\colon J^{\langle 1\rangle}\times G\to J^{\langle 1,1\rangle}\times G$} given by
\[ \mb c^{\langle 0, 1\rangle}_f(t, z, z_1, \e)=\bigl(t, z, z_1, \e, w=f(t), w_1=\dot{f}(t), w_2=\ddot{f}(t)\bigr) \]
for an arbitrary smooth function $f\colon \B R\to\B R$ has the property that
\[
\bigl(\mb c^{\langle 0, 1\rangle}_f\bigr)^*\bsy{\g}^G:=\bsy{\bar{\g}}^G=\operatorname{span}\bigl\{{\rm d}z-z_1{\rm d}t,\,
{\rm d}\e-z\bigl(1-\ddot{f}(t)\bigr){\rm d}t\bigr\}
\]
This system is static feedback linearizable for arbitrary $f$, since,
letting $\bar{\CH}_G:=\ker\bsy{\bar{\g}}^G$, we have
\[
\bar{\CH}_G=\operatorname{span}\{\P t+z_1\P z+z\bigl(1-\ddot{f}(t)\bigr)\P \e, \P {z_1}\},
\]
and we compute that
\[ \mathfrak{d}_r\bigl(\bar{\CH}_G\bigr)=[[2, 0], [3, 1, 1], [4, 4]]. \]
Therefore, we have
$\text{decel}\bigl(\bar{\CH}_G\bigr)=\langle 0, 1\rangle$. Since ${\rm d}t$ annihilates $\ch\bar{\CH}_G^{(1)}=\operatorname{span}\{\P {z_1}\}$, we conclude that $\bar{\CH}_G$ is \stf\
linearizable (see Theorem \ref{Goursat SFL}), hence explicitly integrable.
This proves that \eqref{CharletExample} is \cfl\ according to Definition \ref{cfl-def}.

We can apply cascade feedback linearization to construct an explicit solution to the original system as follows.
Since we can compute an explicit solution $\varsigma\colon \B R\to J^{\langle 1\rangle}\times G$ of
$\bsy{\bar{\g}}^G$,
we have that
\[ \mb c^{\langle 0, 1\rangle}_f\circ\varsigma\colon\ \B R\to J^{\langle 1,1\rangle}\times G \]
is an explicit solution of
$\bsy{\g}^G$. Therefore,
\[ \bigl(\wt{\varphi}\bigr)^{-1}\circ\mb c^{\langle 0, 1\rangle}_f\circ\varsigma\colon\ \B R\to M \]
is an explicit solution of $(M, \bsy{\o})$ . Thus one obtains an explicit solution of the non-\sfl\ system
$(M, \bsy{\o})$ by a ``cascade'' of the explicit solutions of a pair of \sfl\ systems, $\bsy{\o}/G$ and
$\bsy{\bar{\g}}^G$.

Let us explicitly demonstrate these calculations; this will also provide another illustration of {\it contact} to find the explicit solution of $\bsy{\bar{\g}}^G$, hence of
$\bsy{\g}^G$ and finally of $\bsy{\o}$. Since $\r_2=1$, we compute the highest order bundle
\[ \Pi^2=\operatorname{span}\{\P {z_1}, [\P {z_1}, Z]\}=\operatorname{span}\{\P {z_1}, \P z\}, \]
which has first integrals $\{t, \e\}$, where $Z$ is the first vector field in the given basis for
$\bar{\CH}_G$.
Hence~$\e$ spans the fundamental functions of order 2, from which we compute contact coordinates by Lie differentiation by $Z$:
\begin{equation}\label{reducedSolution}
a=\e, \qquad a_1=\LieD_Z\e=z\bigl(1-\ddot{f}\bigr), \qquad a_2=\LieD^2_Z\e=-z\dddot{f}+z_1\bigl(1-\ddot{f}\bigr).
\end{equation}
These are the components of the \stf\ transformation
\[
\bar{\varphi}\colon\ J^{\langle 1\rangle}\times G\to J^{\langle 1,1\rangle}
\]
that identifies $\bar{\CH}_G$ with $\CB_{\langle 0, 1\rangle}:=\operatorname{span}\{\P t+a_1\P a+a_2\P {a_1},\, \P {a_2}\}$.
Since for any smooth real-valued function $h$, the functions
\[ a=h(t),\qquad a_1=\dot{h}(t),\qquad a_2=\ddot{h}(t) \]
form an explicit solution of the image $\CB_{\langle 0, 1\rangle}$ of $\bar{\CH}_G$ under the map $\bar{\varphi}(t,z,z_1,\e)=(t, a, a_1, a_2)$ defined by \eqref{reducedSolution}, we can invert this to obtain the explicit solution $\varsigma$ of $\bar{\CH}_G$,
\[
\e=h, \qquad z=\frac{\dot{h}}{1-\ddot{f}},\qquad z_1=\frac{\rm d}{{\rm d}t}\left(\frac{\dot{h}}{1-\ddot{f}}\right).
\]
Thus, the explicit solution of $\bsy{\g}^G$ is $\mb c^{\langle 0, 1\rangle}_f\circ\varsigma$, given by
\[
\e=h, \qquad z=\frac{\dot{h}}{1-\ddot{f}},\qquad z_1=\frac{\rm d}{{\rm d}t}\left(\frac{\dot{h}}{1-\ddot{f}}\right),\qquad
w=f,\qquad w_1=\dot{f},\qquad w_2=\ddot{f},
\]
from which we obtain explicit solution $\wt{\varphi}^{-1}\circ\mb c^{\langle 0, 1\rangle}_f\circ\varsigma$ of
$\CV$, given by
\[
x^1=f,\qquad x^2=\dot{f}, \qquad x^3=\frac{\dot{h}}{1-\ddot{f}},\qquad x^4=h,\qquad u^1=\ddot{f},\qquad u^2=\frac{\rm d}{{\rm d}t}\left(\frac{\dot{h}}{1-\ddot{f}}\right),
\]
upon referring to $\wt{\varphi}$ in \eqref{sigmaEq-Exmp} from Example \ref{contactConnectionExample}.

\end{Example}

\subsection{How hard is it to implement cascade linearization?}
\label{implementingCFLsect}
Cascade linearization relies on constructing a Lie group $G$ of control symmetries for a given control system
$\CV$, and then a subgroup $H\subset G$ which is control admissible.
A reader may well ask: How difficult is it to construct
$G$ and $H$, and then to check that $H$ is control admissible? Finally, how does then one show that $\CV/H$ is \sfl\ without computing this quotient explicitly?


One of Lie's most significant contributions is that one does not need to work with the symmetries themselves but rather with the {\it infinitesimal symmetries}. These are vector fields on the ambient manifold whose flows are symmetries. The infinitesimal symmetries are obtained from solutions of linear, homogeneous partial differential equations called the {\it determining equations}. These are very often easy to solve. Indeed, as a practical matter the general solution of the determining equations can often be constructed using software such as {\sc Maple} on a laptop computer. This leads to the Lie algebra $\bsy{\G}$ of the control symmetry transformation group $G$, consisting of vector fields whose flows are symmetry transformations of the control system. However, the flows themselves are not required for the analysis of control systems. Only the Lie algebra $\bsy{\G}$ is needed.\looseness=-1

Once a basis $\{X_1,\dots,X_r\}$ for $\bsy{\G}$ has been constructed, it is very easy to write out the structure equations of $\bsy{\G}$ by computing the Lie brackets $[X_i,X_j]$. From these it is often possible to spot Lie subalgebras $\bsy{\Theta}\subset\bsy{\G}$, corresponding to Lie subgroups $H \subset G$, just by inspection.
This approach, as a practical matter, is sufficient to find \sfl\ quotients
and most of our calculations so far have been done this way. This arises from the observation that the most significant property of a subalgebra leading to a \sfl\ quotient is its dimension. Its structure equations or the action it generates do not seem to play a significant role.

However, one can be much more systematic if one is interested in getting as much information as possible from a given symmetry group. A Lie algebra usually has infinitely many subalgebras but they will not all be ``different''. A common way to classify Lie subalgebras is to use the differential of the conjugation map $C\colon G\to G$, defined by $C_g(h)=ghg^{-1}$
for all~${h\in G}$, and some $g\in G$. Two subgroups $H_1$, $H_2$ will be conjugate if there is a $g\in G$ such that~${H_2=C_g(H_1)=gH_1g^{-1}}$. This leads to a classification of subalgebras by the differential of~$C$ and is a well-known, widely used method for producing a finite list of distinct subalgebras and hence distinct Lie subgroups (see \cite[Chapter~3]{OlverLieBook}). One subgroup $H$ from each class on the list, with
Lie algebra $\bsy{\Theta}$, may then be used to check whether
$H$ is control admissible by checking that~${\CV^{(1)}\cap\bsy{\Theta}=0}$ (see Definition \ref{controlAdmissibleDefn}). If so, then to determine if
$\CV/H$
is \sfl\ one checks whether $\CV\oplus\bsy{\Theta}$ is a \stf\ relative Goursat bundle (see Section~\ref{relativeGoursatSect}). There will only be a finite number of subalgebras to check. Each of these can then be checked for the existence of a dynamic feedback linearization, as outlined in Section \ref{dynamicExtensionsSection}.
Essentially, all checking reduces to linear algebra because we are always working ``infinitesimally'' and the number of choices is finite.
In all such calculations, the burden is very significantly eased by the use of dedicated software which often makes them very easy. For this paper, we used the {\sc Maple} package {\sc DifferentialGeometry}~\cite{AndersonTorre}.

Finally, checking contact curve reductions is straightforward. For instance, in the case of~2-input systems there are two possible reductions, one for each input. In the case of 3-input systems there are 3 single input reductions and 3 reductions in pairs, etc. Each reduction only involves linear algebra and differentiation and is therefore algorithmic.

\subsection{Relation between cascade linearization \& differential flatness}

Cascade linearization is an approach to finding dynamic feedback linearizations of a control system by exploiting the existence of its Lie group of symmetries. A flat control system (see Definition~\ref{flatness-def}) is \dfl. Conversely, according to \cite{ChetverikovDFL, LevineDFL}, a \dfl\ control system is flat. Granting this, one could also characterize
cascade linearization as a method for finding the flat outputs of a flat control system with symmetry. It provides a~sufficient condition for the existence of a complete set of {\it flat outputs} (see Definition~\ref{flatness-def}) and it can be checked algorithmically once the Lie algebra of infinitesimal symmetries has been found. The latter problem is discussed above in Section \ref{implementingCFLsect}. The problem of construction then relies on the Frobenius theorem. However, in developing the theory of cascade linearization, we make crucial use of Definitions \ref{dynamic-feedback-def} and \ref{DFLdefn1} of \df\ linearization because that is the framework that most appropriately serves the goals of this paper, as will be shown in Section~\ref{DFLbySymmetrySection}.

\section[Framework for dynamic feedback linearization of invariant control systems]{Framework for dynamic feedback linearization\\ of invariant control systems}\label{CommutativitySection}

In Example \ref{cflExample}, it was possible to find an explicit solution for the system \eqref{CharletExample} after
finding those of $\bsy{\o}/G$ and $\bsy{\bar{\g}}^G$ by applying {\it contact} to each of them. While the system
 \eqref{CharletExample} is sufficiently simple that this application of cascade feedback linearization succeeds, for more complicated examples, finding the explicit solution to the SFL partial contact curve reduction $\bsy{\bar{\g}}^G$
 often proves to be impractical due to the complexity of the computations arising from the arbitrary functions of time present in the coefficients of $\bsy{\bar{\g}}^G$. In the following sections, it will be shown how the contact sub-connection can be exploited to overcome this problem. Even when the solution of $\bsy{\bar{\g}}^G$ cannot be computed explicitly, it is straightforward to determine an upper bound for the {\em signature} of its solution. This upper bound then determines a specific partial prolongation of the original contact sub-connection $\bsy{\g}^G$ which is guaranteed to be static feedback linearizable. In practice, the explicit solution to this prolonged sub-connection is far more straightforward to compute than that for the reduced sub-connection, and it realizes an explicit dynamic feedback linearization for the original system.

This framework will be developed in Section \ref{dynamicExtensionsSection} to establish verifiable criteria for the existence and construction of dynamic feedback linearizations as well as explicit solutions. Subsequently, in running Examples \ref{pvtol sub-connection ex}, \ref{CharletExCont'd}, \ref{pvtol dfl 1}, and in
Sections \ref{tVTOL Ex} and \ref{PVTOLsection}, all features of the general theory will be illustrated.

\subsection{Signature of an explicit solution}\label{solutionSignature}

Recall that to an explicit solution $s\colon\B R \to M$ of $(M, \bsy{\o})$, we can associate the notion of a~{\it signature},
$
\nu=\langle \r_1, \r_2, \dots, \r_k\rangle$,
where $\r_j$ is the number of arbitrary functions occurring to highest order $j$ in the explicit solution~$s$ (cf.\ Definition \ref{EI solSignature}).

We use precisely the same notation to denote the signature of the jet space $J^\nu$, as explained in
Section \ref{cascadeIntegrationSect.}. The dimension of $J^\nu$ and rank of $\bsy{\b}^\nu$ are
\[
N_\nu=\dim J^\nu=1+\sum_{j=1}^k(1+j)\r_j, \qquad \operatorname{rank} \bsy{\b}^\nu=\sum_{j=1}^k j\r_j.
\]
Now a given explicit solution $s$ of signature $\nu$ naturally factors through a map $\psi\colon J^\nu\to M$. That is,
$
s=\psi\circ j^\nu f,
$
where $j^\nu f\colon\B R\to J^\nu$ is the $\nu$-jet of an arbitrary function $f\colon\B R\to \B R^m$ and~${m=\sum_{j=1}^k\r_j}$. Thus we have
\[
0=s^*\bsy{\o}=(j^\nu f)^*\psi^*\bsy{\o},\qquad \forall f.
\]
It follows that the elements of $\psi^*\bsy{\o}$ are contact forms since they are annihilated by the $\nu$-jet of an arbitrary function, and hence
$\psi^*\bsy{\o}\subseteq \bsy{\b}^\nu$.

In the special case when $\dim M=\dim J^\nu$, $\operatorname{rank} \bsy{\o}=\operatorname{rank} \bsy{\b}^\nu$, we prove that $\psi$ is a local diffeomorphism and we conclude that $\psi^*\bsy{\o}=\bsy{\b}^\nu$.
Indeed, we have
\begin{Proposition}\label{equivalence2brunovsky}
Let $s\colon\B R\to M$ be the explicit solution of a control system $(M,\bsy{\o})$ of signature~${\nu=\langle\r_1, \r_2,\dots,\r_k\rangle}$, and suppose that $\dim M=N_\nu$ and
$\operatorname{rank} \bsy{\o}=\operatorname{rank} \bsy{\b}^\nu$. Let $\psi\colon J^\nu\to M$ be the smooth map that locally factors $s$ as
$s=\psi\circ j^\nu f$, where $f\colon\B R\to\B R^m$ is an arbitrary smooth function and
$m=\sum_{i=1}^k\r_i$. Then $\psi$ is a local diffeomorphism and $\psi^*\bsy{\b}^\nu = \bsy{\o}$. Furthermore, $\psi$~is a \stf\ transformation.
\end{Proposition}
\begin{proof}
Suppose the derivative map of $\psi$ is singular in some open set ${\rm U}\subseteq J^\nu$, and that
${s(t)=(\psi_{|_{\rm U}}\circ j^\nu f)(t)}$.
If $\psi$ has components
$\psi^i$, $1\leq i\leq N_\nu:=n$, then ${\rm d}\psi^1\wedge {\rm d}\psi^2\wedge\cdots\wedge {\rm d}\psi^n\equiv 0$ on ${\rm U}$.
Hence there is a regular function $F$ at $x\in\B R^n$ (i.e., ${\rm d}F$ is non-zero at $x$) such that
$F\bigl(\psi^1,\dots,\psi^n\bigr)\equiv 0$ in a neighborhood ${\rm N}_x$ of $x$.
By the regularity of $F$, $F(\psi)\circ j^\nu f=0$
can be expressed as a locally solvable ordinary differential equation for $f$, if necessary by shrinking~${\rm N}_x$. This contradicts the hypothesis that $f$ is arbitrary. Hence, by the inverse function theorem, we deduce that $\psi$ is a local diffeomorphism. Since time $t$ is a parameter along trajectories, it follows from
\cite[Theorem~3.11]{DTVa} that $\psi$ is a \stf\ transformation.
\end{proof}

\begin{Proposition}\label{EqualityOfSignatures}
Let $(M, \bsy{\o})$ be an explicitly integrable control system on a manifold $M$ whose explicit solution has
signature $\u$. Then the explicit solution determines a static feedback linearization of $(M, \bsy{\o})$ if and only if $\dim M=\dim J^\u$. Consequently, $(M, \bsy{\o})$ is SFL if and only it admits an explicit solution whose signature $\u$ satisfies $\dim M=\dim J^\u$.
\end{Proposition}
\begin{proof}
Since the explicit solution $s(t)$ has signature $\u$, we can express $s$
by a formula of the form~${\bsy{x}(t)=H(j^\u f(t))}$, where $\bsy{x}$ denotes all the dynamical variables of the control system~${(M, \bsy{\o})}$. If $\dim M=\dim J^\u$ then we can invert $\bsy{x}=H(z^\u)$, defining a local diffeomorphism $\th\colon M\to J^k$ such that $z^\u=\th(\bsy{x})$ and satisfying $\th^*\bsy{\b}^\u=\bsy{\o}$. For if the map $z^\u\mapsto H(z^\u)$ is not full rank, it implies dependencies among the variables $\bsy{x}$. By Proposition \ref{equivalence2brunovsky}, $\bsy{\o}$ is \sfl.

Conversely, if $(M, \bsy{\o})$ is \sfl\ then it has some explicit solution with the signature $\u'=\langle\r_1,\r_2,\dots,\r_k\rangle$ of some
Brunovsk\'y normal form $\bsy{\b}^{\u'}$,
and therefore $\bsy{\o}$ has an explicit solution of signature $\u'$. But
the explicit solution of $\bsy{\o}$ has signature $\u$ and hence $\u'=\u$. Thus $\dim M=1+\sum_{j=1}^k(1+j)\r_j=\dim J^\u$.
\end{proof}

We now study the general case $\dim M<N_\u$.

\begin{Definition}[partial prolongation of a control system]\label{pp-def}
Let $(M, \bsy{\o})$ be a DFL control system. Suppose $s\colon\B R\to M$ is an explicit solution of $\bsy{\o}$; let $\u$ denote the signature of $s$ and assume that $\dim M<N_\u$. Consider the augmented control system $\bsy{\o'}:=\widetilde{\pi}^*\bsy{\o}\oplus\bsy{\alpha}$,
where
\[
\bsy{\alpha}=\operatorname{span}\{{\rm d}u^a-p_1^a {\rm d}t,\, {\rm d}p_1^a-p_2^a {\rm d}t,\,\dots,\,
{\rm d}p_{k_a-1}^a-p_{k_a}^a {\rm d}t\},\qquad a\in \{1, 2,\dots,m\},
\]
on a manifold ${M\bsy{'}}$, where $\wt{\pi}\colon M'\to M$ is the submersion. If the explicit solution $s$ can be augmented to an explicit solution $s'$ of $\bsy{\o}'$ without changing the signature of $s$ (i.e., the signatures of $s$ and $s'$ are both equal to $\u$), then we call $(M\bsy{'}, \bsy{\o'})$ the {\em $\u$-prolongation} of $(M, \bsy{\o})$ and $s'$ the~$\u$-prolonged explicit solution.
\end{Definition}

\begin{Example}
Let $\bsy{\o}$ be a rank 5 control system on the manifold $M$ with states and inputs~${x^1,\dots, x^5}$, $u^1$, $u^2$.
Suppose $\bsy{\o}$ has explicit solution $s$ given by
\begin{align*}
&x^1=\dot{g}f^2/\dot{f},\qquad x^2=\bigl(\dot{f}g+f\dot{g}\bigr)/\dot{f},\qquad x^3=f,\qquad u^1=\dot{f},\cr
&x^4=\bigl(\ddot{f}g+f\ddot{g}\bigr)/\dot{f}^3=:F,\qquad x^5={\rm d}F/{\rm d}t,\qquad u^2={\rm d}^2F/{\rm d}t^2,
\end{align*}
in terms of arbitrary functions $f$, $g$. The signature of $s$ is $\u=\langle 0,0,0,2\rangle$ and we have
$8=\dim M<\dim J^\u=11$. We construct a $\u$-prolongation of $s$ by adjoining the equations
\[
p_1^1=\ddot{f},\qquad p^1_2=\dddot{f},\qquad p^1_3=\ddddot{f}.
\]
The resulting augmented solution $s'$ also has signature $\u = \langle 0,0,0,2\rangle$ and factors through a~\stf\ transformation $\psi\colon J^\u\to M{\bsy{'}}$, where
$M{\bsy{'}}$ has local coordinates
$t$, $x^1,\dots, x^5$, $u^1$, $u^2$, $p_1^1$, $p^1_2$, $p^1_3$,
and the $\u$-prolongation of $\bsy{\o}$ is
\[
\bsy{\o'}=\widetilde{\pi}^*\bsy{\o}\oplus\bigl\{{\rm d}u^1-p^1_1 {\rm d}t,\, {\rm d}p^1_1-p^1_2 {\rm d}t,\, {\rm d}p_2^1-p^1_3 {\rm d}t\bigr\}.
\]
Thereby we have carried out a {\em partial prolongation} by differentiating $u^1$ three times, while $u^2$ is left undifferentiated. We sometimes call this a 3-fold ``{\em partial prolongation along $u^1$}''. Then we have
$\operatorname{rank} \bsy{\o'}=8=\operatorname{rank} \bsy{\b}^{\langle 0,0,0,2\rangle}$. Hence $\psi^*\bsy{\o'}=\bsy{\b}^{\langle 0,0,0,2\rangle}$ by Proposition~\ref{equivalence2brunovsky}, where $\psi$ is a~\stf\ transformation.
\end{Example}

\begin{Example}\label{pvtol sub-connection ex}
We work through an additional example of how the principal bundle isomorphism~${\bigl(\varphi, \wt{\vf}\bigr)}$ of Figure \ref{diagramBasic} and the contact sub-connection
$\CH_G$, defined in \eqref{contactConnection}, can be constructed in practice by studying the well-known PVTOL control system \cite{VTOL,DynamicVTOL} given by
\begin{align}
&\dot{x}=x_1,\qquad
\dot{x}_1=-u^1\sin\th+h u^2\cos\th,\qquad
\dot{z}=z_1,\nonumber\\
&\dot{z}_1=u^1\cos\th+h u^2\sin\th-1,\qquad
\dot{\th}=\th_1,\qquad
\dot{\th}_1=u^2,\label{tVTOL1}
\end{align}
in 6 states and 2 controls, where $h$ is a parameter. We denote the ambient manifold of the system by $M$. The Lie group of all control symmetries of \eqref{tVTOL1} is studied in detail in forthcoming work. One of its subgroups consists of Galilean transformations in $x$ and $z$, generated infinitesimally by the abelian Lie algebra
\begin{equation*}
\bsy{\G'}=\operatorname{span} \{t\P x+\P {x_1},\, \P x,\, t\P z+\P {z_1},\, \P z \}.
\end{equation*}
It turns out that $\bsy{\G'}$ is not strongly transverse to $\CV$ and we instead use the subalgebra
\begin{equation*}
\bsy{\G}=\operatorname{span} \{t\P x+\P {x_1},\P x, \P z \}\subset\bsy{\G'}
\end{equation*}
to study \eqref{tVTOL1}. The subalgebra $\bsy{\G}$ generates Galilean transformations in $x$ and translations in $z$ and we denote this group by~$G$.

We calculate that $\wh{\CV}:=\CV\oplus\bsy{\G}$ is a static feedback relative Goursat bundle of
signature $\langle 1, 1\rangle$. This means $\CV/G$ is \stf\ equivalent to $\CB_{\langle 1, 1\rangle}$.
The $G$-invariant functions are the first integrals of $\bsy{\Gamma}$ and are generated by
\begin{align*}
&t=t,\qquad q_1=\th,\qquad q_2=\th_1,\qquad q_3=z_1,\qquad v^1=u^1,\qquad v^2=u^2,
\end{align*}
and thus $t$, $q^1$, $q^2$, $q^3$, $v^1$, $v^2$ form local coordinates on $M/G$. The quotient control system is then easily computed and given by
\begin{equation*}
\CV/G=\operatorname{span}\bigl\{\P t+q^2\P {q^1}+v^2\P {q^2}+\bigl(v^1\cos q^1+hv^2\sin q^1-1\bigr)\P {q^3},
\, \P {v^1},\, \P {v^2}\bigr\}.
\end{equation*}
Applying the procedure {\it contact} to $\CV/G$, we obtain the transformation
$\varphi\colon M/G \to J^{\langle 1, 1\rangle}$ given~by
\begin{align*}
&t=t,\qquad m=q^1,\qquad m_1=q^2,\qquad m_2=v^2,\qquad n=q^3,\\ &n_1=u^1\cos q^1+hv^2\sin q^1-1,
\end{align*}
where $t$, $m$, $m_1$, $m_2$, $n$, $n_1$ denote the standard contact coordinates on $J^{\langle 1, 1\rangle}$. Then
\begin{align*}
&\varphi_*(\CV/G)=
\operatorname{span}\{\P t+m_1\P m+m_2\P {m_1}+n_1\P n, \,\P {m_2}, \,\P {n_1}\}=\CB_{\langle 1, 1\rangle},
\end{align*}
as claimed. To construct $\wt{\vf}$, we need the coordinates of the local group $G$ in terms of $t$ and the state variables of $\mcal{V}$. As per the discussion following equation \eqref{bundle sfl}, we find that
\begin{equation*}
\e^1=x_1,\qquad \e^2=x-t x_1,\qquad \e^3=z.
\end{equation*}
Then as in Figure \ref{diagramBasic}, we have the \stf\ transformation
\[
\wt{\varphi}=\varphi\circ\pi\times \bigl(\e^1=x_1,\, \e^2=x-t x_1,\, \e^3=z\bigr)
\]
giving,
\begin{gather}
\wt{\varphi}_*\CV= \operatorname{span}\biggl\{\P t+m_1\P m+m_2\P {m_1}+n_1\P n+
\frac{(n_1+1)\sin m-hm_2}{\cos m}(t\P {\e^2}-\P {\e^1})+n\P {\e^3}, \nonumber\\
\hphantom{\wt{\varphi}_*\CV= \operatorname{span}\biggl\{}{} \P {m_2}, \,\P {n_1}\biggr\}.\label{connectionVTOLgalilean}
\end{gather}
The distribution $\CH_G:=\wt{\varphi}_*\CV$ is the contact sub-connection for the control system~\eqref{tVTOL1},
with respect to the control admissible symmetry group~$G$.
For completeness and later use, we also record the annihilator $\bsy{\g}^G$ of $\CH_G$,
\begin{align}
\bsy{\g}^G={}&\{{\rm d}m-m_1{\rm d}t,\, {\rm d}m_1-m_2{\rm d}t,\, {\rm d}n-n_1{\rm d}t\}\nonumber\\
&\oplus\!\left\{\!{\rm d}\e^1\!+\frac{(n_1+1)\sin m-hm_2}{\cos m}{\rm d}t, \,{\rm d}\e^2\!-t\frac{(n_1+1)\sin m-hm_2}{\cos m}{\rm d}t,\,
{\rm d}\e^3\!-n{\rm d}t\!\right\}.\!\label{connectionFormVTOLgalilean}
\end{align}
\end{Example}

\section{Dynamic linearization via symmetry reduction}\label{dynamicExtensionsSection}\label{DFLbySymmetrySection}

We saw in Example \ref{cflExample} that to construct an explicit solution of a \cfl\ control system it is sufficient to construct the integral submanifolds of the static feedback linearizable reduced sub-connection
$\bar{\CH}_G$. We called this procedure cascade feedback linearization. While this direct approach is straightforward {\em in principle}, the presence of arbitrary functions of time in
 $\bar{\CH}_G$ usually makes the construction of an explicit solution, as a practical matter, computationally challenging due to expression swell. One of the purposes of this section is to prove (in Theorem \ref{prolongationPredictor}) that, in fact, for the purpose of constructing an explicit solution the integration of $\bar{\CH}_G$ can be completely avoided. The only information we shall require of $\bar{\CH}_G$ is its {\em derived length}.

To construct an explicit solution we require a dynamic feedback linearization.
We will establish a procedure for constructing \df\ linearizations for control systems with symmetry in accordance with the standard definition (see Definition \ref{dynamic-feedback-def}).

\subsection{Explicit solution via the contact sub-connection}
In Example \ref{cflExample}, we produced a formula for the explicit trajectories of the control system \eqref{CharletExample}. We did this by computing a {reduced} contact sub-connection $\bar{\CH}_G$ which is SFL and led to the explicit trajectories of $\CH_G$ as a control system in its own right. In this section, we apply the main results of Section \ref{solutionSignature} to $\CH_G$, showing that the linearization of $\bar{\CH}_G$ can be avoided by proving that one can always achieve a dynamic feedback linearization of $\CH_G$ via an integrator chain (partial prolongation) when $\bar{\mathcal{H}}_G$ is SFL.

The information required to determine such a DFL of $\CH_G$ via partial prolongation need not come from the contact coordinates of $\bar{\CH}_G$. These may be computationally unwieldy (though are manageable in Example \ref{cflExample}). Indeed, all we need is the signature of the explicit solution of~${\CH}_G$, which can be determined by carefully inspecting the highest order derivatives of the family of partial contact curves featured in its fundamental functions, together with the derived length of~$\bar{\CH}_G$, without being required to explicitly compute the remaining contact coordinates. But we will further show that, in practice, even the fundamental functions need not be computed.

\begin{Example}\label{CharletExCont'd}
Recall that the partial contact curve reduction of $\CH_G$ for the control system~\eqref{CharletExample} in Example \ref{cflExample} is given by
\[
\bar{\CH}_G=\operatorname{span}\bigl\{\partial_t+z_1\partial_z+z\bigl(1-\ddot{f}(t)\bigr)\partial_\e,\partial_{z_1} \bigr\}
\simeq_{\text{SFE}}\CB_{\langle 0,1 \rangle},
\]
for arbitrary, smooth, real-valued functions $f$.
As we saw, the fundamental function of $\bar{\CH}_G$ turned out to be $\e$. Letting $t$, $a$, $a_1$, $a_2$ be the contact coordinates for the Brunovsky normal form of $\bar{\CH}_G$, we let $a=\e$ and the remaining contact coordinates are determined from successive Lie differentiation of the fundamental function by the vector field
\[
\bar{Z}=\partial_t+z_1\partial_z+z\bigl(1-\ddot{f}(t)\bigr)\partial_\e.
\]
Without carrying out the computations to determine the remaining contact coordinates \smash{$a_1\!=\!\LieD_{\bar{Z}} a$} and \smash{$a_2=\LieD^2_{\bar{Z}} a$}, we notice that the highest order derivative of $f(t)$ in the family of partial contact curves that appears in $\bar{Z}$ is two. As such, $a_1$ must depend on \smash{$\ddot{f}(t)$} and $a_2$ must depend on \smash{$\dddot{f}(t)$}. Additionally, since the signature of the Brunovsk\'y form equivalent to $\bar{\CH}_G$ is $\langle 0,1 \rangle$, then we can deduce that the signature of an explicit solution $s$ of $\CH_G$ must be
$\nu=\langle 0, 1, 1 \rangle$. The dimension of the associated jet space is
$N_\nu=8>7=\dim \bigl(J^{\langle 1,1 \rangle} \times G\bigr)$. As such, we confirm again (this time via Proposition \ref{EqualityOfSignatures}) that the system \eqref{CharletExample} is not SFL; however, if we perform a prolongation of~$\bsy{\g}^G$ along the $w$ jet coordinate by adjoining the 1-form ${\rm d}w_2-w_3 {\rm d}t$ to $\bsy{\g}^G$
(see Definition \ref{pp-def})
in equation \eqref{connectionFormCharlet}, we get
\smash{$\bsy{\widecheck{\bsy{\g}}}^G:=\bsy{\g}^G\oplus\{{\rm d}w_2-w_3 {\rm d}t\}$}, which turns out to be a DFL of $\bsy{\g}^G$,
since~\smash{$\bsy{\widecheck{\bsy{\g}}}^G$} is SFL by Proposition \ref{EqualityOfSignatures}.
\end{Example}
This example is an instance of a general result concerning prolongation of the jet coordinates involved in the partial contact curve reduction. In fact, the order counting technique demonstrated above can be generalized to any contact sub-connection and is the subject of Theorem~\ref{prolongationPredictor} below. To prove this, we first introduce some helpful notation and establish a technical lemma concerning the fundamental functions of a reduced contact sub-connection.

\begin{notate}
Given a signature $\u=\langle \r_1,\dots,\r_k \rangle$,
we will use subscripts to denote an increase in order of all jets so that
\begin{equation*}
\u_\ell=\langle \underbrace{0, \dots, 0}_{\ell\text{ entries }},\r_1,\dots,\r_k \rangle,\qquad \ell\geq0,
\end{equation*}
and $\u_\ell=\u$ if $\ell<0$.
\end{notate}
\begin{notate}\label{max sig sum}
Consider a collection of signatures $S=\{\nu_\ell\}_{\ell=1}^k$ and their corresponding jet spaces $J^{\nu_\ell}(\mathbb{R},\mathbb{R}^{m_{\nu_\ell}})$, where $m_{\nu_\ell}=m_\nu$ for all $1\leq \ell \leq k$ and $m_\nu$ is the sum of the entries of $\nu$.
We denote by $\nu^*_{S}$ the signature of the smallest jet space $J^{\nu^*_S}(\mathbb{R},\mathbb{R}^{m_\u})$ for which the collection of jet spaces $\{J^{\nu_\ell}(\mathbb{R},\mathbb{R}^{m_{\u_\ell}})\}_{\ell=1}^k$ may be embedded into $J^{\u^*_{S}}(\mathbb{R},\mathbb{R}^{m_{\u}})$.
\end{notate}

For example, if $S=\{\u_0,\u_1\}$ where $\u_0=\langle 1,0,1,1\rangle$ and $\u_1=\langle 0,3 \rangle$ then $\u^*_S=\langle 0,1,1,1 \rangle$ and~${J^{\u_\ell}\subset J^{\u^*_S}}$, $\ell=0,1$; $J^{\u^*_S}$ is the smallest such jet space.
\begin{Definition}\label{partialOrderOnSignatures}
If $J^{\u_1}$ and $J^{\u_2}$ are jet spaces, we say that $\u_1$ {\it is less than or equal to} $\u_2$ if~${J^{\u_1}\subseteq J^{\u_2}}$ and we write $\u_1\leq \u_2$.
\end{Definition}

\begin{Lemma}\label{fund func sig}
Let $\CH_G=\ker\bigl(\bsy{\b}^{\u^\perp}\oplus\bsy{\b}^\u\oplus\Theta^G\bigr)$ be the contact sub-connection of a $G$-invariant control system $\bsy{\omega}$, with $G$ being control admissible. Assume further that $\bar{\CH}_G$ is an SFL contact curve reduction of $\CH_G$, where the family of partial contact curves $\bsy{c}^\nu_f$ are induced by integral manifolds of $\bsy{\b}^\nu$ with $\nu=\langle \r_1,\dots,\r_{k_\u}\rangle$ and \smash{$\bsy{\b}^\k=\bsy{\b}^{\u^\perp}\oplus\bsy{\b}^\u$}. Each order $\ell$ fundamental function~${\phi^{a_\ell,\ell}}$,~${1\leq a_\ell\leq \rho_\ell}$, $1\leq \ell\leq \bar{k}$ of $\bar{\CH}_G$ defining its \stf\ linearization depends on derivatives of $\bsy{f}(t)$ having signature no larger than $\u_{\ell-1}$, where $\u_{\ell-1}$ is the signature of the order~${\ell-1}$ time derivative of
$j^\u \bsy{f}(t)$ and
$\bar{k}$ is the derived length of $\bar{\CH}_G$.
\end{Lemma}
\begin{proof}
Let $\bar{Z}\in \bar{\CH}_G$ denote the vector field
\[
\bar{Z}=\mathbf{\bar{D}}_t+\sum_{a=1}^r\lambda^a(\bsy{z}, j^{\nu} \bsy{f}(t))R_a,
\]
where the $\{ R_a\} _{a=1}^r$ form a basis of right invariant vector fields on $G$ and
$\mathbf{\bar{D}}_t$ is the total derivative operator on \smash{$J^{\u^\perp}$}. Notice that, depending on the functions $\lambda^a$, the vector field $\bar{Z}$ may depend on at most as many derivatives of $\bsy{f}(t)$ as are defined by the Brunovsk\'y normal form $\bsy{\b}^\u$. By the procedure {\it contact}, there exist fundamental functions $\phi^{a_\ell,\ell}$ for each $\ell<\bar{k}-1$ such that the quotient bundle \smash{$\ch\bar{\CH}_G^{(\ell)}/\ch\bigl(\bar{\CH}_G\bigr)^{(\ell)}_{\ell-1}$} is nontrivial. In particular, the fundamental functions are the invariants of the integrable quotient bundle \smash{$\bar{\Xi}^{(\ell)}_{\ell-1}/\bar{\Xi}^{(\ell)}$} where
\[
\bar{\Xi}^{(\ell)}=\operatorname{ann} \bigl(\ch\bar{\CH}_G^{(\ell)}\bigr),\qquad \bar{\Xi}^{(\ell)}_{\ell-1}=\operatorname{ann} \bigl(\ch\bigl(\bar{\CH}_G\bigr)^{(\ell)}_{\ell-1}\bigr).
\]
 Each Cauchy bundle
\smash{$\ch\bar{\CH}_G^{(\ell)}$} satisfies
\begin{equation*}
\ch\bigl(\bar{\CH}_G\bigr)^{(\ell)}_{\ell-1}\subseteq \ch\bar{\CH}_G^{(\ell)}\subseteq
\ch\bigl(\bar{\CH}_G\bigr)^{(\ell+1)}_{\ell},
\end{equation*}
by definition of the intersection bundle; see \eqref{intersectionBundleDefn}. Furthermore, as per equation \eqref{ad char},
\begin{align*}
\begin{aligned}
&\ch\bigl(\bar{\CH}_G\bigr)^{(\ell)}_{\ell-1}=\operatorname{span}\bigl\{\bar{C}_0, \operatorname{ad}(Z)\bar{C}_0,\dots,\operatorname{ad}(Z)^{\ell-1}\bar{C}_0\bigr\},\\
&\ch\bigl(\bar{\CH}_G\bigr)^{(\ell+1)}_{\ell}=\operatorname{span}\bigl\{\bar{C}_0, \operatorname{ad}(Z)\bar{C}_0,\dots,\operatorname{ad}(Z)^{\ell}\bar{C}_0\bigr\},\qquad
\bar{C}_0:=\ch\bigl(\bar{\CH}_G\bigr)^{(1)}_0.
\end{aligned}
\end{align*}
In particular, we notice that for each $1\leq i \leq \ell$ the coefficients of the vector fields in $\operatorname{ad}^i(Z)\bar{C}_0$ will have at most $i-1$ time derivatives of $j^\u \bsy{f}(t)$. As such, with $i=\ell$, the first integrals $\phi^{a_\ell,\ell}$ of the non-trivial quotient bundles \smash{$\bar{\Xi}^{(\ell)}_{\ell-1}/\bar{\Xi}^{(\ell)}$} have at most dependence on $\ell-1$ time derivatives of $j^\u \bsy{f}(t)$, which we denote by $j^{\u_{\ell-1}} \bsy{f}(t)$. Finally, in case $\bar{\Delta}_{\bar{k}}=1$ the associated fundamental bundle \smash{$\bar{\Pi}^{\bar{k}}$} as given by equation \eqref{ad char} allows us to make the same conclusion when $\ell=\bar{k}$. In case $\bar{\Delta}_{\bar{k}}>1$, a similar argument applies.
\end{proof}

\begin{Theorem}[explicit solution via symmetry]\label{prolongationPredictor}
Suppose the control system $(M,\bsy{\o})$ of derived length $k$ is cascade feedback linearizable with respect to a Lie group~$G$.
Let $\CH_G$ be its contact sub-con\-nection on $J^\kappa\times G$, where $J^\kappa$ is decomposed as $\bigl(J^\u\times J^{\u^\perp}\bigr)/{\sim}$ and $\bar{\CH}_G$ is the static feedback linearizable reduced contact sub-connection with respect to a generic family of partial contact curves $\bsy{c}^\u_f$. Finally, let $\operatorname{pr}\CH_G$ be the partial prolongation of $\CH_G$ along its $J^\u$ components with the resulting signature
of $\operatorname{pr}\CH_G$ given by \smash{$\kappa'=\u^\perp+\u'$} where $\u'=\u^*_S$,
with \[S=\bigl\{\u^{a_\ell,\ell}_\ell \colon a_\ell\neq0,\, 1\leq \ell \leq\bar{k}\bigr\}\] where $\u^{a_\ell,\ell}$ denotes the signature of the derivatives of $\bsy{f}(t)$ appearing in each fundamental function $\phi^{a_{\ell},\ell}$. Then $\operatorname{pr}\CH_G$ is SFL.
\end{Theorem}
\begin{proof}
 Let us decompose the contact coordinates on $J^\k= \bigl(J^{\u^\perp}\times J^\u\bigr)/{\sim}$ such that $J^{\u^\perp}$ has coordinates $(t,\bsy{z})$ and $J^\u$ has coordinates $(t,\bsy{w})$. Once we restrict $\CH_G$ to a family of partial contact curves $\bsy{c}^\u_{\bsy{f}}$ we will have
\begin{equation}\label{reducedHGgeneral}
\bar{\CH}_G=\operatorname{span}\biggl\{\mathbf{\bar{D}}_t+\sum_{a=1}^r\lambda^a(\bsy{z}, j^{\nu} \bsy{f}(t))R_a,\, \P {z^{a_j}_{j}}\biggr\},
\end{equation}
where $z^{a_j}_j$ are the highest order contact coordinates on $J^{\u^\perp}$.
By hypothesis, $\bar{\CH}_G$ is static feedback linearizable of some signature $\bar{\k}$. Consequently, the static feedback transformation that implements the linearization is determined by the procedure {\it contact}~\cite{VassiliouGoursatEfficient}. According to this, one computes the fundamental functions $\phi^{a_\ell, \ell}$ of all orders, after which the contact coordinates are determined by Lie differentiation by the total differential operator $\bar{Z}$, which in this case is given by the first vector field in~\eqref{reducedHGgeneral}.\footnote{This is because we assume that $\bar{\CH}_G$ has passed the test of being static feedback linearizable.}
Each fundamental function $\phi^{a_\ell,\ell}$ depends upon the contact coordinates
$\bsy{z}$, the local group coordinates $\bsy{\epsilon}$ and the jet \smash{$ j^{\u^{a_\ell,\ell}} \bsy{f}(t)$} of the family of functions~$\bsy{f}(t)$ that define the contact curve reduction. Let $\bsy{\bar{z}}$ denote the contact coordinates for the static feedback linearization of $\bar{\CH}_G$. Then,
\begin{equation*}
\bsy{\bar{z}}^{a_\ell, \ell}_{i_\ell}=\LieD^{i_\ell}_{\bar{Z}}\phi^{a_\ell, \ell}\bigl(\bsy{z}, j^{\u^{a_{\ell},\ell}} \bsy{f}(t), \bsy{\e}\bigr),\qquad
0\leq i_\ell\leq \ell,\quad 1\leq \ell \leq \bar{k}.
\end{equation*}
As such, each highest order contact coordinate $\bsy{\bar{z}}^{a_\ell,\ell}_\ell$ of the linearizing static feedback transformation depends upon, at most, the \smash{$\u^{a_\ell,\ell}_\ell$} jet of $\bsy{f}(t)$. That is,
\begin{equation}\label{HbarContactCoords}
\phi^{a_\ell,\ell}_\ell=\phi^{a_\ell,\ell}_\ell\bigl(\bsy{z},j^{\u^{a_\ell,\ell}_{\ell}}\bsy{f}(t),\bsy{\e}\bigr),\qquad 0\leq l \leq \bar{k}.
\end{equation}
Let $\bar{\k}$ be the signature of $\bar{\CH}_G$ as a SFL control system, that is, $\bar{\CH}_G\cong\operatorname{ann} \bsy{\b}^{\bar{\k}}$. Note the list of non-negative integers $\bar{\k}$ does not contain information about the jet orders of $\bsy{f}(t)$ that occur in the fundamental functions $\phi^{a_\ell,\ell}$. Let
\begin{equation*}
\u'=\u^*_S,\qquad S=\bigl\{\u^{a_\ell,\ell}_\ell \colon a_\ell\neq0,\, 1\leq \ell \leq\bar{k}\bigr\}.
\end{equation*}
Then solving for $\bsy{z}$ and $\e$ in \eqref{HbarContactCoords} obtains the explicit solution of $\bar{\CH}_G$
in the form $(\bsy{z}, \e)=\zeta\bigl(j^{\bar{\k}}\bsy{g}, j^{\nu'} \bsy{f}\bigr)$, for an arbitrary, smooth function
$\bsy{g}(t)$, valued in $\B R^{\bar{m}}$, where $\bar{m}$ is the number of inputs in $\bar{\CH}_G$. Combining this with the partial contact curve reduction in which $\bsy{w}=j^\u {\bsy{f}}$, obtains the explicit solution of $\CH_G$. Thus we see that the signature of the {\it explicit solution} of~$\CH_G$~-- say~$\k'$~-- is given by, $\k'=\bar{\k}+\u'$.
We note in passing that $\k'$ is generally not the signature of~$\CH_G$ itself. Only if $\CH_G$ is \stf\ linearizable will this be so.

To finish the proof, we invoke Proposition~\ref{equivalence2brunovsky} by constructing a $\k'$-prolongation as in Definition \ref{pp-def}. Toward this goal, we need some dimension counts:
\begin{align*}
&\dim  J^\k\times G=r+1+\sum_{\r_i\in \k}(1+i)\r_i=r+1+\sum_{\r_i\in \u^\perp}(1+i)\r_i+\sum_{\r_i\in \u}(1+i)\r_i,\\
&\dim  J^{\u^\perp}\times G=r+1+\sum_{\r_i\in \u^\perp}(1+i)\r_i,\\
&\dim  J^{\k'}=1+\sum_{\r_i\in \bar{\k}}(1+i)\rho_i+\sum_{\r_i\in \u'}(1+i)\rho_i,
\end{align*}
where $\dim  J^{\u^\perp}\times G=\dim  J^{\bar{\k}}$. Let $N_{\u'-\u}:=\dim  J^{\u'}-\dim  J^\u$. Then it follows that
\begin{equation*}
\dim  J^{\k'}=N_{\u'-\u}+\dim  J^\k\times G.
\end{equation*}
Let us then partially prolong the contact sub-connection $\bsy{\g}^G=\bsy{\b}^\kappa\oplus\Theta^G$ so that the $\bsy{\b}^\u$ Brunovsk\'y form summand in $\bsy{\b}^\k$ is prolonged to $\bsy{\b}^{\u'}$. In that case the Pfaffian system
$\bsy{\g}^G$ thus prolonged will be denoted by $\bsy{\widecheck{\bsy{\g}}}^G$ and we let $\operatorname{pr}\CH_G:=\ker\bsy{\widecheck{\bsy{\g}}}^G$. Matters being so, this partial prolongation is exactly the
$\kappa'$-prolongation of Definition \ref{pp-def}, and therefore by Proposition \ref{EqualityOfSignatures} we conclude that~$\operatorname{pr}\CH_G$ is SFL with signature $\k'$.
\end{proof}

We point out that Theorem \ref{prolongationPredictor} can also be expressed by the diagram in Figure \ref{prolongPredict-diagram}.
\begin{figure}[ht]\centering
\begin{tikzcd}[every label/.append style = {font = \normalsize}]
{\bigl(J^{\nu'+\nu^\perp}\times G,\boldsymbol{\beta}^{\nu'}\oplus\boldsymbol{\beta}^{\nu^\perp}\oplus\Theta^G\bigr)} \arrow[r, "\varphi'"] \arrow[d, "\tilde{\pi}'"] & {\bigl(J^{\bar{\kappa}+\nu'},\boldsymbol{\beta}^{\bar{\kappa}+\nu'}\bigr)} \\
{\bigl(J^\kappa\times G,\boldsymbol{\gamma}^G\bigr)=\bigl(J^{\nu+\nu^\perp}\times G,\boldsymbol{\gamma}^G\bigr)} &
\end{tikzcd}
\caption{Here $\varphi'$ is the static feedback map that linearizes the prolonged contact sub-connection.}
\label{prolongPredict-diagram}
\end{figure}

The number of partial prolongations of $\CH_G$ described in Theorem \ref{prolongationPredictor} to achieve the \stf\ linearization $\operatorname{pr}\CH_G$ is exact. However, as far as possible we want to avoid the task of computing the fundamental functions of $\bar{\CH}_G$ explicitly by deducing the {\it maximum number} of prolongations of
$\CH_G$ required to achieve a \stf\ linearization. Such a maximum is provided by Corollary \ref{prolongationPredictor-cor}, below. Before that we need to state the following lemma.
\begin{Lemma}\label{prolongationOFsflConnection}
Let $\bigl(J^{\u^\perp+\u'}\times G, \operatorname{pr}\CH_G\bigr)$ be the minimally prolonged \stf\ linearizable contact connection provided by Theorem {\rm\ref{prolongationPredictor}}. Performing any partial prolongation along any of the components corresponding to the nonzero entries in $\u'$ preserves the \stf\ linearizability of $\operatorname{pr}\CH_G$.
\end{Lemma}
\begin{proof}
This follows easily, for instance by using the Gardner--Shadwick test for \stf\ linearizability \cite{GSalgorithm}.
\end{proof}

\begin{Corollary}\label{prolongationPredictor-cor}
Let $\CH_G$ be as in the hypotheses of Theorem~{\rm\ref{prolongationPredictor}} and let $\bar{k}$ be the derived length of $\bar{\CH}_G$. Then, it is enough to prolong all the $J^\u$ components of $\CH_G$ by $2\bar{k}-1$ differentiations, so as to achieve a SFL control system $\operatorname{pr} \CH_G$.
\end{Corollary}
\begin{proof}
 Recall that to compute an optimal set of contact coordinates for $\operatorname{pr}\CH_G$ we examined the fundamental functions $\phi^{a_\ell, \ell}$ of $\bar{\CH}_G$ which must be differentiated $\ell$-times and then we minimized over the set of signatures
\smash{$S=\bigl\{\u^{a_\ell,\ell}_\ell \bigr\}_{\ell=1}^{\bar{k}}$}.
By Lemma \ref{fund func sig}, the $\phi^{a_\ell, \ell}$ depend upon $j^\l \bsy{f}$, where~${\l\leq \u_{\ell-1}}$ and hence the $\ell^{\text{th}}$ derivatives of $\phi^{a_\ell, \ell}$ depend upon
$j^\mu \bsy{f}$, where $\mu\leq \u_{2\ell-1}$. Replacing
$S$ by $S'=\{\u_{2\ell-1}\}_{\ell=1}^{\bar{k}}$
we have $\u^*_{S'}=\u_{2\bar{k}-1}$. By Lemma \ref{prolongationOFsflConnection}, we conclude that since the partial prolongation of the summand $\b^\u$ in $\bsy{\g}^G$ prolonged to $\b^{\u_S^*}$ leads to an SFL partial prolongation $\bsy{\widecheck{\bsy{\g}}}^G$, the
partial prolongation of $\b^\u$ to $\b^{\u_{2\bar{k}-1}}$ instead, will ensure that this partial prolongation of
$\bsy{\g}^G$ will also be SFL.
\end{proof}

Thus, the only data we require from $\bar{\CH}_G$ to obtain an SFL prolonged contact connection~$\operatorname{pr}\CH_G$ is its derived length, which is algorithmically obtained. The number of prolongations described in Corollary~\ref{prolongationPredictor-cor} agrees with the well-known sharp Sluis--Tilbury
bound \cite{SluisTilbury} in the case of 2-input systems. In fact, we can immediately find an entire class of examples that satisfy the sharp Sluis--Tilbury bound of $2n-3$ prolongations of a 2-input system where $n$ is the number of states.

Indeed, consider any contact sub-connection $\CH_G$ with 2 controls whose Brunovsk\'y component has signature $\langle 1,0,\dots, 0,1\rangle$, where the final `1' is in position $k$
for some $k>1$. If the group~$G$ has dimension~$r$, then the contact sub-connection, as a control system, has $1+k+r$ state variables. If a \stf\ linearizable partial contact curve reduction is performed along the contact coordinate of order 1 then the signature of the Brunovsk\'y component of the reduced sub-connection
$\bar{\CH}_G$ will have signature $\langle 0,0,\dots, 0,1\rangle$, where the final `1' continues to be in position~$k$. This means that $\bar{\CH}_G$
will be a rank 2 Goursat bundle and its derived length will be $\bar{k}=k+r$, where upon Corollary \ref{prolongationPredictor-cor} tells us that $\CH_G$ requires at most
\[
2\bar{k}-1=2(k+r)-1=2(1+k+r)-3=2n-3
\]
differentiations to achieve a SFL system.

Actually, the bound in Corollary \ref{prolongationPredictor-cor} is far more efficient
than the best known bounds \cite{ImprovedST-2006} on the required number of partial prolongations; there they are polynomial and/or exponential in the number of states and controls, in contrast to the linear bounds provided by Corollary \ref{prolongationPredictor-cor}. However, in \cite{ImprovedST-2006} the authors permit themselves no additional information on the control systems in question. In particular, unlike in the present paper no assumption on admitted symmetry is made.\looseness=1

\begin{Example}\label{4c prolong predict ex}
Consider the following 4-input contact sub-connection on $J^{\langle 1,2,1 \rangle}\times G$, where $G$ is an abelian Lie group with local coordinates $\bigl(\e^1,\e^2,\e^3\bigr)$,
\begin{equation*}
\CH_G=\operatorname{span}\bigl\{\mathbf{D}_t+w^1_2z^1_3\partial_{\e^1}+w^1_2z^2_1\partial_{\e^2}+w^2_2z^2_1\partial_{\e^3}, \, \partial_{z^1_3},\,\partial_{z^2_1},\,\partial_{w^1_2},\,\partial_{w^2_2}\bigr\},
\end{equation*}
where $\mathbf{D}_t$ is the total derivative operator on $J^{\langle 1,2,1 \rangle}$. Furthermore, the coordinates of $J^\kappa=J^{\langle 1,2,1\rangle}$ are such that if we decompose $\kappa=\u^\perp+\u$ such that $\u^\perp=\langle 1,0,1 \rangle$ and $\u=\langle 0,2,0\rangle$ then the~$\bsy{z}$ coordinates correspond to \smash{$J^{\u^\perp}$} and the $\bsy{w}$ coordinates correspond to $J^{\u}$. The refined derived type of $\CH_G$ is
\begin{equation*}
[[5,0],[9,0,0],[15,8,13],[16,16]],
\end{equation*}
which is clearly not that of a Brunovsk\'y form and hence this system is not SFL. However, let us illustrate Theorem \ref{prolongationPredictor} by reducing $\CH_G$ by a family of partial contact curves \smash{$\bsy{c}^\u_{\bsy{f}}$} that entail integral manifolds of $\bsy{\b}^{\u}$. Indeed, let $\bsy{f}(t)=(f_1(t),f_2(t))$ such that $\bsy{c}^\u_{\bsy{f}}$ imposes \[
\bigl(w^i_0,w^i_1,w^i_2\bigr)=\bigl(f_i(t),\dot{f}_i(t),\ddot{f}_i(t)\bigr),\qquad i=1,2.\]
 The reduced contact sub-connection is therefore
\begin{equation*}
\bar{\CH}_G=\operatorname{span}\bigl\{\mathbf{\bar{D}}_t+z^1_3\ddot{f}_1(t)\partial_{\e^1}+z^2_1\ddot{f}_1(t)\partial_{\e^2}+z^2_1\ddot{f}_2(t)\partial_{\e^3},\, \partial_{z^1_3},\, \partial_{z^2_1}\bigr\},
\end{equation*}
where $\mathbf{\bar{D}}_t$ is the total derivative operator on $J^{\u^\perp}$.
By calculation, we find that the reduced system~$\bar{\CH}_G$ is SFL with signature $\bar{\k}=\langle 0,0,1,1 \rangle$. The fundamental functions of $\bar{\CH}_G$ are easy to compute,
\begin{align*}
\phi^{1,3}&=\frac{\dddot{f}_2(t)}{\dddot{f}_1(t)}\bigl(\ddot{f}_1(t)z^2_0-\e^2\bigr)+\bigl(\e^3-\ddot{f}_2(t)z^2_0\bigr),\\
\phi^{1,4}&=\e^1-\ddddot{f}_1(t)z^1_0+\dddot{f}_1(t)z^1_1-\ddot{f}_2(t)z^1_2.
\end{align*}
Note that $j^\u \bsy{f}(t)$ has signature $\u=\langle 0,2 \rangle$ and the jet $j^{\u^{1,3}} \bsy{f}(t)$ appearing in $\phi^{1,3} (\ell=3)$ is of signature $\u^{1,3}=\langle 0,0,2\rangle$, which is predicted by the maximal bound of Lemma \ref{fund func sig}. For $\phi^{1,4}(\ell=4)$ the signature of the jet of $\bsy{f}(t)$ that appears has signature $\u^{1,4}=\langle 0,1,0,1\rangle$, which is not the maximum signature predicted by Lemma \ref{fund func sig}, since the lemma says the signature is no greater than $\u_3=\langle 0,0,0,0,2 \rangle$. We have
$\u^{1,3}_3=\langle 0,0,0,0,0,2 \rangle$ and $\u^{1,4}_4=\langle 0,0,0,0,0,1,0,1 \rangle$, so by Theorem \ref{prolongationPredictor} we set $S=\bigl\{\u^{1,3}_3,\u^{1,4}_4\bigr\}$ and obtain $\u'=\u^*_S=\langle 0,0,0,0,0,1,0,1 \rangle$. Finally, recalling that our partial contact curve reduction occurred along the two second order jet coordinates~$w^1$,~$w^2$, we prolong $\bsy{\b}^\u$ to $\bsy{\b}^{\u'}$ so that
\[
\u^\perp+\u=\langle 1, 2, 1\rangle=\langle 1,0,1\rangle+\langle0,2,0\rangle\mapsto
\langle1,0,1\rangle+\langle 0,0,0,0,0,1,0,1\rangle=\u^\perp+\u',
\]
and thus we manage to construct $\operatorname{pr}\CH_G$ on $J^{\langle 1,0,1,0,0,1,0,1 \rangle}\times G$ as
\begin{equation*}
\operatorname{pr} \CH_G=\operatorname{span}\bigl\{\operatorname{pr} \mathbf{D}_t+w^1_2z^1_3\partial_{\e^1}+w^1_2z^2_1\partial_{\e^2}+w^2_2z^2_1\partial_{\e^3},\,\partial_{z^1_3},\partial_{z^2_1},\,\partial_{w^1_8}, \,\partial_{w^2_6}\bigr\},
\end{equation*}
where $\operatorname{pr} \mathbf{D}_t$ is the total derivative operator on $J^{\langle 1,0,1,0,0,1,0,1 \rangle}$. One can then apply Theorem~\ref{Goursat SFL} or the GS algorithm to discover that $\operatorname{pr}\CH_G$ is indeed SFL.
The signature of $\operatorname{pr}\CH_G$ turns out to be $\langle 0,0,1,1,0,1,0,1 \rangle$. The GS algorithm or {\it contact} can now be used to construct an explicit solution of $\operatorname{pr}\CH_G$ and hence that of
$\CH_G$. This verifies Theorem \ref{prolongationPredictor} in that the signature of $\operatorname{pr} \CH_G$ is
$\langle 0,0,1,1,0,1,0,1 \rangle=\bar{\k}+\u'$, this being the signature of the {\em explicit solution} of
$\CH_G$.

Finally, if we did not have the expressions for the fundamental functions we could nevertheless obtain an upper bound for the number of prolongations using Corollary \ref{prolongationPredictor-cor}. Since in this case~${\bar{k}\!=\!4}$, we would prolong
$\u=\langle 0,2\rangle$ by $2\bar{k}-1=7$ derivatives giving ${\nu_{S'}^*=\langle 0,0,0,0,0,0,0,0,2\rangle}$ instead of ${\nu'=\langle 0,0,0,0,0,\allowbreak 1,0,1\rangle}$. Lemma \ref{prolongationOFsflConnection} says that a prolongation of $\nu$ to $\nu_{S'}^*$ instead of $\u'$ ensures that the resulting prolonged sub-connection
$\operatorname{pr} \CH_G$ is also SFL.
\end{Example}

\begin{Remark}
Theorem \ref{prolongationPredictor} and Corollary \ref{prolongationPredictor-cor} do more than provide a bound on the number of partial prolongations required to produce a SFL control system $\operatorname{pr}\CH_G$. They also tell us which controls need to be prolonged.
Moreover, as is by now clear, cascade linearization gives a~general canonical procedure for expressing a control system that is not DFL by differentiation of given inputs to one that is.
\end{Remark}

\begin{Example}[Example \ref{pvtol sub-connection ex} continued]\label{pvtol dfl 1}
This makes use of the contact connection worked out in \eqref{connectionVTOLgalilean} and \eqref{connectionFormVTOLgalilean}.

We now check the cascade feedback linearization of \eqref{tVTOL1}. We can perform a contact curve reduction along
$m$ or $n$. For this example we fix a family of contact curve reductions by defining the
map
$\mathbf{c}_f^{\langle 0, 1\rangle}\colon J^{\langle 1\rangle}\times G\to J^{\langle 1, 1\rangle}\times G$
given by
\begin{align*}
\begin{aligned}
\mathbf{c}_f^{\langle 0, 1\rangle}\bigl(t, n, n_1,\e^1, \e^2, \e^3\bigr)&=\bigl(t, n, n_1, \e^1, \e^2, \e^3, m, m_1, m_2\bigr)\\
&=\bigl(t, n, n_1, \e^1, \e^2, \e^3, f(t), \dot{f}(t),\ddot{f}(t)\bigr).
\end{aligned}
\end{align*}
Then
\begin{align*}
\bsy{\bar{\g}}^G={}&\bigl(\mathbf{c}_f^{\langle 0, 1\rangle}\bigr)^*\bsy{\g}^G=
\{{\rm d}n-n_1{\rm d}t\}\cr
&\oplus\left\{{\rm d}\e^1+\frac{(n_1+1)\sin f-h\ddot{f}}{\cos f}{\rm d}t, \, {\rm d}\e^2-t\frac{(n_1+1)\sin f-h\ddot{f}}{\cos f}{\rm d}t,\,
{\rm d}\e^3-n{\rm d}t\right\}.
\end{align*}
It can be checked that $\bsy{\bar{\g}}^G$ is \sfl\ of signature $\langle 0, 0, 0, 1\rangle$ and therefore
\stf\ equivalent to the Brunovsk\'y normal form with that signature. Moreover, by Theorem \ref{prolongationPredictor}, we perform a~7 \big($=2\bar{k}-1$\big)-fold prolongation in $m$, in order that $\operatorname{pr}\CH_G$ is \sfl, giving a dynamic feedback linearization of $\CH_G$. In fact, it turns out that it is sufficient to prolong to order 6 in $m$. That is,
\begin{gather*}
\operatorname{pr}\CH_G= \operatorname{span}\biggl\{\P t+m_1\P m+m_2\P {m_1}+m_3\P {m_2}+
m_4\P {m_3}+m_5\P {m_4}+m_6\P {m_5}+n_1\P n\\
\hphantom{\operatorname{pr}\CH_G= \operatorname{span}\biggl\{ }{} +\frac{(n_1+1)\sin m-hm_2}{\cos m}(t\P {\e^2}-\P {\e^1})+n\P {\e^3}, \P {m_6}, \P {n_1}\biggr\}
\end{gather*}
is \sfl\ of signature $\langle 0, 0, 0, 1, 0, 1\rangle$, which is also the signature of the explicit solution.
\end{Example}

\subsection{From cascade to dynamic feedback linearization}\label{standardDFL}

So far we have shown how to construct the explicit solution of a smooth, invariant control system \eqref{controlSystem} via the contact sub-connection, $\CH_G$.
The goal of this subsection is to give the precise relationship between cascade feedback linearization and \df\ linearization and show how to construct a dynamic feedback linearization of an invariant control system \eqref{controlSystem} in local coordinates, as expressed by Definition \ref{dynamic-feedback-def}. The basic idea is that the map $\widetilde{\varphi}$ should lift to a map \smash{$\wt{\wt{{\varphi}}}$} between
\[
\bigl(J^{\nu'+\nu^\perp}\times G, \bsy{\b}^{\nu'+\nu^\perp}\oplus \Theta^G\bigr)
\]
 and some manifold and Pfaffian system \smash{$\bigl(\widecheck{M},\widecheck{\bsy{\o}}\bigr)$} with $\widecheck{\bsy{\o}}$ encoding a dynamic extension of
$\bsy{\o}$.

\medskip

\emph{Constructing a \df\ linearization:}
\begin{enumerate}\itemsep=0pt
\item[(1)] Construct $\operatorname{pr}\CH_G$ and $\bsy{\widecheck{\bsy{\g}}}^G:=\operatorname{ann} (\operatorname{pr}\CH_G)$ as prescribed in Theorem \ref{prolongationPredictor}.
\item[(2)] Compute the Goursat bundle \smash{$\theta=\bigl(\wt{\wt{\varphi}}\bigr)^*\bsy{\widecheck{\bsy{\g}}}^G$} and notice that the state and control variables are not yet well defined.
\item[(3)] For $\mcal{T}=\operatorname{ann} \theta$, choose any $t$-preserving diffeomorphism \smash{$\chi\colon\widecheck{M}\to\widecheck{M}$} such that the integrable bundle \smash{$\ch\mcal{T}^{(1)}_0$} is spanned by coordinate vector fields
$\{\partial_{W^1},\dots,\partial_{W^m}\}$ according to the Frobenius theorem. The diffeomorphism $\chi$ defines the dynamic compensator~${u=\b(t, x, y, W)}$ as in Definition
\ref{dynamic-feedback-def} and thus defines the new control variables~$W$ and new additional state variables $y$ for the SFL control system \smash{$\widecheck{\mcal{V}}:=\vartheta^{-1}_*\left(\operatorname{pr}\CH_G\right)$} where~\smash{$\vartheta=\wt{\wt{\bsy{\varphi}}}\circ \chi^{-1}$}.
\end{enumerate}
\begin{figure}[ht]
\centering
\begin{tikzcd}[every label/.append style = {font = \normalsize}]
{\bigl(\widecheck{M},\widecheck{\boldsymbol{\omega}}\bigr)} \arrow[d, "\chi^{-1}"'] \arrow[rd, "\vartheta"] & \\
{\bigl(\widecheck{M},\widecheck{\boldsymbol{\omega}}\bigr)} \arrow[r, "\widetilde{\widetilde{\boldsymbol{\varphi}}}"] & {\bigl(J^{\nu^\perp+\nu'}\times G,\widecheck{\boldsymbol{\gamma}}^G\bigr)}
\end{tikzcd}
\caption{Deformation of \smash{$\wt{\wt{\varphi}}$} to a static feedback transformation $\vartheta$.}
\label{diagramDeformation}
\end{figure}
\begin{Theorem}\label{cfl to dfl}
If a control system $(M,\bsy{\o})$ is cascade feedback linearizable with respect to a control admissible symmetry group $G$, then there exists a dynamic compensator $\bsy{u}=\bsy{\b}(t,\bsy{x},\bsy{y},\bsy{W})$ together with linear dynamics in new state variables $\bsy{y}$ such that the dynamic augmentation is static feedback linearizable.
\end{Theorem}
\begin{proof}
Let \smash{$\widecheck{M}$} be a manifold that locally has the form of an open subset of $P\times M$, where $P$ is some manifold with $\dim P=\dim J^{\nu'}-\dim J^\nu$ and local coordinates denoted by $(\bsy{p})$. Moreover, let \smash{$\widetilde{\pi}\colon \widecheck{M} \to M$} be the natural projection with typical fiber given by $P$. Let $(t,\bsy{z},\bsy{w},\bsy{w}')$ be local coordinates on \smash{$J^{\nu'+\nu^\perp}$} with $\bsy{w}'$ denoting the coordinates along the fibers of \smash{$J^{\nu'+\nu^\perp}\to J^{\nu+\nu^\perp}$}. Then we define the map \smash{$\wt{\wt{\varphi}}\colon \widecheck{M} \to J^{\nu'+\nu^\perp}\times G$} by
\begin{equation*}
\wt{\wt{\varphi}}=\wt{\varphi}\circ\wt{\pi} \times \{\bsy{w}'=\bsy{p}\}.
\end{equation*}
Furthermore, we will define the Pfaffian system and dual distribution by
\begin{align*}
\bsy{\th}&=\bigl(\wt{\wt{\varphi}}\bigr)^* \widecheck{\bsy{\g}}^G,\qquad
\mcal{T}=\operatorname{ann} \bsy{\th}.
\end{align*}
Notice that \smash{$\wt{\wt{\varphi}}$} is \textit{not} a static feedback transformation since the state and control variables for~$\mcal{T}$ are not yet well defined; however it does leave $t$ unchanged. Additionally, and again because the state and control variables are not well defined, it is not obvious that $\mcal{T}$ is static feedback equivalent to $\operatorname{pr} \CH_G$ and hence static feedback linearizable. However, since
$\wt{\wt{\varphi}}
$ is a diffeomorphism preserving the independent variable $t$, and since the refined derived type is a diffeomorphism invariant, we do have that $\mcal{T}$ is still a Goursat bundle whose integral curves are parameterized by $t$. The intersection bundle \smash{$\ch\mcal{T}^{(1)}_0$} will allow us to make a well defined and canonical choice of state and control variables via a diffeomorphism, as we now demonstrate.

Let $\widecheck{\Xi}^{(1)}_0$ be the annihilator of \smash{$\ch (\operatorname{pr} \CH_G)^{(1)}_0$}. Then,
\begin{equation*}
\widecheck{\Xi}^{(1)}_0=\operatorname{span}\bigl\{{\rm d}t,{\rm d}\e,{\rm d}\bsy{z}^{\floor{\nu^\perp}},{\rm d}\bsy{w}^{\floor{\nu}},{\rm d}\bsy{w}^{\ceil{\nu}},{\rm d}\bsy{w}^{\floor{\nu'-\nu}}\bigr\},
\end{equation*}
where $\bsy{w}^{\floor{\nu'-\nu}}$ denotes all the $\bsy{w'}$ variables of order strictly smaller than the highest order $\bsy{w'}$ variables. The first integrals of this bundle are precisely all the state variables for the control system $\operatorname{pr} \CH_G$ and time $t$. Additionally, the flow box coordinates for \smash{$\ch (\operatorname{pr} \CH_G)^{(1)}_0$} provide the control variables, which are precisely \smash{$\bigl(\bsy{z}^{\ceil{\nu^\perp}},\bsy{w}^{{\ceil{\nu'}}}\bigr)$} for $\operatorname{pr} \CH_G$.
With this in mind, we now~find~the pullback of the bundle $\widecheck{\Xi}^{(1)}_0$ by \smash{$\wt{\wt{\varphi}}$} to be
\begin{equation*}
\begin{aligned}
\bigl(\wt{\wt{\varphi}}\bigr)^*\widecheck{\Xi}^{(1)}_0=\operatorname{span}\bigl\{{\rm d}t,{\rm d}\bsy{x},{\rm d}\bigl(\bsy{w}^{\ceil{\nu}}(t,\bsy{x},\bsy{u})\bigr),{\rm d}\bsy{p}^{\floor{\nu'-\nu}}\bigr\},
\end{aligned}
\end{equation*}
since
\begin{equation*}
\bigl(\wt{\wt{\varphi}}\bigr)^*\operatorname{span}\bigl\{{\rm d}t,{\rm d}\e,{\rm d}\bsy{z}^{\floor{\nu^\perp}},{\rm d}\bsy{w}^{\floor{\nu}}\bigr\}=\operatorname{span}\{{\rm d}t,{\rm d}\bsy{x}\},
\end{equation*}
by construction of $\wt{\wt{\varphi}}$ from $\wt{\varphi}$. Without loss of generality, let us relabel the indices of the control variables $\bigl(u^1,\dots, u^m\bigr)$ of $(M,\bsy{\o})$ so that of the $m=m_{\nu^\perp}+m_{\nu}$ controls, the last $m_{\nu}$ control variables $\bigl(u^{m_{\nu^\perp}+1},\dots,u^{m_\nu}\bigr)$ may be determined from the equations $\bsy{y}=\bigl(\bsy{w}^{\ceil{\nu}}(t,\bsy{x},\bsy{u}),\bsy{p}^{\floor{\nu'-\nu}}\bigr)$ in terms of $\bsy{x}$, $\bsy{y}$, and $\bigl(u^1,\dots,u^{m_{\nu^\perp}}\bigr)$. Let \smash{$\chi\colon\widecheck{M}\to\widecheck{M}$} be the diffeomorphism whose inverse is given by
\begin{equation*}
(t,\bsy{x},\bsy{u},\bsy{p})\mapsto(t,\bsy{x},\bsy{y},\bsy{W}),
\end{equation*}
where
\begin{align}
\bsy{y}&=\bigl(\bsy{w}^{\ceil{\nu}}(t,\bsy{x},\bsy{u}),\bsy{p}^{\floor{\nu'-\nu}}\bigr),\qquad
\bsy{W}=\bigl(u^1,\dots,u^{m_{\nu^\perp}},\bsy{p}^{\ceil{\nu'}}\bigr).\label{diffeo compensator}
\end{align}
Using this diffeomorphism, we find that
\begin{equation*}
\bigl(\wt{\wt{\varphi}}\circ \chi^{-1}\bigr)^*\widecheck{\Xi}^{(1)}_0=\operatorname{span}\{{\rm d}t,{\rm d}\bsy{x},{\rm d}\bsy{y}\}.
\end{equation*}

It follows that if $X$ is a vector field on $J^{\u^\perp+\u'}\times G$ such that
\[
\bigl(\wt{\wt{\vf}}\circ\chi^{-1}\bigr)_*X\in\operatorname{ann} \widecheck{\Xi}^{(1)}_0=
\ch(\operatorname{pr}\CH_G)^{(1)}_0
\]
 then
 \[\chi^{-1}_*X\in\chi^{-1}_*\ch (\operatorname{pr}\CH_G)^{(1)}_0=\ch\CT^{(1)}_0.\]
 Therefore, $\chi_*\ch\CT^{(1)}_0$ is annihilated by $\operatorname{span}\{{\rm d}t,{\rm d}\bsy{x},{\rm d}\bsy{y}\}$, from which we have
\begin{equation*}
\chi_*\bigl(\ch \mcal{T}^{(1)}_0\bigr)=\operatorname{span}\{\partial_{W^1},\dots,\partial_{W^m}\}.
\end{equation*}
Now we have that
\[
\chi_*\CT=\chi_*\wt{\wt{\vf}}^{-1}_*(\operatorname{pr}\CH_G):=
\psi_*(\operatorname{pr}\CH_G)
\]
 and since $\operatorname{pr}\CH_G$ is SFL and \smash{$\psi=\chi\circ\wt{\wt{\vf}}^{-1}$} is a diffeomorphism we have that $\chi_*\CT$ is a Goursat bundle and so the latter has a Cauchy bundle~${\ch\bigl(\chi_*\CT^{(k-1)}\bigr)}$ with an annihilator which we will denote by $\L^{(k'-1)}$, where $k'$ is the derived length
of $\CT$; therefore, we have that
\[
\L^{(k'-1)}=\bigl(\psi^{-1}\bigr)^*\widecheck{\Xi}^{(k'-1)}.
\]
 In the case~${\Delta_{k'}=1}$, we have that \smash{${\rm d}t\in\widecheck{\Xi}^{(k'-1)}$} and hence since $\psi$ preserves $t$
$(t\circ\psi=t)$, we see that~${{\rm d}t\in\L^{(k'-1)}}$. Therefore, in the case $\Delta_{k'}=1$, the hypotheses of Theorem \ref{Goursat SFL} are satisfied and we can conclude that $\chi_*\CT$ is SFL. A similar argument holds in the case
$\Delta_{k'}>1$.

Let \smash{$\widecheck{\CV}:=\chi_*\mcal{T}$} and let \smash{$\widecheck{\bsy{\o}}:=\operatorname{ann} \widecheck{\CV}$}. Then, $\widecheck{\bsy{\o}}$ is precisely the dynamic feedback augmentation of~$\bsy{\o}$ with dynamic compensator $\bsy{u}=\bsy{\beta}(t,\bsy{x},\bsy{y},\bsy{W})$ given by solving equations \eqref{diffeo compensator} for $\bsy{u}$. We have that the map \smash{$\vartheta:=\wt{\wt{\varphi}}\circ \chi^{-1}\colon \widecheck{M}\to J^{\nu^\perp+\nu'}\times G$} is a static feedback transformation, and since~\smash{$\widecheck{\bsy{\g}}^G$} is SFL, then $\widecheck{\bsy{\o}}$ must also be SFL.
Lastly, we show that the associated dynamics for $\bsy{y}$ are linear. Denote the $\nu'$ prolongation of
$\bsy{\gamma}^G$ by
\[
\alpha=\operatorname{span}\bigl\{{\rm d}\bsy{w}^{\ceil{\nu}}-\bsy{w}_1^{\floor{\nu'-\nu}} {\rm d}t,\dots,{\rm d}\bsy{w}_{l_{\nu'}-1}^{\floor{\nu-\nu'}}-\bsy{w}^{\ceil{\nu'}}{\rm d}t\bigr\},
\]
where $l_{\nu'}$ is a place holder for the highest order jet variable determined by Theorem \ref{prolongationPredictor}. Then we merely observe that a pullback by $\wt{\wt{\varphi}}$ of $\alpha$ results in
\[
\bigl(\wt{\wt{\varphi}}\bigr)^*\alpha=\operatorname{span}\bigl\{{\rm d}\bigl(\bsy{w}^{\ceil{\nu}}(t,\bsy{x},\bsy{u})\bigr)-\bsy{p}_1^{\floor{\nu'-\nu}} {\rm d}t,\dots,
{\rm d}\bsy{p}_{l_{\nu'}-1}^{\floor{\nu-\nu'}}-\bsy{p}^{\ceil{\nu'}}{\rm d}t\bigr\}.
\]
Thus, one final pullback by $\chi^{-1}$ results in
\[
\bigl(\wt{\wt{\varphi}}\circ \chi^{-1}\bigr)^*\alpha=\operatorname{span}\bigl\{{\rm d}\bsy{y}_0-\bsy{y}_1 {\rm d}t,\dots,{\rm d}\bsy{y}_{l_\nu'-1}-\bsy{W}^{\ceil{\nu'}} {\rm d}t\bigr\},
\]
where the subscripts on $\bsy{y}_i$ for $1\leq i\leq l_{\nu}$ are used to keep track of how the $\bsy{y}$ variables are assigned to the $\bsy{p}$ variables and $\bsy{W}^{\ceil{\nu'}}$ denotes those control variables arising from the highest order $\bsy{p}$ variables. Thus, we can see that the above Pfaffian system is that of a linear system.
\end{proof}

In Example \ref{pvtol dfl 1}, the coordinate change $\chi$ in Step 3 turns out to be the identity transformation, as we demonstrate at the end of this section. This can always be arranged for systems which are \dfl\ by differentiation of the {\it given} inputs. We note also that sometimes $\ch\mcal{T}^{(1)}_0$ is already in a basis of coordinate vector fields, which defines the control variables
for a dynamic feedback linearization of $\CV$ without the need for the map $\chi$ (which can be chosen to be the identity). This choice can be made at the cost of having nonlinear dynamics for the variables $\bsy{p}$; see Section \ref{PVTOLsection} for an example. We also have the following corollary concerning prolongation by differentiation.

\begin{Corollary}\label{state-sym-CFL}
Let $(M,\bsy{\o})$ be a CFL control system with respect to a state-space symmetry group $G$. Then $(M,\bsy{\o})$ is DFL by differentiation of the given inputs.
\end{Corollary}
\begin{proof}
The proof of Theorem \ref{cfl to dfl} essentially shows that the dynamic compensator consists of components of the inverse of the map $\wt{\wt{\varphi}}$ which, in turn, is constructed from $\wt{\varphi}$. The inverse of $\wt{\varphi}$ may be written in terms of the inverse of $\varphi$ which gives expressions for the control variables $\bsy{v}$ on $M/G$. That is, $\varphi^{-1}$ acts as the identity on $\bsy{v}$. However, if $G$ is only a state-space symmetry group, then the controls of $(M,\bsy{\o})$ descend to $(M/G,\bsy{\o}/G)$ since they are invariant functions of the action of $G$, i.e., $\bsy{v}=\bsy{u}$. As a result, we have $\bsy{u}=\bigl(\bsy{W},\bsy{w}^{\ceil{\nu}}\bigr)$ which is a dynamic compensator by partial prolongation.
\end{proof}

The need for -- and the construction of -- a nontrivial $\chi$ in the general case will be illustrated by the example in Section~\ref{tVTOL Ex}. For now, we illustrate Corollary \ref{state-sym-CFL}.
\begin{Example}\label{CharletFinal}
We continue with Example \ref{CharletExCont'd}. Recall that the contact sub-connection can be prolonged to the SFL system
\[
\widecheck{\bsy{\g}}^G=\bsy{\g}^G\oplus\{{\rm d}w_2-w_3 {\rm d}t\}
\] and that the relevant symmetry of the original control system $\bsy{\omega}$ is translation in $x^4$, i.e., a state-space symmetry. We construct \smash{$\widecheck{M}$} to be $M\times P$ where $P=\mathbb{R}$ and has coordinate $p_1$. Then the map $\wt{\wt{\bsy{\varphi}}}$ is
\begin{align*}
\begin{aligned}
&t=t,\qquad \varepsilon=x^4,\qquad z=x^3,\qquad z_1=u^2,\\
& w=x^1, \qquad w_1=x^2,\qquad w_2=u^1,\qquad w^3=p_1.
\end{aligned}
\end{align*}
It is immediate that $\widecheck{\bsy{\g}}^G$ pulls back by $\wt{\wt{\bsy{\varphi}}}$ to
$\bsy{\theta}=\bsy{\omega} \oplus \operatorname{span}\bigl\{{\rm d}u^1-W^2 {\rm d}t \bigr\}$,
and since
\[
\operatorname{Char}\mcal{T}^{(1)}_0=\operatorname{span}\{\partial_{u^2},\partial_{p_1}\}
\]
we find that $\chi$ is an identity map and so the dynamic compensator is $u^2=W^1$ and $u^1=y_0$ and the new controls are $W^2=p_1$ and $W^1=u^2$. Therefore, the control system $\bsy{\omega}$ of Example~\ref{ExmpResolventRelativeGoursat} is DFL by one partial prolongation of~$u^1$.
\end{Example}

Finally, we summarize the relationship between CFL control systems and DFL control systems by the diagram
 in Figure~\ref{cfl to dfl diag}.
\begin{figure}[ht]\centering
\begin{tikzcd}[every label/.append style = {font = \normalsize}]
{\bigl(\widecheck{M},\widecheck{\boldsymbol{\omega}}\bigr)} \arrow[d, "\chi^{-1}"'] \arrow[rd, "\vartheta"] & & & \\
{\bigl(\widecheck{M},\widecheck{\boldsymbol{\omega}}\bigr)} \arrow[d, "\wt{\pi}"'] \arrow[r, "\wt{\wt{\varphi}}"] & {\bigl(J^{\nu^\perp+\nu'}\times G, \widecheck{\boldsymbol{\gamma}}^G\bigr)} \arrow[d, "\wt{\pi}'"] \arrow[r, "\varphi'"] & {\bigl(J^{\bar{\kappa}+\nu'}, \boldsymbol{\beta}^{\bar{\kappa}+\nu'}\bigr)} & \\
{(M, \boldsymbol{\omega})} \arrow[r, "\wt{\varphi}"] \arrow[d, "\pi"'] \arrow["G"', loop, distance=2em, in=215, out=145] & {\bigl(J^{\kappa}\times G, \boldsymbol{\gamma}^G\bigr)}\arrow[d, "\pi'"] & {\bigl(J^{\nu^\perp}\times G, \boldsymbol{\bar{\gamma}}^G\bigr)} \arrow[l, "\boldsymbol{c}^\nu_{\boldsymbol{f}}", hook'] \arrow[r, "\bar{\varphi}"] \arrow[r, "(2)", phantom, shift right=3] & {\bigl(J^{\bar{\kappa}}, \boldsymbol{\beta}^{\bar{\kappa}}\bigr)} \\
{(M/G, \boldsymbol{\omega}/G)} \arrow[r, "\varphi"] \arrow[r, "(1)", phantom, shift right=3] & { (J^\kappa, \boldsymbol{\beta}^\kappa)} & &
\end{tikzcd}
\caption{In this diagram the signature $\kappa$ decomposes as $\kappa=\nu+\nu^\perp$. Additionally, the maps $\varphi$, $\bar{\varphi}$, $\wt{\varphi}$, and $\vartheta$ are all SFTs.}
\label{cfl to dfl diag}
\end{figure}

\begin{Example}\label{pvtolExample0}
We continue with Example \ref{pvtol dfl 1}. According to the outlined method, we now compute $\wt{\wt{\varphi}}$ given by
\begin{align*}
&t=t,\qquad m=\th,\qquad m_1=\th_1,\qquad m_2=v^2,\qquad n=z_1,\\
 &n_1=v^1\cos\th+hv^2\sin\th-1,\qquad \e^1=x_1,\qquad \e^2=x-tx_1,\\& \e^3=z,\qquad m_3=p_1,\qquad m_4=p_2,\qquad m_5=p_3,\qquad m_6=p_4.
\end{align*}
We then find that
\begin{align*}
\mcal{T}={}&\bigl(\wt{\wt{\varphi}}\bigr)^{-1}_*\operatorname{pr}\CH_G\\
={}&
\bigl\{\P t+x_1\P x-\bigl(u^1\sin\th-h u^2\cos\th\bigr)\P {x_1}+z_1\P z+\bigl(u^1\cos\th+h u^2\sin\th-1\bigr)\P {z_1}\\
&+\th_1\P {\th}+u^2\P {\th_1}+p_1\P {u^2}+p_2\P {p_1}+p_3\P {p_2}+p_4\P {p_3}, \,\P {u^1}, \,\P {p_4}\Big\}.
\end{align*}
This is precisely a dynamic feedback linearization by differentiation of the PVTOL by a 4-fold prolongation of the control $u^2$, which can be compared to the bound of $2\bar{k}-1=7$ established in Corollary \ref{prolongationPredictor-cor}. Here $\chi$ is the identity transformation, but we still use the fact that
\[
\ch\mcal{T}^{(1)}_0=\operatorname{span}\{\partial_{u^1},\partial_{p_4}\}
\]
to guarantee that the new control variables are $u^1$ and $p_4$, respectively.
Thus $\vartheta=\wt{\wt{\varphi}}\circ \chi^{-1}=\wt{\wt{\varphi}}$ is the \stf\ transformation
\[
\vartheta\colon \ \widecheck{M} \to J^{\k'}\times G,
\]
where $\k'=\langle 1, 0, 0, 0, 0, 1\rangle$. The fundamental functions of $\CT$ are flat outputs of the PVTOL system \eqref{tVTOL1} with respect to the Galilean group generated by $\bsy{\G}$ in Example~\ref{pvtol dfl 1}.
\end{Example}

\begin{Remark}Example \ref{pvtolExample0} highlights the fact that an invariant flat control system will generally have distinct sets of flat outputs depending upon the Lie subgroup $G$ being considered. In each case the flat outputs are generated by the fundamental functions of the corresponding dynamic linearization of $\CV$ (i.e., the static feedback linearization of \smash{$\widecheck{\CV}$}). In Section \ref{PVTOLsection}, we will derive the well-known flat outputs of \eqref{tVTOL1} using a different subgroup of the control symmetry group for this system.
\end{Remark}

\section[A non-'integrator chain' dynamic feedback linearization example]{A non-`integrator chain' dynamic feedback\\ linearization example}\label{tVTOL Ex}

In this section, we illustrate the dynamic feedback linearization framework developed in this paper via an elementary example. From \cite{SluisTilbury}, we deduce that this system does not possess a~linearization by the differentiation of the {\it given} inputs.
The goal of this section is to systematically derive a \df\ linearization in local coordinates using only the framework established in this paper.

While the system in question cannot be linearized by integrator chains, it can nevertheless be transformed to such a system by a preliminary feedback transformation which can, in this case, be guessed by inspection. However, the theory developed in this paper does not rely on such {\it ad hoc} steps. In the search for a dynamic linearization, a preliminary feedback transformation, if one is required, {\em arises} canonically from our general theory, as will be demonstrated in this case.

The control system to be studied is
\begin{align}
&\dot{x}_1=x^2,\qquad
\dot{x}_2=-x^2u^1\bigl(x^1+x^5\bigr)+ u^2,\qquad
\dot{x}_3=\bigl(x^3\bigr)^2x^4,\nonumber\\
&\dot{x}_4=-x^5\bigl(x^1x^2u^1-u^2\bigr)+x^2u^1-1,\qquad
\dot{x}_5=x^6,\qquad
\dot{x}_6=u^2-x^1x^2u^1,\label{tVTOL}
\end{align}
in 6 states and 2 controls. We denote the ambient manifold of the system by $M$
and calculate that the abelian Lie algebra spanned by
\[
\G=\operatorname{span}\left\{ Y_1= t\P {x^1}+\P {x^2}-\frac{u^1}{x^2}\P {u^1}+tx^2u^1\P {u^2},\,
Y_2=\P {x^1}+x^2u^1\P {u_2},\,
 Y_3=\bigl(x^3\bigr)^2\P {x^3}\right\}
\]
generates a control admissible symmetry group.
Let $\CV\subset TM$ denote the given control system for~\eqref{tVTOL}. Away from the submanifold $x^2=0$ the functions
\[
t, \qquad q_1=x^4,\qquad q_2=x^5,\qquad q_3=x^6,\qquad v^1=x^2u^1,\qquad v^2=u^2-x^1x^2u^1
\]
are $G$-invariant, serve as local coordinates on $M/G$ and determine the projection $\pi$ in local coordinates.

By Theorem \ref{relStatFeedbackLin}, to determine the structure of the quotient of $\CV$ by $G$ we need only study the refined derived type of the distribution
$\widehat{\CV}=\CV\oplus\bsy{\G}$. We find that $\widehat{\CV}$ is a static feedback relative Goursat bundle (Definition \ref{StFrelGoursatDefn}) of signature $\langle 1,1\rangle$. By Theorem \ref{relStatFeedbackLin}, it follows that
$\CV/G$ is SFL with signature, $\langle 1,1\rangle$; that is, $(M/G, \CV/G)$ is static feedback equivalent to Brunovsk\'y normal form $\left(J^{\langle 1,1\rangle}, \CB_{\langle 1,1\rangle}\right)$. An easy calculation shows that the quotient control system
$\CV/G$ has the local form
\begin{equation}\label{quotient_tVTOL}
\CV/G=\operatorname{span}\bigl\{\P t+\bigl( v^2q_2+v^1-1\bigr)\P {q_1}+q_3\P {q_2}+v^2\P {q_3},\, \P {v^1},\, \P {v^2}\bigr\},
\end{equation}
where $v^1$, $v^2$ are the controls on $M/G$, and a further calculation verifies the above mentioned properties of the quotient.\footnote{The $G$-invariant functions $v^1$, $v^2$ are precisely those components of the preliminary feedback transformation one would ordinarily choose to simplify system~\eqref{tVTOL}. This highlights the comments made at the beginning of this section.}

Indeed, applying procedure {\it contact} \cite{VassiliouGoursatEfficient} to \eqref{quotient_tVTOL} determines the local diffeomorphism
$\varphi\colon M/G\allowbreak\to J^{\langle 1,1\rangle}$ given by
\begin{equation*}
\varphi=\bigl(t=t, \,z=q_1,\, z_1= v^2q_2+v^1-1,\, w=q_2,\, w_1=q_3,\, w_2=v^2\bigr),
\end{equation*}
where $(t, z, z_1, w, w_1, w_2)$ are the standard contact coordinates on $J^{\langle 1,1\rangle}$. Next we consider the static feedback linearization $\wt{\vf}$ constructed as per the discussion following equation \eqref{bundle sfl}. The transformation group $G$ generated by $\bsy{\G}$ consists of
\begin{align*}
&\bar{x}_1=x^1+t\e^1+\e^2,\qquad \bar{x}_2=x^2+\e^1,\qquad \bar{x}_3=\frac{x^3}{1-x^3\e^3},\cr
& \bar{u}^1=\frac{x^2u^1}{x^2+\e^1}, \qquad \bar{u}^2=u^2+\bigl(t\e^1+\e^2\bigr)x^2u^1,
\end{align*}
with the remaining variables fixed by $G$. Routine calculation leads to the section
\[
\varphi\colon\ \bigl(t,q^1,q^2,q^3,v^1,v^2\bigr)\mapsto \bigl(t,0,1,0,x^4,x^5,x^6,u^1,u^2\bigr)
\]
yielding $\e^1=x^2-1$, $\e^2=x^1-t\bigl(x^2-1\bigr)$, $\e^3=1-\frac{1}{x^3}$
and so
\begin{align*}
\wt{\varphi}=\biggl(&t=t,\, z=x^4,\, z_1= x^5\bigl(u^2-x^1x^2u^1\bigr)+x^2u^1-1, \,w=x^5,\, w_1=x^6,\\
 & w_2=u^2-x^1x^2u^1,\, \e^1=x^2-1,\, \e^2=x^1-t(x^2-1),\, \e^3=1-\frac{1}{x^3}\biggr).
\end{align*}
In particular, we compute the contact sub-connection $\CH_G=\wt{\varphi}_*\CV$ in this case to be
\[
\CH_G=\operatorname{span}\{\P t+z_1\P z+w_1\P w+w_2\P {w_1}+
\lambda\P {\e^1}+(1-t\lambda)\P {\e^2}+z\P {\e^3},\, \P {z_1},\, \P {w_2}\},
\]
where $\lambda=w_2\bigl(w^2+1\bigr)-w(z_1+1)$.

To determine the existence of a dynamic feedback linearization of the system \eqref{tVTOL}, we carry out a partial contact curve reduction (see Section \ref{contactReduction}, Definition \ref{contactReductionDef}). As there are 2 inputs, we can do this along either the order 1 variable $z$ or the order 2 variable $w$. In this case we choose $w$ and form the map
\smash{$\mathbf{c}^{\langle 0, 1\rangle}_{f}=
j^{\langle 0, 1\rangle} f\times \text{Id}_{J^{\langle 1\rangle}}\times\text{Id}_{G\to G}\colon
J^{\langle 1\rangle}\times G\to J^{\langle 1, 1\rangle}\times G$},
where~${j^{\langle 0, 1\rangle} f=\bigl(t,f(t), \dot{f}(t),\ddot{f}(t)\bigr)}$ and $f$ is an arbitrary, smooth real-valued function of $t$. As in Section~\ref{contactReduction}, we form the Pfaffian system
\smash{$\bigl(\mathbf{c}^{\langle 0, 1\rangle}_{f}\bigr)^*\bsy{\g}^G=:\bsy{\bar{\g}}^G$}.
Then
\[
\bar{\CH}_G=\operatorname{span}\bigl\{\P t+z_1\P z+\bar{\lambda}\P {\e^1}+\bigl(1-t\bar{\lambda}\bigr)\P {\e^2}+z\P {\e^3},\,
\P {z_1}\bigr\},
\]
where $\bar{\CH}_G$ is the kernel of $\bsy{\bar{\g}}^G$ on $J^{\langle 1\rangle}\times G$ and $\bar{\lambda}=\ddot{f}\bigl(f^2+1\bigr)-f(z_1+1)$. By the results in Sections~\ref{brun-form-subsect} and \ref{Goursat-subsect}, we find that
$\bar{\CH}_G$ is static feedback linearizable with signature $\langle 0, 0, 0, 1\rangle$. Thus, the derived length of $\bar{\CH}_G$ is $\bar{k}=4$ and by Corollary \ref{prolongationPredictor-cor}, a maximum 7-fold partial prolongation of $\CH_G$ along the $w$-series of contact coordinates is sufficient for a static feedback linearizable control system. In fact, in this case it is enough to perform a 4-fold partial prolongation to obtain,
\begin{align*}
\operatorname{pr}\CH_G= \operatorname{span}\{&\P t+z_1\P z+w_1\P w+w_2\P {w_1}+w_3\P {w_2}+w_4\P {w_3}+w_5\P {w_4}
+w_6\P {w_5}\\
&+\lambda\P {\e^1}+(1-t\lambda)\P {\e^2}+z\P {\e^3}, \,\P {z_1},\, \P {w_6}\},
\end{align*}
on $J^{\langle 1, 0, 0, 0, 0, 1\rangle}\times G$, and a further calculation shows that
\[
\operatorname{pr}\CH_G\simeq_{\text{SF}}\CB_{\langle 0, 0, 0, 1, 0, 1\rangle}.
\]
We have therefore proven that \eqref{tVTOL} is \cfl\ and hence \dfl\ with a dynamic extension of signature
$\langle 0, 0, 0, 1, 0, 1\rangle$. Our next goal is to derive the linearizable dynamic extension of \eqref{tVTOL},
according to Theorem \ref{cfl to dfl}.

For this we apply the results of Section \ref{cfl to dfl} by computing
$\wt{\wt{\vf}}=\wt{\vf}\circ\wt{\pi}\times \bigl\{w_3=p_1^1,\,w_4=p_2^1,\,w_5=p_3^1,\,w_6=p_4^1\bigr\}$. Indeed,
\begin{align*}
\wt{\wt{\vf}}= \biggl(&t=t,\, z=x^4, \,z_1= x^5\bigl(u^2-x^1x^2u^1\bigr)+x^2u^1-1,\, w=x^5,\, w_1=x^6,\\
 &w_2=u^2-x^1x^2u^1,\, w_3=p^1_1,\, w_4=p^1_2,\, w_5=p^1_3,\, w_6=p^1_4,\\
 & \e^1=x^2-1, \,\e^2=x^1-t(x^2-1), \,\e^3=1-\frac{1}{x^3}\biggr).
\end{align*}
Calculating $\mcal{T}:=\wt{\wt{\vf}}^{-1}_*\operatorname{pr}\CH_G$ shows that this distribution has the form
\begin{equation*}
\mcal{T}=\operatorname{span}\bigl\{Y, \P {u^1}+x^1x^2\P {u^2},\, \P {p^1_4}\bigr\},
\end{equation*}
where $Y$ is the image of the first vector field in $\operatorname{pr}\CH_G$; in particular,
\begin{equation}\label{FirstAdaptation}
\ch\mcal{T}^{(1)}_0=\operatorname{span} \bigl\{\P {u^1}+x^1x^2\P {u^2},\, \P {p^1_4}\bigr\}.
\end{equation}
Because of the first vector field in \eqref{FirstAdaptation}, $\mcal{T}$ does not yield a well-defined control system in its current form.
Therefore, we compose \smash{$\wt{\wt{\vf}}$} with a map $\chi$ as in Figure~\ref{diagramDeformation}, embodying the change of variable
$u=\b(t,x,y,W)$ of Definition \ref{dynamic-feedback-def}, chosen so that the flowbox coordinates of \smash{$\ch\mcal{T}^{(1)}_0$} are the new controls.

In this case, it is easy to see that \eqref{FirstAdaptation} is transformed by
\begin{align*}
\chi^{-1}={} &\bigl(t=t,\, x^i=x^i,\, y^1=u^2-x^1x^2u^1,\, y^2=p^1_1,y^3=p^1_2,\, y^4=p^1_3,\, W^1=u^1,\, W^2=p^1_4\bigr),
\end{align*}
where $i=1,\dots,6.$
We see that the new controls are $W^1=u^1,\ W^2=p^1_4$ with additional states given by $y^1=u^2-x^1x^2u^1$, $y^2=p^1_1$, $y^3=p^1_2$, $y^4=p^1_3$. Thus the map $u=\b(t,x,y,W)$ is given~by%
\begin{equation}\label{chi diff}
u^1=W^1,\qquad u^2=y^1+x^1x^2W^1.
\end{equation}
Hence the dynamic feedback linearization of \eqref{tVTOL}, written symbolically in the form \[\dot{x}=f\bigl(x,u^1,u^2\bigr),\] is given by
\begin{align}
&\dot{x}=f\bigl(\bsy{x}, W^1, y^1+x^1x^2W^1 \bigr)\qquad \dot{y}^i=y^{i+1},\quad i=1,\dots, 3,\qquad\dot{y}^4=W^2.\label{DFL4Example3}
\end{align}
It can be checked directly that this control system is SFL. Its trajectories uniquely determine those of \eqref{tVTOL}.
If $s'$ is a solution of \smash{$\widecheck{\CV}$} then $\wt{\pi}\circ\chi\circ s'\colon {\rm I}\subseteq\B R\to M$ solves
$\CV$. The distribution \smash{$\widecheck{\CV}$} is defined by \eqref{DFL4Example3}.

It remains only to identify the fundamental functions which must be Lie differentiated by a vector field
\smash{$Z\in\widecheck{\CV}$} satisfying $\LieD_Z(t)=1$.
By procedure {\it contact} (see Section~\ref{procContact}), we deduce that there is one fundamental function~$\phi^4$ of order~4 and one function $\phi^6$ of order~6. Indeed,
\[
\phi^6=x^5,\qquad \phi^4=\frac{\bigl(\bigl(2x^4 x^5 + 2x^2\bigr) x^6 + x^1 y^1\bigl) x^3 - y^1 x^5 + 2\bigl(x^6\bigr)^2)}{2x^6 x^3 x^5}.
\]
If $(t, \a_a,\ \b_b)$, $0\leq a\leq 4$, $0\leq b\leq 6$ are the contact coordinates on
$J^{\langle 0, 0, 0, 1, 0, 1\rangle}$, then the static feedback transformation
$\psi\colon P\times\pi^{-1}({\rm U})\to J^{\langle 0, 0, 0, 1, 0, 1\rangle}$ that identifies \smash{$\widecheck{\CV}$} with its Brunovsk\'y normal form $\CB_{\langle 0, 0, 0, 1, 0, 1\rangle}$ is given by
\[
\a_a(t, x, y,W)=(\LieD_Z)^a\phi^4,\quad 0\leq a\leq 4,\qquad
\b_b(t, x, y, W)=(\LieD_Z)^b\phi^6,\quad0\leq b\leq 6,
\]
from which the explicit solution is readily deduced, but too lengthy to record here.

Note that the functions
$\phi^4$ and $\phi^6$ are {\it flat outputs} for \eqref{tVTOL} and that they arise canonically as the {\it fundamental functions} of the dynamic linearization \eqref{DFL4Example3}, once $y^1$ is replaced by $u^2-x^1x^2u^1$, according to \eqref{chi diff}. Thus cascade feedback linearization provides a geometric method based on symmetries for constructing flat outputs in the case of invariant flat control systems; see also Section \ref{PVTOLsection}.

\section{Yet another look at the PVTOL}\label{PVTOLsection}

In \cite{DynamicVTOL}, the authors made a study of the well-known PVTOL control system \eqref{tVTOL1} by writing down a complete set of flat outputs and thereby derived a dynamic linearization feedback equivalent to the Brunovsk\'y form
of signature $\langle 0, 0, 0, 2\rangle$. According to their account the authors achieved this by inspired guesswork. In this section, we show how these flat outputs can be derived systematically using only the symmetry considerations of the present paper. We remark that additional recent work \cite{Gstottner_Kolar_Schoberl_2022_Flat,Gstottner_Kolar_Schoberl_2023_Flat} also give a method to systematically develop dynamic feedbacks of the PVTOL system.

We denote the ambient manifold of the control system \eqref{tVTOL1} by $M$
and calculate that the abelian Lie algebra spanned (over $\B R$) by
\begin{equation}\label{pvtolSymmGens}
\bsy{\G}=\operatorname{span}\{X_1, X_2\},
\end{equation}
where
\begin{gather*}
X_1 = h\sin^2(\theta)\cos(\theta)\partial_x+ h\theta_1\bigl(3\cos^2(\theta)-1\bigr)\sin(\theta)\partial_{x_1}+\bigl(x-h\sin(\theta)\cos^2(\theta)\bigr)\partial_z\\
\hphantom{X_1 =}{} +\bigl(x_1+2h\theta_1\cos(\theta)-3h\theta_1\cos^3(\theta)\bigr)\partial_{z_1}+\sin^2(\theta)\partial_\theta+\theta_1\sin(2\theta)\partial_{\theta_1}\\
\hphantom{X_1 =}{} +\cos(\theta)\sin(\theta)\bigl(5h\theta_1^2-u_1\bigr)\partial_{u_1}+\bigl(2\theta_1^2\cos(2\theta)+u_2\sin(2\theta)\bigr)\partial_{u_2},\\
X_2= \partial_z
\end{gather*}
generates a control admissible symmetry group acting on an open subset of $M$. The action of the control admissible symmetry group $G$ generated by $\bsy{\G}$ can be found in the appendix. For the present purposes, we will not record further details here about this subgroup of the maximal Lie group of control symmetries of \eqref{tVTOL1} and instead focus on its application to the dynamic feedback linearization issue at hand.

Letting $\CV$ denote the distribution for system \eqref{tVTOL1}, we have
\begin{align*}
\CV=\operatorname{span}\bigl\{&\P t+x_1\P x-\bigl(u^1\sin\th-h u^2\cos\th\bigr)\P {x_1}+z_1\P z+\bigl(u^1\cos\th+h u^2\sin\th-1\bigr)\P {z_1}\cr
&+\th_1\P {\th}+u^2\P {\th_1},\, \P {u^1}, \,\P {u^2}\bigr\}.
\end{align*}
We find that $\wh{\CV}:=\CV\oplus\bsy{\G}$ is a \stf\ relative Goursat bundle of signature
$\langle 0,2\rangle$, and hence $\CV$ has a contact sub-connection $\wt{\varphi}_*\CV=\CH_G$, which we calculate to be of the form\footnote{See the appendix for the local trivialization $\wt{\varphi}$, in this case. Here $t$, $k$, $k_1$, $k_2$, $g$, $g_1$, $g_2$ denote the standard contact coordinates on $J^{\langle 0,2\rangle}$.}
\begin{align*}
\begin{aligned}
\CH_G= \operatorname{span}\biggl\{&\P t+k_1\P k+k_2\P {k_1}+g_1\P g+g_2\P {g_1}-
\frac{1+k_1}{g_1}\P {\e^1}\\
& {}+\frac{g_1(k-g_1)+g(1+k_1)}{g_1}\P {\e^2}, \,\P {k_2},\, \P {g_2}\biggr\}.
\end{aligned}
\end{align*}
To determine the existence of a dynamic feedback linearization of the system \eqref{tVTOL1}, we carry out a partial contact curve reduction (Section \ref{contactReduction}, Definition \ref{contactReductionDef}). We can do this along either of the variables $k$ or $g$. It turns out that we can choose $g$ and form the map
\[
\mathbf{c}^{\langle 0, 1\rangle}_{f}=
j^{\langle 0, 1\rangle} f\times \text{Id}_{J^{\langle 0,1\rangle}\times G}\colon \
J^{\langle 0,1\rangle}\times G\to J^{\langle 0, 2\rangle}\times G,
\]
where $j^{\langle 0, 1\rangle} f=\bigl(t, f(t), \dot{f}(t), \ddot{f}(t)\bigr)=(t, g, g_1, g_2)$ and $f$ is an arbitrary, smooth, real-valued function of $t$. As in Section \ref{contactReduction}, we form the Pfaffian system \smash{$\bigl(\mathbf{c}^{\langle 0, 1\rangle}_{f}\bigr)^*\bsy{\g}^G=:\bsy{\bar{\g}}^G$}.
Then
\[
\bar{\CH}_G=\operatorname{span}\left\{\P t+k_1\P k+k_2\P {k_1}-\frac{1+k_1}{\dot{f}}\P {\e^1}+
\frac{\dot{f}(k-\dot{f})+f(1+k_1)}{\dot{f}}\P {\e^2},\,
 \P {k_2}\right\},
\]
where $\bar{\CH}_G$ is the kernel of $\bsy{\bar{\g}}^G$ on $J^{\langle 0,1\rangle}\times G$. By Theorems \ref{Generalized Goursat Normal Form} and \ref{Goursat SFL}, we find that
$\bar{\CH}_G$ is static feedback linearizable with signature $\langle 0, 0, 0, 1\rangle$. Thus, the derived length of $\bar{\CH}_G$ is $\bar{k}=4$ and by Corollary~\ref{prolongationPredictor-cor}, the 7-fold partial prolongation of $\CH_G$ along the $g$-series of contact coordinates is sufficient to achieve a static feedback linearizable control system. However, this over-estimates the number of prolongations and in fact, it is easy to discover that only 2 prolongations by differentiation of $g_2$ are required; this can be easily deduced by computing the fundamental function of $\bar{\CH}_G$. With only 2 prolongations, we get
\begin{align*}
\operatorname{pr} \CH_G= \operatorname{span}\biggl\{&\P t+k_1\P k+k_2\P {k_1}+g_1\P g+g_2\P {g_1}+g_3\P {g_2}+
g_4\P {g_3}-\frac{1+k_1}{g_1}\P {\e^1} \\
& +\frac{g_1(k-g_1)+g(1+k_1)}{g_1}\P {\e^2},\, \P {k_2}, \,\P {g_4}\biggr\}
\end{align*}
on $J^{\langle 0, 1, 0, 1\rangle}\times G$ and a further calculation shows that $\operatorname{pr}\CH_G\simeq\CB_{\langle 0, 0, 0, 2\rangle}$. We have therefore proven that \eqref{tVTOL1} is \cfl\ and hence \dfl\ with a dynamic linearization of signature
$\langle 0, 0, 0, 2\rangle$. Our next goal is to derive the linearizable dynamic extension of~\eqref{tVTOL1}.

The map \smash{$\wt{\wt{\varphi}}$} is given by \smash{$\wt{\wt{\varphi}}=\wt{\varphi}\circ\wt{\pi}\times\{g_3=p_1, \ g_4=p_2\}$}. We do not record the transformation~\smash{$\wt{\wt{\varphi}}^{-1}$} here; however, the relevant part of \smash{$\wt{\wt{\vf}}$} needed to construct the dynamic feedback is
\begin{gather}\label{dcomp pvtol}
g_2=\sin(\th)\bigl(h\th_1^2-u^1\bigr).
\end{gather}
The distribution $\mcal{T}:=\wt{\wt{\varphi}}^{-1}_*\operatorname{pr}\CH_G$ has the form
\begin{align*}
\mcal{T}= \operatorname{span}\bigl\{&\P t+x_1\P x-\bigl(u^1\sin\th-h u^2\cos\th\bigr)\P {x_1}+z_1\P z\\
&+\bigl(u^1\cos\th+h u^2\sin\th-1\bigr)\P {z_1}+\th_1\P {\th}+u^2\P {\th_1}\\
&+\bigl(\th_1\bigl(\th_1^2h-u^1\bigr)\cot\th+2hu^2\th_1 -p_1\bigr)\P {u^1}+p_2\P {p_1},\, \P {u^2},\, \P {p_2}\bigr\}.
\end{align*}
We notice that
\[
\ch\mcal{T}^{(1)}_0=\operatorname{span}\{\partial_{u^2},\partial_{p_2}\},
\]
which ensures that $\mcal{T}$ is a dynamic extension of the original PVTOL control system $\mcal{V}$. Notice, however, that we can still use equation \eqref{dcomp pvtol} to construct a dynamic extension as described in Theorem~\ref{cfl to dfl}. Indeed, let $\chi^{-1}$ be the diffeomorphism such that
\[
\bigl(y^1,y^2\bigr)=\bigl(\sin(\th)\bigl(h\th_1^2-u^1\bigr),p_1\bigr),\qquad
\bigl(W^1,W^2\bigr)=\bigl(u^2,p_2\bigr),
\]
together with the identity on all other variables. Then the dynamic compensator $\bsy{u}=\bsy{\beta}(t,\bsy{x},\allowbreak\bsy{y},\bsy{W})$ is $u^1=h\theta_1^2-\csc(\th)y^1$, $u^2=W^1$.
As such, the control system \eqref{tVTOL1} is dynamically extended to
\begin{align}
&\dot{x}=x_1\qquad
\dot{x}_1=-h\theta_1^2\sin\th-y^1+h W^1\cos\th,\qquad
\dot{z}=z_1,\nonumber\\
&\dot{z}_1=h\theta_1^2\cos(\th)-\cot(\th)y^1+h W^1\sin\th-1,\qquad
\dot{\th}=\th_1,\qquad
\dot{\th}_1=W^1,\nonumber\\
&\dot{y^1}=y^2,\qquad
\dot{y^2}=W^2,\label{dfl pvtol}
\end{align}
which has linear dynamics along the fibers of the submersion, with the corresponding distribution denoted by
\smash{$\widecheck{\CV}$}.

To compute the flat outputs of \eqref{tVTOL1}, we compute the fundamental functions of either $\mcal{T}$ or~\smash{$\widecheck{\CV}$} in the usual way. These turn out to be $
x-h\sin\th$\ and $z+h\cos\th$,
the so called {\it center of oscillation} flat outputs of \eqref{tVTOL1} first deduced by Martin, Devasia and Paden in \cite{DynamicVTOL}
but here obtained canonically from the fundamental functions of the dynamic linearization \eqref{dfl pvtol}. The authors in \cite{DynamicVTOL} rather proceed in the opposite direction: the flat outputs are firstly ingeniously guessed and then verified by deriving the dynamic linearization \eqref{dfl pvtol} from them.

\section{Conclusion}
We have studied generic, smooth control systems $(M,\bsy{\o})$ that are invariant under a Lie group~$G$ of control admissible transformations satisfying the further property that they admit \sfl\ quotients by $G$ (symmetry reductions). An {\it infinitesimal test} for rapidly identifying \sfl\ quotients was given, based on the Lie algebra of infinitesimal generators of the action of $G$.
The requirement that an invariant control system possesses a \sfl\ quotient appears to be very mild. In practice, it appears to hold quite generally, in contrast to the \stf\ linearizability of a control system itself, which is rare.

We proved that a $C^\infty$ control system invariant under a Lie group $G$ of control admissible transformations and that has a \sfl\ quotient, admits the local normal form
$\bsy{\g}^G=\bsy{\b}^\k\oplus\Theta$
on a local trivialization $\pi'\colon J^\k\times G\to J^\k$ of the principal bundle~${\pi\colon (M, \bsy{\o})\to (M/G, \bsy{\o}/G)}$ in which $(J^\k, \bsy{\b}^\k)$ is the Brunovsk\'y normal form of $(M/G, \bsy{\o}/G)$ and $\Theta$ is a differential system for the reconstruction of the trajectories of $\bsy{\o}$ from those of the quotient $\bsy{\o}/G$.
The present paper built on this result to construct \df\ linearizations of control systems with symmetry.

A widely applicable procedure was established showing that any cascade feedback linearizable control system can be dynamically linearized with linear dynamics in the fiber coordinates. It was proven that this can be achieved by a prolongation by differentiation of {\em canonically chosen} inputs of the associated contact sub-connection $\CH_G$. Theorem \ref{prolongationPredictor} and Corollary \ref{prolongationPredictor-cor} established a bound on the number of differentiations that must be performed to achieve a \stf\ linearizable dynamically extended control system that appears to outperform the best currently known bounds \cite{ImprovedST-2006}, being linear in the derived length of the reduced sub-connection $\bar{\CH}_G$. For instance, the sharp bound of Sluis and Tilbury \cite{SluisTilbury} in the case of 2-input control systems is achieved by our general bound. However, we have not, in this paper, made a rigorous general comparison of efficiency. On the other hand, a procedure was given whereby the {\it exact} number of differentiations can be predicted.

We have given numerous examples to show how the general theory can be applied in practice. In particular,
a symmetry-based construction of the {centre of oscillation flat outputs} of the PVTOL control system was described. These flat outputs together with the corresponding dynamic linearization were first presented in \cite{DynamicVTOL}. It was also shown that the PVTOL also has a dynamic linearization by differentiation of a {\em given} input by studying the control system's invariance under a subgroup of the Galilean transformations in $x$ and $z$. These are instances of a general symmetry-based procedure for deriving the flat outputs of flat control systems with symmetry.\looseness=-1

Throughout we have emphasized intrinsic geometric structures underlying control systems with symmetry to derive coordinate-free, canonical procedures that inform the local coordinate calculations that lead to their dynamic feedback linearizations and explicit solutions.

\appendix

\section{PVTOL group action}\label{PVTOL-group-action}

The local Lie group action generated by the 2-dimensional control admissible symmetry algebra~\eqref{pvtolSymmGens} of Section \ref{PVTOLsection} is given by
\begin{align*}
&t\mapsto t,\qquad
\theta\mapsto F(\theta,\varepsilon_1),\qquad
\theta_1\mapsto \theta_1\csc^2(\theta)\sin^2(F(\theta,\varepsilon_1)),\\
&x\mapsto x+h(\sin(F(\theta,\varepsilon_1))-\sin(\theta)),\\
&x_1\mapsto x_1-h\theta_1\cos(\theta)+h\theta_1\csc^2(\theta)\cos(F(\theta,\varepsilon_1))\sin^2(F(\theta,\varepsilon_1)),\\
&z\mapsto \varepsilon_1x+z+\varepsilon_2+h\cos(\theta)-h\varepsilon_1\sin(\theta)-h\cos(F(\theta,\varepsilon_1)),\\
&z_1\mapsto h\theta_1\csc^2(\theta)\sin^3(F(\theta,\varepsilon_1))+(x_1-h\theta_1\cos(\theta))\varepsilon_1+z_1-h\theta_1\sin(\theta),\\
&u_1\mapsto \bigl(u_1\sin(\theta)-h\theta_1^2\sin(\theta)\bigr)\csc(F(\theta,\varepsilon_1))+h\theta_1^2\csc^4(\theta)\sin^4(F(\theta,\varepsilon_1)),\\
&u_2\mapsto \bigl(u_2\csc(\theta)-2\theta_1^2\cos(\theta)\csc^3(\theta)\bigr)\sin^2(F(\theta,\varepsilon_1))\\
&\phantom{u_2\mapsto}+2\theta_1^2\csc^4(\theta)\cos(F(\theta,\varepsilon_1))\sin^3(F(\theta,\varepsilon_1)),
\end{align*}
where $F(\theta,\varepsilon_1)=\arccot (\cot(\theta)-\varepsilon_1)$, and $\e_1$ and $\e_2$ form the local coordinates of the control admissible symmetry group $G$. With $t$, $k$, $k_1$, $k_2$, $g$, $g_1$, $g_2$ denoting the standard contact coordinates on
$J^{\langle 0,2\rangle}$, the local trivialization map $\wt{\vf}\colon M\to J^{\langle 0,2\rangle}\times G$ is given by
\begin{gather*}
t =t,\qquad
g =x-h\sin(\theta),\qquad
g_1 =x_1-h\theta_1\cos(\theta),\qquad
g_2 =\bigl(h\theta_1^2-u_1\bigr)\sin(\theta),\\
\e_1 =1-\cot(\theta),\qquad
k =(x_1\cos(\theta)-h\theta_1)\csc(\theta)+z_1,\\
k_1 =\bigl(h\theta_1^2\cos(\theta)-\theta_1x_1\bigr)\csc^2(\theta)-1,\\
k_2 =-\csc^3(\theta)\theta_1^2\bigl(\theta_1h-x_1\cos(\theta)\bigr)+u_1\theta_1\csc(\theta)\\
\phantom{k_2 =}{} -\csc^3(\theta)\bigl(\cos(\theta)\theta_1h-x_1\bigr)\bigl(\theta_1^2\cos(\theta)-u_2\sin(\theta)\bigr), \\
\e_2 =\left(z-x-\frac{h}{\sqrt{2}}+h\sin(\theta)\right)\sin(\th)+x\cos(\theta).
\end{gather*}

\subsection*{Acknowledgements}

We are grateful to the Simons Foundation for its support of the first author via a Collaboration Grant for Mathematicians. We would also like to thank the anonymous referees for their careful reviews and helpful suggestions; this paper is much improved thanks to their efforts.

\pdfbookmark[1]{References}{ref}
\LastPageEnding


\begin{thebibliography}{99}
\footnotesize\itemsep=0pt

\bibitem{AndersonTorre}
Anderson I., Torre C., The DifferentialGeometry package, available at
 \url{http://digitalcommons.usu.edu.au/dg}.

\bibitem{AF}
Anderson I.M., Fels M.E., Exterior differential systems with symmetry,
 \href{https://doi.org/10.1007/s10440-005-1136-y}{\textit{Acta Appl. Math.}} \textbf{87} (2005), 3--31.

\bibitem{ArandaMoogPomet}
Aranda-Bricaire E., Moog C.H., Pomet J.-B., A linear algebraic framework for
 dynamic feedback linearization, \href{https://doi.org/10.1109/9.362886}{\textit{IEEE Trans. Automat. Control}}
 \textbf{40} (1995), 127--132.

\bibitem{AP}
Aranda-Bricaire E., Pomet J.-B., Some explicit conditions for a control system
 to be feedback equivalent to extended {G}oursat form, \href{https://doi.org/10.1016/S1474-6670(17)46877-7}{\textit{IFAC Proc.
 Vol.}} \textbf{28} (1995), 489--494.

\bibitem{AP07}
Avanessoff D., Pomet J.-B., Flatness and {M}onge parameterization of two-input
 systems, control-affine with~4 states or general with 3 states, \href{https://doi.org/10.1051/cocv:2007011}{\textit{ESAIM
 Control Optim. Calc. Var.}} \textbf{13} (2007), 237--264,
 \href{https://arxiv.org/abs/math.OC/0505443}{arXiv:math.OC/0505443}.

\bibitem{BC3control}
Battilotti S., Califano C., A geometric approach to dynamic feedback
 linearization, in Analysis and {D}esign of {N}onlinear {C}ontrol {S}ystems,
 \href{https://doi.org/10.1007/978-3-540-74358-3_23}{Springer}, Berlin, 2008, 397--411.

\bibitem{Brunovsky}
Brunovsk\'y P., A classification of linear controllable systems,
 \textit{Kybernetika (Prague)} \textbf{6} (1970), 173--188.

\bibitem{BryantThesis}
Bryant R.L., Some aspects of the local and global theory of {P}faffian systems,
 Ph.D.~Thesis, {T}he University of North Carolina at Chapel Hill, 1979.

\bibitem{CharletLevineMarino2}
Charlet B., L\'evine J., Marino R., Sufficient conditions for dynamic state
 feedback linearization, \href{https://doi.org/10.1137/0329002}{\textit{SIAM~J. Control Optim.}} \textbf{29} (1991),
 38--57.

\bibitem{Chetverikov}
Chetverikov V.N., New flatness conditions for control systems, \href{https://doi.org/10.1016/S1474-6670(17)35172-8}{\textit{IFAC
 Proc. Vol.}} \textbf{34} (2001), 191--196.

\bibitem{ChetverikovDFL}
Chetverikov V.N., Flatness of dynamically linearized systems, \href{https://doi.org/10.1007/s10625-005-0106-5}{\textit{Differ.
 Equ.}} \textbf{40} (2004), 1747--1756.

\bibitem{ClellandHu}
Clelland J.N., Hu Y., On absolute equivalence and linearization~{I},
 \href{https://doi.org/10.1007/s10711-023-00811-0}{\textit{Geom. Dedicata}} \textbf{217} (2023), 78, 31~pages,
 \href{https://arxiv.org/abs/2005.00643}{arXiv:2005.00643}.

\bibitem{DTVa}
De~Don\'a J., Tehseen N., Vassiliou P.J., Symmetry reduction, contact geometry,
 and partial feedback linearization, \href{https://doi.org/10.1137/15M1046538}{\textit{SIAM~J. Control Optim.}}
 \textbf{56} (2018), 201--230, \href{https://arxiv.org/abs/1510.05761v1}{arXiv:1510.05761v1}.

\bibitem{TVFL1}
D'Souza R.S., Nielsen C., An algorithm for local transverse feedback
 linearization, \href{https://doi.org/10.1137/21M1444588}{\textit{SIAM~J. Control Optim.}} \textbf{61} (2023),
 1248--1272, \href{https://arxiv.org/abs/2109.00641}{arXiv:2109.00641}.

\bibitem{Elkin}
Elkin V.I., Reduction of nonlinear control systems, \textit{Math. Appl.}, Vol.
 472, \href{https://doi.org/10.1007/978-94-011-4617-3}{Kluwer Academic Publishers}, Dordrecht, 1999.

\bibitem{FLRM1}
Fliess M., L\'evine J., Martin P., Rouchon P., Flatness and defect of
 non-linear systems: introductory theory and examples, \href{https://doi.org/10.1080/00207179508921959}{\textit{Internat.~J.
 Control}} \textbf{61} (1995), 1327--1361.

\bibitem{FLRM2}
Fliess M., L\'evine J., Martin P., Rouchon P., A {L}ie--{B}\"acklund approach
 to equivalence and flatness of nonlinear systems, \href{https://doi.org/10.1109/9.763209}{\textit{IEEE Trans.
 Automat. Control}} \textbf{44} (1999), 922--937.

\bibitem{ImprovedST-2006}
Franch J., Fossas E., Linearization by prolongations: new bounds on the number
 of integrators, \href{https://doi.org/10.3166/ejc.11.171-179}{\textit{Eur.~J. Control}} \textbf{11} (2005), 171--179.

\bibitem{GSalgorithmExample}
Gardner R.B., Shadwick W.F., An algorithm for feedback linearization,
 \href{https://doi.org/10.1016/0926-2245(91)90028-8}{\textit{Differential Geom. Appl.}} \textbf{1} (1991), 153--158.

\bibitem{GSalgorithm}
Gardner R.B., Shadwick W.F., The {GS} algorithm for exact linearization to
 {B}runovsky normal form, \href{https://doi.org/10.1109/9.121623}{\textit{IEEE Trans. Automat. Control}} \textbf{37}
 (1992), 224--230.

\bibitem{Grizzle}
Grizzle J.W., Marcus S.I., The structure of nonlinear control systems
 possessing symmetries, \href{https://doi.org/10.1109/TAC.1985.1103927}{\textit{IEEE Trans. Automat. Control}} \textbf{30}
 (1985), 248--258.

\bibitem{Gstottner_Kolar_Schoberl_2022_Flat}
Gst\"ottner C., Kolar B., Sch\"oberl M., A structurally flat triangular form
 based on the extended chained form, \href{https://doi.org/10.1080/00207179.2020.1841302}{\textit{Internat.~J. Control}} \textbf{95}
 (2022), 1144--1163, \href{https://arxiv.org/abs/2007.09935}{arXiv:2007.09935}.

\bibitem{Gstottner_Kolar_Schoberl_2023_Flat}
Gst\"ottner C., Kolar B., Sch\"oberl M., Necessary and sufficient conditions
 for the linearisability of two-input systems by a two-dimensional endogenous
 dynamic feedback, \href{https://doi.org/10.1080/00207179.2021.2015542}{\textit{Internat.~J. Control}} \textbf{96} (2023), 800--821,
 \href{https://arxiv.org/abs/2106.14722}{arXiv:2106.14722}.

\bibitem{Guay}
Guay M., McLellan P.J., Bacon D.W., A condition for dynamic feedback
 linearization of control-affine nonlinear systems, \href{https://doi.org/10.1080/002071797223749}{\textit{Internat.~J.
 Control}} \textbf{68} (1997), 87--106.

\bibitem{VTOL}
Hauser J., Sastry S., Meyer G., Nonlinear control design for slightly
 nonminimum phase systems: application to {V}/{STOL} aircraft,
 \href{https://doi.org/10.1016/0005-1098(92)90029-F}{\textit{Automatica~J. IFAC}} \textbf{28} (1992), 665--679.

\bibitem{HuntSuMeyerLin}
Hunt L.R., Su R.J., Meyer G., Design for multi-input nonlinear systems, in
 Differential {G}eometric {C}ontrol {T}heory ({H}oughton, {M}ich., 1982),
 \textit{Progr. Math.}, Vol.~27, Birkh\"auser, Boston, MA, 1983, 268--298.

\bibitem{RespondekLin}
Jakubczyk B., Respondek W., On linearization of control systems, \textit{Bull.
 Acad. Polon. Sci. S\'er. Sci. Math.} \textbf{28} (1980), 517--522.

\bibitem{KlotzThesis}
Klotz T.J., Geometry of cascade feedback linearizable control systems,
 \href{https://doi.org/10.1016/j.difgeo.2023.102044}{\textit{Differential Geom. Appl.}} \textbf{90} (2023), 102044, 37~pages,
 \href{https://arxiv.org/abs/2102.08521}{arXiv:2102.08521}.

\bibitem{KobNomizu}
Kobayashi S., Nomizu K., Foundations of differential geometry. {V}ol~{I},
 Interscience Publishers, New York, 1963.

\bibitem{LevineDFL}
L\'evine J., On the equivalence between differential flatness and dynamic
 feedback linearizability, \href{https://doi.org/10.3182/20071017-3-BR-2923.00056}{\textit{IFAC Proc. Vol.}} \textbf{40} (2007),
 338--343.

\bibitem{Levine}
L\'evine J., Analysis and control of nonlinear systems: {A} flatness-based
 approach, \textit{Math. Eng.}, \href{https://doi.org/10.1007/978-3-642-00839-9}{Springer}, Berlin, 2009.

\bibitem{LNR}
Li S., Nicolau F., Respondek W., Multi-input control-affine systems static
 feedback equivalent to a triangular form and their flatness,
 \href{https://doi.org/10.1080/00207179.2015.1056232}{\textit{Internat.~J. Control}} \textbf{89} (2016), 1--24, \href{https://arxiv.org/abs/1411.6282}{arXiv:1411.6282}.

\bibitem{Marino}
Marino R., On the largest feedback linearizable subsystem, \href{https://doi.org/10.1016/0167-6911(86)90130-1}{\textit{Systems
 Control Lett.}} \textbf{6} (1986), 345--351.

\bibitem{MarinoBoothby85}
Marino R., Boothby W.M., Elliott D.L., Geometric properties of linearizable
 control systems, \href{https://doi.org/10.1007/BF01699463}{\textit{Math. Systems Theory}} \textbf{18} (1985), 97--123.

\bibitem{DynamicVTOL}
Martin P., Devasia S., Paden B., A different look at output tracking: control
 of a {VTOL} aircraft, \href{https://doi.org/10.1016/0005-1098(95)00099-2}{\textit{Automatica~J. IFAC}} \textbf{32} (1996),
 101--107.

\bibitem{MontgomeryGoursatFlags2001}
Montgomery R., Zhitomirskii M., Geometric approach to {G}oursat flags,
 \href{https://doi.org/10.1016/S0294-1449(01)00076-2}{\textit{Ann. Inst. H.~Poincar\'e C Anal. Non Lin\'eaire}} \textbf{18} (2001),
 459--493.

\bibitem{Nic_Res_2016_Flatness}
Nicolau F., Respondek W., Flatness of two-input control-affine systems
 linearizable via a two-fold prolongation, in 2016 {IEEE} 55th {C}onference on
 {D}ecision and {C}ontrol {(CDC)}, \href{https://doi.org/10.1109/CDC.2016.7798852}{IEEE}, 2016, 3862--3867.

\bibitem{Nic_Res_2017_Flatness}
Nicolau F., Respondek W., Flatness of multi-input control-affine systems
 linearizable via one-fold prolongation, \href{https://doi.org/10.1137/140999463}{\textit{SIAM~J. Control Optim.}}
 \textbf{55} (2017), 3171--3203.

\bibitem{OlverLieBook}
Olver P.J., Applications of {L}ie groups to differential equations, 2nd ed.,
 \textit{Grad. Texts Math.}, Vol. 107, \href{https://doi.org/10.1007/978-1-4612-4350-2}{Springer}, New York, 1993.

\bibitem{OlverSymmetryBook}
Olver P.J., Equivalence, invariants, and symmetry, \href{https://doi.org/10.1017/CBO9780511609565}{Cambridge University Press},
 Cambridge, 1995.

\bibitem{palais}
Palais R.S., A global formulation of the {L}ie theory of transformation groups,
 \href{https://doi.org/10.1090/memo/0022}{\textit{Mem. Amer. Math. Soc.}} \textbf{22} (1957), iii+123~pages.

\bibitem{PdS08}
Pereira~da Silva P.S., Some remarks on static-feedback linearization for
 time-varying systems, \href{https://doi.org/10.1016/j.automatica.2008.10.001}{\textit{Automatica~J. IFAC}} \textbf{44} (2008),
 3219--3221.

\bibitem{Pomet95}
Pomet J.-B., A differential geometric setting for dynamic equivalence and
 dynamic linearization, in Geometry in {N}onlinear {C}ontrol and
 {D}ifferential {I}nclusions ({W}arsaw, 1993), \textit{Banach Center Publ.},
 Vol.~32, \href{https://doi.org/10.4064/-32-1-319-339}{Polish Acad. Sci. Inst. Math.}, Warsaw, 1995, 319--339.

\bibitem{Pomet97}
Pomet J.-B., On dynamic feedback linearization of four-dimensional affine
 control systems with two inputs, \href{https://doi.org/10.1051/cocv:1997107}{\textit{ESAIM Control Optim. Calc. Var.}}
 \textbf{2} (1997), 151--230.

\bibitem{R2}
Respondek W., Symmetries and minimal flat outputs of nonlinear control systems,
 in New {T}rends in {N}onlinear {D}ynamics and {C}ontrol, and their
 {A}pplications, \textit{Lect. Notes Control Inf. Sci.}, Vol. 295, \href{https://doi.org/10.1007/978-3-540-45056-6_5}{Springer},
 Berlin, 2003, 65--86.

\bibitem{Shadwick90}
Shadwick W.F., Absolute equivalence and dynamic feedback linearization,
 \href{https://doi.org/10.1016/0167-6911(90)90041-R}{\textit{Systems Control Lett.}} \textbf{15} (1990), 35--39.

\bibitem{SPR}
Silveira H.B., Pereira~da Silva P.S., Rouchon P., A flat triangular form for
 nonlinear systems with two inputs: necessary and sufficient conditions,
 \href{https://doi.org/10.1016/j.ejcon.2015.01.001}{\textit{Eur.~J. Control}} \textbf{22} (2015), 17--22, \href{https://arxiv.org/abs/1312.3527}{arXiv:1312.3527}.

\bibitem{SluisTilbury}
Sluis W.M., Tilbury D.M., A bound on the number of integrators needed to
 linearize a control system, in Proceedings of 1995 34th {IEEE} {C}onference
 on {D}ecision and {C}ontrol, \href{https://doi.org/10.1109/CDC.1995.478961}{IEEE}, 1995, 602--607.

\bibitem{vdS84}
van~der Schaft A.J., Linearization and input-output decoupling for general
 nonlinear systems, \href{https://doi.org/10.1016/0167-6911(84)90005-7}{\textit{Systems Control Lett.}} \textbf{5} (1984), 27--33.

\bibitem{NMR}
van Nieuwstadt M., Rathinam M., Murray R.M., Differential flatness and absolute
 equivalence of nonlinear control systems, \href{https://doi.org/10.1137/S0363012995274027}{\textit{SIAM~J. Control Optim.}}
 \textbf{36} (1998), 1225--1239.

\bibitem{VassiliouGoursat}
Vassiliou P.J., A constructive generalised {G}oursat normal form,
 \href{https://doi.org/10.1016/j.difgeo.2005.12.001}{\textit{Differential Geom. Appl.}} \textbf{24} (2006), 332--350,
 \href{https://arxiv.org/abs/math.DG/0404377}{arXiv:math.DG/0404377}.

\bibitem{VassiliouGoursatEfficient}
Vassiliou P.J., Efficient construction of contact coordinates for partial
 prolongations, \href{https://doi.org/10.1007/s10208-004-0148-8}{\textit{Found. Comput. Math.}} \textbf{6} (2006), 269--308,
 \href{https://arxiv.org/abs/math.DG/0406234}{arXiv:math.DG/0406234}.

\bibitem{VassiliouCascade1}
Vassiliou P.J., Cascade linearization of invariant control systems,
 \href{https://doi.org/10.1007/s10883-017-9389-0}{\textit{J.~Dyn. Control Syst.}} \textbf{24} (2018), 593--623.

\end{thebibliography}
\end{document}